\documentclass[authoryear]{elsarticle}

\usepackage{multirow}
\usepackage[ruled]{algorithm2e}
\usepackage{graphicx}  
\usepackage[font=small,format=plain,labelfont=bf,up,textfont=it,up]{caption}

\usepackage{todonotes}
\usepackage{subfigure}
\usepackage{geometry}
\usepackage{amsmath}
\usepackage{amssymb}
\usepackage{amsfonts}
\usepackage{amsthm}

\usepackage{fancyhdr}
\definecolor{berkeley_lab_blue}{RGB}{0,57,90}

\usepackage{float}

\usepackage{pgf,tikz}
\usetikzlibrary{shapes,arrows,calc,automata}

\usepackage{ctable}

\let\oldref\ref

\renewcommand{\ref}[1]{(\oldref{#1})}

\title{\textbf{Parallel-in-Time Multi-Level Integration \\ of the Shallow-Water Equations  on the Rotating Sphere}}
\date{}

\author[lbl]{Fran\c cois P. Hamon\corref{cor1}}
\ead{fhamon@lbl.gov}
\author[uexeter]{Martin Schreiber}

\author[lbl2]{Michael L. Minion}

\cortext[cor1]{Corresponding author, now with TOTAL Exploration \& Production}
\address[lbl]{Center for Computational Sciences and Engineering, Lawrence Berkeley National Laboratory, Berkeley, USA}
\address[uexeter]{Chair of Computer Architecture and Parallel Systems, Technical University of Munich, Germany}
\address[lbl2]{Department of Applied Mathematics, Lawrence Berkeley National Laboratory, Berkeley, USA}

\begin{document}

\begin{abstract} 
  The modeling of atmospheric processes in the context of
  weather and climate simulations is an important and
  computationally expensive challenge. The temporal integration
  of the underlying PDEs requires a very large number of time
  steps, even when the terms accounting for the propagation
  of fast atmospheric waves are treated implicitly. Therefore,
  the use of  parallel-in-time integration schemes to reduce
  the time-to-solution is of increasing interest, particularly
  in the numerical weather forecasting field.

  We present a multi-level parallel-in-time integration
  method combining the Parallel Full Approximation Scheme in
  Space and Time (PFASST) with a spatial discretization
  based on Spherical Harmonics (SH). The iterative algorithm
  computes multiple time steps concurrently by
  interweaving parallel high-order fine corrections and serial
  corrections performed on a coarsened problem. To do that, we design a
  methodology relying on the spectral basis of the SH
  to coarsen and interpolate the problem in space. 

  The methods are evaluated on the shallow-water
  equations on the sphere using a set of tests
  commonly used in the atmospheric flow community.
  We assess the convergence of PFASST-SH upon refinement
  in time. We also investigate the impact of the coarsening
  strategy on the accuracy of the scheme, and specifically
  on its ability to capture the high-frequency modes
  accumulating in the solution. Finally, we study the
  computational cost of PFASST-SH to demonstrate that
  our scheme resolves the main features of the solution
  multiple times faster than the serial schemes. 
\end{abstract}


\maketitle

\section{\label{section_introduction}Introduction}

The accurate simulation of atmospheric flows over long periods
of time is one of the critical components in the fields of
numerical weather prediction and climate modeling.  Despite
a rich history of numerical methods in these fields, the development
of more accurate and efficient temporal integration methods
for weather and climate simulations is still an ongoing
challenge. 
One difficulty is the
presence of a wide range of time scales in the equations,
including the propagation of fast atmospheric waves, which
imposes a severe stability restriction on the time step size
of fully explicit integration schemes. Implicit-explicit
schemes partly overcome this limitation by treating the
stiff terms implicitly, but still involve a very large number
of time steps for high-resolution long-range simulations
(up to a thousand years for paleoclimate studies). Therefore,
parallel-in-time methods are an attractive approach to reduce
the time-to-solution by accessing  an additional  axis of
parallelism in the temporal direction. 

In this work, we study the performance and accuracy of an
implicit-explicit, iterative, multi-level, parallel-in-time
integration scheme based on the Parallel Full Approximation
Scheme in Space and Time \citep[PFASST,][]{emmett2012toward}
in the context of atmospheric flows. The method is applied
to the Shallow-Water Equations (SWE) on the rotating sphere,  a
common two-dimensional proxy that incorporates the horizontal
features of the full three-dimensional atmospheric equations.
In the classification of \cite{burrage1997parallel}, the
parallelization  strategy in PFASST relies on
\textit{parallelization across the steps}, which consists in
solving multiple time steps concurrently on different processing
units to accelerate the simulations. Other examples
of this type of parallelism include Parareal
\citep{lions2001resolution},
the Parallel Implicit Time-integration
Algorithm \citep[PITA,][]{farhat2003time}
and the MultiGrid Reduction in Time
\citep[MGRIT,][]{falgout2014parallel}. This approach differs
from other classes of methods not considered
in this work, namely, \textit{parallelization across the method},
in which intermediate stage values are computed in parallel 
\citep[e.g.,][]{butcher1997order}, and
\textit{parallelization across the problem},
in which the full problem is split into subproblems coupled
iteratively \citep[e.g.,][]{gander1999waveform}. 
Parallel-in-time methods based on exponential integration have
also been proposed
\citep{gander2013paraexp,haut2014asymptotic}.
They have been applied to the SWE on the plane
\citep{schreiber2017beyond} and on the rotating sphere
\citep{schreiber2018sph} using a REXI (Rational Approximation
of Exponential Integrators) approach to compute the exponential. 
Unlike PFASST, these schemes integrate one time step at a time,
but perform the integration of a given time step by summing the
solutions of fully decoupled systems solved in parallel. This
yields a highly accurate 
and stable integration of the linear(ized) terms that requires an
efficient solver for complex-valued linear systems.
For nonlinear problems, exponential integration must be combined with
other time integrators handling the nonlinear terms using a splitting
scheme \citep{schreiber_2019_expnonlinearswe_sphere}.

The structure of PFASST is similar to that of Parareal as both
algorithms combine parallel, relatively expensive updates based
on approximate initial conditions and serial inexpensive updates
propagating the new initial conditions on the time interval.
However, the PFASST solution updates are computed iteratively
with Spectral Deferred Corrections
\citep[SDC,][]{dutt2000spectral}
instead of direct ODE
integrators such as Runge-Kutta schemes in standard Parareal.
This makes it possible to amortize the cost of the fine-level
SDC updates  over multiple PFASST iterations
\citep{minion2011hybrid}, thereby reducing the computational cost
of the fine propagator. The flexibility inherent in the SDC
approach also allows for the straightforward construction of
implicit-explicit (IMEX) methods
of very high order, which has not yet been attained with methods
like linear-multistep or Runge-Kutta methods \citep{minion2003semi}.
IMEX schemes are preferred in the current
setting over explicit or fully implicit methods.

In PFASST the  SDC corrections are interwoven on a hierarchy of
levels representing coarse space-time approximations of the fine
problem under consideration. This iterative procedure allowing 
for spatial coarsening results in a larger theoretical parallel
efficiency than with the Parareal algorithm. An efficient 
implementation of the PFASST algorithm has been proposed
\citep{emmett2014efficient}, and the parallel efficiency of the
scheme has been demonstrated for some applications, including
the N-body problem \citep{speck2012massively} and the heat 
equation \citep{speck2014space}. Recently, in
\cite{bolten2016multigrid,bolten2017asymptotic}, the PFASST
algorithm has been cast as a multigrid method and its convergence
has been studied for synthetic diffusion-dominated and
advection-dominated problems.  PFASST has also been used in the
context of optimal control problems \citep{Gotschel2019-ka}.

The space-time coarsening and interpolation strategy is one of
the critical determinants of the accuracy and efficiency
of the PFASST algorithm. We consider the approach
proposed in the serial Multi-Level Spectral
Deferred Corrections (MLSDC) scheme in \cite{hamon2018multi}.
It takes advantage of the structure of the spatial discretization
based on the Spherical Harmonics (SH) transform, which plays a key
role in our parallel-in-time scheme, thereafter referred
to as PFASST-SH. Specifically, we have shown in previous
work that the SH basis can be truncated, or padded with zeros,
to construct consistent spatial restriction and interpolation
operators between levels. As in \cite{hamon2018multi}, we
exploit the properties of the SH transform to efficiently solve
the implicit systems at each temporal node. This study is
relevant for operational numerical weather prediction systems
as the SH transform is used in the Integrated Forecast System
(IFS) at the European Centre for Medium-Range Weather Forecasts
\citep[ECMWF, ][]{wedi2013fast} and the Global Spectral Model
(GSM) at the Japan Meteorological Agency 
\citep[JMA, ][]{kanamitsu1983description}.

We evaluate the properties of the proposed PFASST-SH algorithm
using a well-known suite of nonlinear shallow-water test
cases that are representative of the horizontal features of
atmospheric flows \citep{williamson1992standard,galewsky2004initial}.
These wave-propagation examples are challenging
for the integration scheme as parallel-in-time methods are
known to suffer from convergence difficulties on advection-dominated
problems \citep[e.g.,][]{ruprecht2018wave}. In addition, the
nonlinear problems considered in this work are characterized by the
progressive amplification over time of high-frequency modes not
captured by the coarse SDC updates. This can undermine the
accuracy of the parallel approach on these modes when spatial
coarsening is aggressive.  However,
the use of damping methods (second-order,
fourth-order, or sixth-order diffusion, spectral viscosity)
to resolve the spectral blocking problem \citep{gelb2001spectral}
brings into question the relevance of accuracy in the highest
frequencies.
We investigate this question by using different measures
of error on the SH discretization in the experiments. In the 
numerical examples, we start by evaluating the convergence
rate of PFASST-SH as the time step size is decreased. This
confirms that spatial coarsening is a key determinant of the
accuracy of the parallel scheme, and demonstrates that PFASST-SH
can capture the small-scale features of the solution when the
spatial coarsening is mild. Then, we measure the computational
cost of PFASST-SH to show that it yields a significant reduction
in the time-to-solution compared to serial single-level and
multi-level SDC schemes for the three numerical examples considered
here. A strong scalability test performed on the benchmark proposed
in \cite{galewsky2004initial} concludes our numerical study of
the properties of PFASST-SH.

In Sections
\oldref{section_governing_equations},
\oldref{section_spatial_discretization}, and
\oldref{section_temporal_splitting},
we review some features of the problem -- namely,
the governing equations, SH transform, and IMEX temporal
splitting --  already discussed in \cite{hamon2018multi}
to ensure the completeness of this paper. We present
the parallel-in-time numerical scheme based on PFASST-SH in
Section \oldref{section_temporal_discretization}. We evaluate
its performance and accuracy on SWE test cases in
Section \oldref{section_numerical_examples}.

\section{\label{section_governing_equations}Mathematical model}

In this work, we study the properties of the parallel integration
scheme applied to the Shallow-Water Equations (SWE) on the rotating
sphere. This simplified mathematical model captures the horizontal
features of realistic atmospheric flows and allows us to evaluate
the main properties of our numerical scheme -- namely, accuracy,
performance, and robustness -- on a set of well-defined test cases
\citep{williamson1992standard,galewsky2004initial}.
To overcome the singularity in the velocity field at the poles,
we use the vorticity-divergence formulation in physical space
\citep{bourke1972efficient,hack1992description}. The prognostic 
variables $\boldsymbol{U} = [ \Phi, \, \zeta, \, \delta ]^T$ are
respectively the geopotential, $\Phi$, the vorticity, $\zeta$,
and the divergence, $\delta$. The system of governing partial
differential equations reads
\begin{equation}
\frac{\partial \boldsymbol{U}}{\partial t} = \boldsymbol{\mathcal{L}}( \boldsymbol{U} )  
                                           + \boldsymbol{\mathcal{N}} ( \boldsymbol{U} ),
\label{linear_nonlinear_decomposition}
\end{equation}
where the first group in the right-hand side of
\ref{linear_nonlinear_decomposition} contains the stiff terms
involved in the  linear wave motion induced by gravitational
forces and also includes the diffusion term
\begin{equation}
\boldsymbol{\mathcal{L}} ( \boldsymbol{U} ) \equiv 
\left[ 
  \begin{array}{c}
    - \bar{\Phi} \delta + \nu \nabla^2 \Phi' \\[5pt]
    \nu \nabla^2 \zeta \\[5pt]
    - \nabla^2 \Phi + \nu \nabla^2 \delta
  \end{array} 
\right].  
\label{linear_wave_motion_induced_by_gravitational_forces}
\end{equation}
The average geopotential, $\bar{\Phi} = g \bar{h}$, is written as 
the product of the gravitational acceleration by the average height,
and $\Phi'$ is defined as $\Phi' = \Phi - \bar{\Phi}$. The diffusion
coefficient is denoted by $\nu$. 
Although spectral viscosity \citep{gelb2001spectral}
is now used in standard codes like IFS, we employ a second-order
diffusion term in the governing equations to stabilize the flow
dynamics and reduce the errors caused by under-resolved nonlinear
interaction modes. Second-order diffusion
remains acceptable here since we are primarily interested in
showing that PFASST-SH can resolve the low-frequency and
mid-frequency features of the solution significantly faster than
serial SDC schemes. This choice also allows for a fair comparison
with the results obtained previously on the same test cases with
SDC methods \citep{jia2013spectral,hamon2018multi}.
The second group in the right-hand side of
\ref{linear_nonlinear_decomposition} contains all the relatively
less stiff, nonlinear terms present in the governing equations
\begin{equation}
\boldsymbol{\mathcal{N}} ( \boldsymbol{U} ) = \left[ 
  \begin{array}{c}
 - \nabla \cdot ( \Phi' \boldsymbol{V} ) \\[5pt]
 - \nabla \cdot ( \zeta + f ) \boldsymbol{V} \\[0.pt]
   \boldsymbol{k} \cdot \nabla \times (\zeta + f) \boldsymbol{V} - \nabla^2 \bigg( \displaystyle \frac{\boldsymbol{V} \cdot \boldsymbol{V}}{2} \bigg) 
  \end{array}
  \right],
 \label{nonlinear_operators}
\end{equation}
where $\boldsymbol{k} = [0,0,1]^T$. The horizontal velocity vector is 
$\boldsymbol{V} \equiv \boldsymbol{i}u + \boldsymbol{j}v$, where $\boldsymbol{i}$ 
and $\boldsymbol{j}$ are the unit vectors in the eastward and northward directions,
respectively.  The Coriolis force is represented by 
$f = 2 \Omega \sin \phi$, where $\Omega$ is the angular rate of rotation, and 
$\phi$ is the latitude. 
In Section \oldref{section_temporal_splitting}, we use this linear-nonlinear decomposition to define an  
IMEX temporal splitting that treats the stiff terms implicitly and the less stiff terms explicitly. 

To express the velocities as a function of the 
prognostic variables, $\zeta$ and $\delta$, we first use the Helmholtz theorem 
which relates $\boldsymbol{V}$ to a scalar stream function, $\psi$, and a 
scalar velocity potential, $\chi$,
\begin{equation}
\boldsymbol{V} = \boldsymbol{k} \times \nabla \psi + \nabla \chi.
\label{helmholtz_theorem}
\end{equation}
Using the identities
\begin{align}
\zeta  &\equiv \boldsymbol{k} \cdot ( \nabla \times \boldsymbol{V} ), \label{identity_defining_zeta} \\
\delta &\equiv \nabla \cdot \boldsymbol{V}, \label{identity_defining_delta} 
\end{align}
the application of the curl and divergence operators to \ref{helmholtz_theorem} yields
$\zeta  =  \nabla^2 \psi$ and $\delta =  \nabla^2 \chi$.
The Laplace operators can be efficiently inverted using the SH transform
to compute the stream function, $\psi$, and the velocity potential, $\chi$, as a function 
of $\zeta$ and $\delta$. This is discussed in the next section, along with the spatial
discretization of $\boldsymbol{\mathcal{L}}$ and $\boldsymbol{\mathcal{N}}$.

\section{\label{section_spatial_discretization}Spectral transform method}

Here, we review the methodology based on the global SH transform
used to discretize the governing equations in space. The parallel
multi-level scheme proposed in this work heavily relies on the
structure of this spatial discretization for the construction
of a hierarchy of space-time levels. In the SH scheme, the
representation of a function of longitude $\lambda$ and Gaussian
latitude $\mu \equiv \sin( \phi )$, $\xi(\lambda,\mu)$, consists
of a sum of spherical harmonic basis functions,
$P^r_s(\mu) e^{i r \lambda}$, weighted by the spectral coefficients 
$\xi^r_s$,
\begin{equation}
\xi(\lambda, \mu) = \sum^R_{r = -R} \sum^{S(r)}_{s = |r|} \xi^r_s P^r_s(\mu) e^{i r \lambda},
\label{spectral_to_physical_space}
\end{equation}
where the indices $r$ and $s$ refer to the
zonal and total wavenumbers, respectively. In
\ref{spectral_to_physical_space}, $P^r_s$ is the normalized
associated Legendre polynomial. Without loss of generality, we
use a triangular truncation with $S(r) = R$,
where $R$ and $S(r)$ are the chosen truncation
limits of zonal and total wavenumbers, respectively. We note that
truncating the modes corresponding
to high-frequency wavenumbers -- i.e., reducing $R$ in
\ref{spectral_to_physical_space} -- 
yields a consistent coarse representation of the discrete problem.
This will be exploited in Section
\oldref{section_temporal_discretization} in the integration scheme.

The transformation from physical to spectral space is achieved
in two steps. The first step consists in taking the discrete
Fourier transform of $\xi(\lambda, \mu)$ in longitude -- i.e.,
over $\lambda$ --, followed by a second step based on the 
application of the discrete Legendre transformation in
latitude. This two-step global transform is applied to
\ref{linear_nonlinear_decomposition} to obtain a system of coupled
ordinary differential equations involving the prognostic variables
in spectral space, 
$\boldsymbol{\Theta}^r_s = [ \Phi^r_s, \, \zeta^r_s, \, \delta^r_s ]$.
Noting that, due to the symmetry of the spectral coefficients,
we only have to consider the indices $r \geq 0$, 
we obtain for $r \in \{ 0, \dots, R \}$ and
$s \in \{ r, \dots, R \}$,
\begin{equation}
\frac{\partial \boldsymbol{\Theta}^r_s}{\partial t} = \boldsymbol{L}^r_s( \boldsymbol{\Theta} )
                                                   + \boldsymbol{N}^r_s ( \boldsymbol{\Theta} ),
\label{semidiscrete_system}
\end{equation}
where $\boldsymbol{L}^r_s$, and $\boldsymbol{N}^r_s$ are the
discrete, spectral representations of the operators defined in 
\ref{linear_wave_motion_induced_by_gravitational_forces}, and
\ref{nonlinear_operators}. The state variable in spectral space,
$\boldsymbol{\Theta}$, is defined as a vector of size $K$ as follows 
\begin{equation}
\boldsymbol{\Theta} = [ \boldsymbol{\Theta}^0_0, \, \boldsymbol{\Theta}^{1}_{0}, \dots, \, \boldsymbol{\Theta}^{R}_{R-1}, \, \boldsymbol{\Theta}^R_R ]^T.
\label{spectral_space}
\end{equation}
An efficient implementation
of the global SH transform is described in 
\cite{temperton1991scalar,rivier2002efficient,schaeffer2013efficient}.
The latter reference served as a basis for the developments
presented here. 

\section{\label{section_temporal_splitting}Temporal splitting}

In atmospheric modeling, the propagation of fast waves --
e.g., sound or gravity waves -- in the system often imposes
a severe stability restriction on the time step size of fully
explicit schemes. Fully implicit schemes overcome this stability
constraint but require solving expensive global nonlinear systems
\citep{evans2010accuracy,jia2013spectral,lott2015algorithmically}. 

Instead, implicit-explicit (IMEX) schemes only evaluate implicitly
the stiff terms involved in the propagation of the fast waves,
while less stiff terms are updated explicitly. This strategy 
reduces the cost of the implicit solves and allows for relatively
large stable time steps. In non-hydrostatic atmospheric modeling,
dimensional splitting is a commonly used IMEX strategy that only
treats implicitly the terms involved in the (fast) vertical dynamics 
\citep{ullrich2012operator,weller2013runge,giraldo2013implicit,lock2014numerical,gardner2018implicit}. In
\cite{smolarkiewicz2014consistent}, an IMEX scheme for the
compressible Euler equations is obtained
by treating all terms pertinent to wave motions
(namely acoustic, gravity, and rotational) implicitly, but advection
explicitly. Alternatively, the approach of
\cite{robert1972implicit,giraldo2005semi} consists in linearizing
the governing PDEs and treating the linearized piece implicitly.
The term treated explicitly is obtained by subtracting the
linearized part from the nonlinear system. 

In this work, we directly discretize the fast linear terms on
the right-hand side of \ref{semidiscrete_system} implicitly,
while the other terms are evaluated explicitly. Specifically,
our integration scheme relies on the following splitting, written
below in semi-discrete form for $r \in \{ 0, \dots, R \}$ and
$s \in \{ r, \dots, R \}$,
\begin{equation}
\frac{\partial \boldsymbol{\Theta}^r_s}{\partial t} = (\boldsymbol{F}_I)^r_s(\boldsymbol{\Theta}) 
                                              + (\boldsymbol{F}_E)^r_s(\boldsymbol{\Theta}),
\label{system_of_odes}
\end{equation}
where 
\begin{align}
(\boldsymbol{F}_I)^r_s &\equiv (\boldsymbol{L})^r_s, \label{system_of_odes_implicit_part} \\
(\boldsymbol{F}_E)^r_s &\equiv \boldsymbol{N}^r_s \label{system_of_odes_explicit_part}.
\end{align}
That is,  the implicit right-hand side, $(\boldsymbol{F}_I)^r_s$,
contains the terms representing linear wave motion induced by
gravitational forces and the diffusion term. The explicit right-hand
side, $(\boldsymbol{F}_E)^r_s$, contains a linear harmonic oscillator
and the nonlinear terms. As  explained in
Section~\oldref{subsection_solver_for_the_implicit_systems}, in the
context of the SH spatial discretization, the implicit solve
necessary in the IMEX time-stepping approach is inexpensive compared
to the cost of the explicit evaluation of nonlinear terms.
Hence, the IMEX approach greatly reduces the computational cost
per step compared to a fully implicit method  by  circumventing
the need for a global nonlinear solver.  The IMEX approach also
allows much larger stable time steps than a fully explicit
method with little additional computational cost.
Next, we describe the parallel 
integration scheme that is used to advance the semi-discrete
system \ref{system_of_odes} in time.

\section{\label{section_temporal_discretization}Parallel-in-time integration}

We review the fundamentals of the Parallel Full Approximation Scheme
in Space and Time (PFASST) algorithm detailed in
\cite{emmett2012toward} after earlier work by \cite{minion2011hybrid}.
In the remainder of the paper, an SDC sweep refers
to the computation of a correction with respect to a previous
iteration.
The multi-level scheme consists in coupling
iteratively serial coarse SDC sweeps propagating the initial
conditions with parallel fine SDC sweeps achieving high-order
accuracy at each time step. We start from two building blocks that
are key to the accuracy and efficiency of the scheme, namely the
SDC sweeps in
Section~\oldref{subsubsection_implicit_explicit_spectral_deferred_correction}
and the Full Approximation Scheme (FAS) in
Section~\oldref{subsubsection_full_approximation_scheme}.
Then we detail the steps of PFASST and its computational cost in
Sections~\oldref{subsection_pfasst_algorithm} to \oldref{subsection_computational_cost}. We consider a system of
coupled ODEs in the generic form
\begin{align}
  \frac{\partial \boldsymbol{\Theta} }{\partial t} (t)
  &= \boldsymbol{F}_I \big( \boldsymbol{\Theta}(t) \big) + \boldsymbol{F}_E \big( \boldsymbol{\Theta}(t) \big),
  \qquad t \in [t^n,t^{n} + \Delta t], 
\label{model_problem}
\\
\boldsymbol{\Theta}(t^n) &= \boldsymbol{\Theta}^n,
\end{align}
where $\boldsymbol{F}_I$ and $\boldsymbol{F}_E$ are the implicit
and explicit right-hand sides, respectively, with
$\boldsymbol{F} = \boldsymbol{F}_I + \boldsymbol{F}_E$ and 
$\boldsymbol{\Theta}(t)$ is the state variable at time $t$. 

\subsection{\label{subsection_method_components}Ingredients of the time integration scheme}

\subsubsection{\label{subsubsection_implicit_explicit_spectral_deferred_correction}IMplicit-EXplicit Spectral Deferred Correction (IMEX SDC)}

SDC methods have been presented in \cite{dutt2000spectral} and then generalized to methods 
with different temporal splittings in \cite{minion2003semi,bourlioux2003high,layton2004conservative}.  
In the context of fast-wave slow-wave problems, the properties of IMEX SDC schemes have been
studied in \cite{ruprecht2016spectral}.
In SDC methods, the interval $[t^n, t^{n+1}]$ is decomposed into $M$ subintervals using $M+1$ Gauss-Lobatto temporal nodes, such that
\begin{equation}
  t^n \equiv t^{n,0} < t^{n,1} < \dots < t^{n,M} = t^n + \Delta t \equiv t^{n+1}.
\label{time_discretization_sdc_nodes}
\end{equation}
In the remainder of this paper, we use the shorthand notation $t^m = t^{n,m}$.
We denote by $\boldsymbol{\Theta}^{m+1,(k+1)}$ the approximate solution at node $m+1$ and at sweep $(k+1)$.
The SDC scheme applied to the implicit-explicit temporal splitting \ref{model_problem} iteratively 
improves the accuracy of the approximation based on the discrete correction equation
\begin{align}
\boldsymbol{\Theta}^{m+1,(k+1)} = \boldsymbol{\Theta}^{n} 
                              &+ \Delta t \sum_{j = 1}^m \tilde{q}^E_{m+1,j} \big[ \boldsymbol{F}_E \big( \boldsymbol{\Theta}^{j,(k+1)} \big) 
                                                                                                - \boldsymbol{F}_E \big( \boldsymbol{\Theta}^{j,(k)}   \big) \big]  \nonumber \\
                              &+ \Delta t \sum_{j = 1}^{m+1} \tilde{q}^I_{m+1,j} \big[ \boldsymbol{F}_I \big( \boldsymbol{\Theta}^{j,(k+1)} \big) 
                                                                                                - \boldsymbol{F}_I \big( \boldsymbol{\Theta}^{j,(k)}   \big) \big] \nonumber \\
&+ \Delta t \sum_{j = 0}^{M} q_{m+1,j} \boldsymbol{F} \big( \boldsymbol{\Theta}^{j,(k)} \big).
\label{update_equation_sdcq_discrete_form} 
\end{align}
The coefficients $q_{m+1,j}$ are chosen to be the Gauss-Lobatto quadrature points, such that
the term on the third line of \ref{update_equation_sdcq_discrete_form} is a high-order approximation of the
integral
\begin{equation}
\int_{t^n}^{t} \boldsymbol{F} \big( \tilde{\boldsymbol{\Theta}}(a) \big)  d a.
\end{equation}
Finally, the coefficients $\tilde{q}^E_{m+1,j}$ correspond to forward-Euler time stepping,
while the weights $\tilde{q}^I_{m+1,j}$ are chosen to be the coefficients of the upper
triangular matrix in the LU decomposition of $\boldsymbol{Q}^T$, where
$\boldsymbol{Q} = \{ q_{ij} \} \in \mathbb{R}^{(M+1)\times(M+1)}$.
We refer to \cite{weiser2015faster} for a proof that this choice leads to fast convergence of the
iterative process to the fixed-point solution for stiff problems, and to \cite{hamon2018concurrent} for numerical examples
illustrating the improved convergence. Each pass of the discrete version of the update 
equation \ref{update_equation_sdcq_discrete_form}, referred to as sweep, increases the formal order 
of accuracy by one until the order of accuracy of the quadrature applied to the third 
integral is reached \citep{christlieb2009comments}. We mention
here that our implementation splits the implicit step of \ref{update_equation_sdcq_discrete_form}
into two substeps, with a substep for the physical linear terms followed by a substep for the artificial diffusion terms.
This additional splitting is described in \cite{bourlioux2003high} and will not be discussed in the present
paper for brevity.

Using the compact notation of \cite{bolten2016multigrid}, we introduce the space-time
vectors $\vec{\boldsymbol{\Theta}} \in \mathbb{C}^{(M+1)K}$ and
$\vec{\boldsymbol{F}} \in \mathbb{C}^{(M+1)K}$ are such that
\begin{align}
\vec{\boldsymbol{\Theta}} &\equiv [\boldsymbol{\Theta}^{n,0}, \dots, \boldsymbol{\Theta}^{n,M}]^T,  \\
\vec{\boldsymbol{F}} &\equiv \vec{\boldsymbol{F}} ( \vec{\boldsymbol{\Theta}} ) = [\boldsymbol{F}( \boldsymbol{\Theta}^{n,0} ), \dots, \boldsymbol{F}( \boldsymbol{\Theta}^{n,M}) ]^T.
\end{align}
We also define the operator 
\begin{equation}
\boldsymbol{A} ( \vec{\boldsymbol{\Theta}} ) \equiv \vec{\boldsymbol{\Theta}} 
                                           - \Delta t ( \boldsymbol{Q} \otimes \boldsymbol{I}_{K} ) \vec{\boldsymbol{F}},
\end{equation}
where $\otimes$ denotes the Kronecker product and $\boldsymbol{I}_{K} \in \mathbb{R}^{K \times K}$
is the identity matrix. $\boldsymbol{1}_{M+1} \in \mathbb{R}^{M+1}$ is a vector of ones.
Using these definitions, the integration scheme \ref{update_equation_sdcq_discrete_form}
can be written as an iterative solution method for the collocation problem defined by
\begin{equation}
\boldsymbol{A} ( \vec{\boldsymbol{\Theta}} ) = \boldsymbol{1}_{M+1} \otimes \boldsymbol{\Theta}^{n,0}. \label{collocation_problem}
\end{equation}
In the following sections, we describe how the SDC sweeps are
applied iteratively on coarse and fine representations of the
problem to obtain a parallel-in-time integration scheme.

\subsubsection{\label{subsubsection_full_approximation_scheme}Full Approximation Scheme (FAS)}

The PFASST algorithm requires the definition of
a coarse approximation in space and in time of the discrete
problem of interest. We refer to the former as
the coarse level $\ell =c$, and to the latter as the fine
level $\ell = f$. The space-time coarsening methodology
used to construct the coarse space is reviewed in
Section~\oldref{subsection_coarsening_strategy}.
On these two levels, the PFASST algorithm computes
multiple time steps concurrently 
by interweaving parallel high-order SDC sweeps on the fine level
and serial SDC sweeps performed on the coarse level to propagate
the updated initial conditions between time steps (see
Section~\oldref{subsection_pfasst_algorithm}).
We do not consider the case in which PFASST involves more than
one coarse level. To ensure that the fine and coarse levels
are properly coupled, we employ the Full Approximation Scheme
\citep[FAS,][]{brandt1977multi}, described next. 
In FAS, the discrete equations solved on the coarse level are
modified with the introduction of a correction term,
$\vec{\boldsymbol{\tau}}_c$. 
The coarse collocation problem with the correction term becomes
\begin{equation}
\boldsymbol{A}_{c} ( \vec{\boldsymbol{\Theta}}_{c} ) - \vec{\boldsymbol{\tau}}_{c}
= \boldsymbol{1}_{M_{c}+1} \otimes \boldsymbol{\Theta}^{n,0}_{c},
\label{collocation_problem_level_l+1}
\end{equation}
where $\boldsymbol{A}_{c}$ is the approximation of $\boldsymbol{A}$
at the coarse level. The FAS correction term
at the coarse level is defined as 
\begin{equation}
\vec{\boldsymbol{\tau}}_{c} \equiv \boldsymbol{A}_{c} ( \boldsymbol{R}^{c}_{f} \vec{\boldsymbol{\Theta}}_{f} ) 
                           - \boldsymbol{R}^{c}_{f} \boldsymbol{A}_{f} ( \vec{\boldsymbol{\Theta}}_{f} )  
                           + \boldsymbol{R}^{c}_{f} \vec{\boldsymbol{\tau}}_{f},
\label{tau_term_level_l+1}
\end{equation}
with, for the two-level case, $\vec{\boldsymbol{\tau}}_f = \boldsymbol{0}$ on the fine level. In \ref{collocation_problem_level_l+1} and
\ref{tau_term_level_l+1}, 
$\vec{\boldsymbol{\Theta}}_{\ell} \in \mathbb{C}^{(M_{\ell}+1) K_{\ell}}$
and
$\vec{\boldsymbol{F}}_{\ell} \in \mathbb{C}^{(M_{\ell}+1) K_{\ell}}$ are
respectively the space-time vector and right-hand side at level $\ell$.
$K_{\ell}$ represents the total number of spectral coefficients
in \ref{spectral_space} on level $\ell$. The matrix
$\boldsymbol{R}^{c}_{f} \in \mathbb{R}^{(M_{c}+1)K_{c} \times  (M_{f} + 1)K_{f}}$
is the linear
restriction operator from the fine level to the coarse level
that consists in truncating the spectral coefficients
corresponding to the high-frequency modes in the SH basis.
With this modification of the coarse equations, we note that 
the restriction of the fine solution, $\boldsymbol{R}^{c}_{f} \vec{\boldsymbol{\Theta}}_{f}$, is a solution of the coarse problem. 
On the coarse problem \ref{collocation_problem_level_l+1}, the modified SDC update for temporal node $m+1$ at sweep $(k+1)$ is
\begin{align}
\boldsymbol{\Theta}^{m+1,(k+1)}_{c} &= \boldsymbol{\Theta}^{n,0}_{c} 
                                + \Delta t \sum_{j = 1}^m (\tilde{q}^E_{m+1,j})_{c} \big[ \boldsymbol{F}_{E,c} \big( \boldsymbol{\Theta}^{j,(k+1)}_{c} \big) 
                                                                                                - \boldsymbol{F}_{E, c} \big( \boldsymbol{\Theta}^{j,(k)}_{c}   \big) \big] \nonumber \\
                                &+ \Delta t \sum_{j = 1}^{m+1} (\tilde{q}^I_{m+1,j})_{c}              \big[ \boldsymbol{F}_{I, c} \big( \boldsymbol{\Theta}^{j,(k+1)}_{c} \big) 
                                                                                                - \boldsymbol{F}_{I, c} \big( \boldsymbol{\Theta}^{j,(k)}_{c}   \big) \big] \nonumber \\
                                &+ \Delta t \sum_{j = 0}^{M} (q_{m+1,j})_c \boldsymbol{F}_c \big( \boldsymbol{\Theta}^{j,(k)}_c \big) 
                                + \boldsymbol{\tau}^{m+1,(k)}_{c} \label{modified_update_equation_sdcq_discrete_form} .
\end{align}

\subsection{\label{subsection_pfasst_algorithm}Parallel Full Approximation Scheme in Space and Time (PFASST) algorithm}

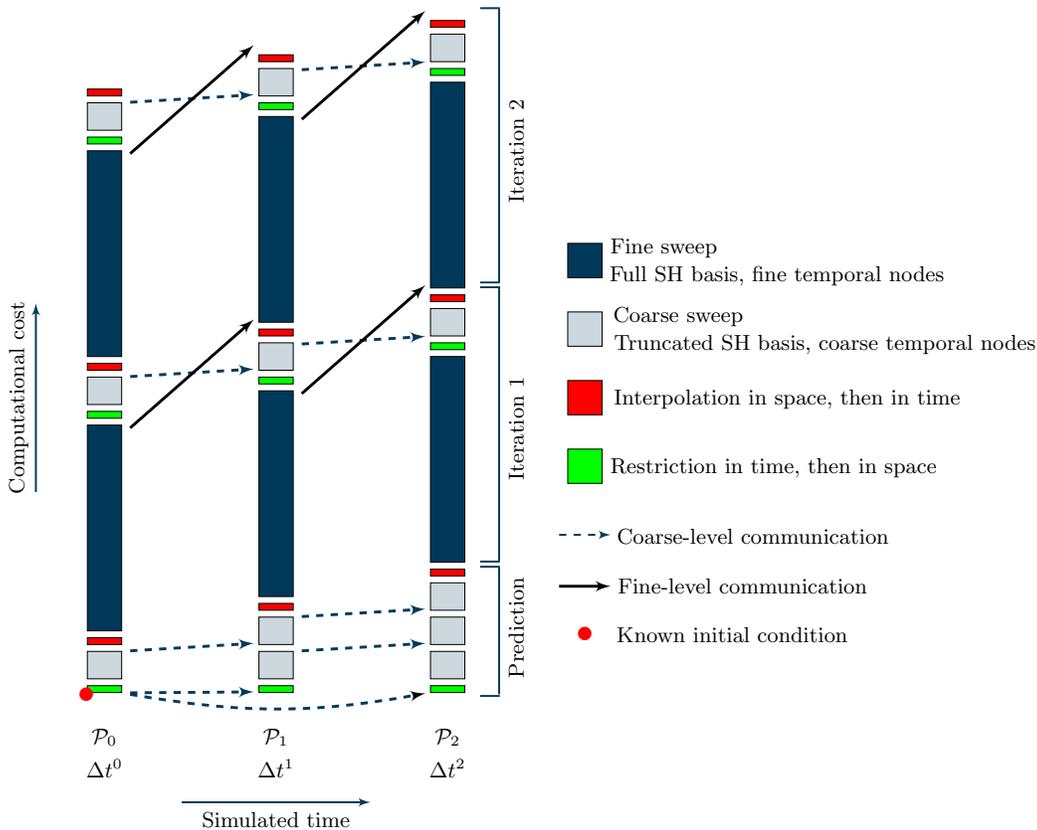
\begin{figure}[ht!]
  \centering
  \resizebox{0.90\linewidth}{!}{
    \tikzstyle{int}=[draw, minimum size=2em]
\tikzstyle{init} = [pin edge={to-,thin,black}]

\begin{tikzpicture}[node distance=3cm,auto,>=latex']

  \node (ib_0) at (-0.75,-1.7) {\small $\mathcal{P}_{0}$};
  \node (ib_0) at (-0.75,-2.15) {\small $\Delta t^{0}$};

  \draw [black,fill=green] (-1,-1.05) rectangle (-0.5,-0.95);
  \draw [black,fill=berkeley_lab_blue!20!white] (-1,-0.85) rectangle (-0.5,-0.45);
  \draw [black,fill=red] (-1,-0.35) rectangle (-0.5,-0.25);
  \draw [black,fill=berkeley_lab_blue] (-1,-0.15) rectangle (-0.5,2.85);
  \draw [black,fill=green] (-1,2.95) rectangle (-0.5,3.05);
  \draw [black,fill=berkeley_lab_blue!20!white] (-1,3.15) rectangle (-0.5,3.55);
  \draw [black,fill=red] (-1,3.65) rectangle (-0.5,3.75);
  \draw [black,fill=berkeley_lab_blue] (-1,3.85) rectangle (-0.5,6.85);   
  \draw [black,fill=green] (-1,6.95) rectangle (-0.5,7.05);
  \draw [black,fill=berkeley_lab_blue!20!white] (-1,7.15) rectangle (-0.5,7.55);
  \draw [black,fill=red] (-1,7.65) rectangle (-0.5,7.75);

  \node (ib_1) at (-0.5,-1.055) {};
  \node (ib_2) at (1.55,-1.04) {};
  \path [draw=berkeley_lab_blue, dashed, very thick,->] (ib_1) -- (ib_2);

  \node (ib_1) at (-0.5,-1.055) {};
  \node (ib_2) at (4.05,-1.04) {};
  \path [draw=berkeley_lab_blue, dashed, very thick,->,bend left] (ib_1) edge [bend right=10] node {} (ib_2);

  \node (ib_1) at (-0.5,-0.45) {};
  \node (ib_2) at (1.55,-0.325) {};
  \path [draw=berkeley_lab_blue, dashed, very thick,->] (ib_1) -- (ib_2);

  \node (ib_1) at (-0.5,3.55) {};
  \node (ib_2) at (1.55,3.675) {};
  \path [draw=berkeley_lab_blue, dashed, very thick,->] (ib_1) -- (ib_2);

  \node (ib_1) at (-0.5,2.7) {};                
  \node (ib_2) at (1.55,4.5) {};        
  \path [draw=black, very thick,->] (ib_1) -- (ib_2);

  \node (ib_1) at (-0.5,7.55) {};
  \node (ib_2) at (1.55,7.675) {};
  \path [draw=berkeley_lab_blue, dashed, very thick,->] (ib_1) -- (ib_2);

  \node (ib_1) at (-0.5,6.7) {};                
  \node (ib_2) at (1.55,8.5) {};        
  \path [draw=black, very thick,->] (ib_1) -- (ib_2);

  \node (ib_0) at (1.75,-1.7) {\small $\mathcal{P}_{1}$};
  \node (ib_0) at (1.75,-2.15) {\small $\Delta t^{1}$};

  \draw [black,fill=green] (1.5,-1.05) rectangle (2.0,-0.95);
  \draw [black,fill=berkeley_lab_blue!20!white] (1.5,-0.85) rectangle (2,-0.45);   
  \draw [black,fill=berkeley_lab_blue!20!white] (1.5,-0.35) rectangle (2,0.05);
  \draw [black,fill=red] (1.5,0.15) rectangle (2,0.25);
  \draw [black,fill=berkeley_lab_blue] (1.5,0.35) rectangle (2,3.35);
  \draw [black,fill=green] (1.5,3.45) rectangle (2,3.55);
  \draw [black,fill=berkeley_lab_blue!20!white] (1.5,3.65) rectangle (2,4.05);
  \draw [black,fill=red] (1.5,4.15) rectangle (2,4.25);
  \draw [black,fill=berkeley_lab_blue] (1.5,4.35) rectangle (2,7.35);
  \draw [black,fill=green] (1.5,7.45) rectangle (2,7.55);
  \draw [black,fill=berkeley_lab_blue!20!white] (1.5,7.65) rectangle (2,8.05);
  \draw [black,fill=red] (1.5,8.15) rectangle (2,8.25);

  \node (ib_1) at (2.0,-0.45) {};                
  \node (ib_2) at (4.05,-0.325) {};        
  \path [draw=berkeley_lab_blue, dashed, very thick,->] (ib_1) -- (ib_2);

  \node (ib_1) at (2.0,0.05) {};                
  \node (ib_2) at (4.05,0.175) {};        
  \path [draw=berkeley_lab_blue, dashed, very thick,->] (ib_1) -- (ib_2);

  \node (ib_1) at (2.0,4.025) {};
  \node (ib_2) at (4.05,4.15) {};
  \path [draw=berkeley_lab_blue, dashed, very thick,->] (ib_1) -- (ib_2);

  \node (ib_1) at (2.0,3.2) {};                
  \node (ib_2) at (4.05,5) {};        
  \path [draw=black, very thick,->] (ib_1) -- (ib_2);      

  \node (ib_1) at (2.0,8.025) {};
  \node (ib_2) at (4.05,8.15) {};
  \path [draw=berkeley_lab_blue, dashed, very thick,->] (ib_1) -- (ib_2);

  \node (ib_1) at (2.0,7.2) {};                
  \node (ib_2) at (4.05,9) {};        
  \path [draw=black, very thick,->] (ib_1) -- (ib_2);

  \node (ib_0) at (4.25,-1.7) {\small $\mathcal{P}_{2}$};
  \node (ib_0) at (4.25,-2.15) {\small $\Delta t^{2}$};

  \draw [black,fill=green] (4,-1.05) rectangle (4.5,-0.95);
  \draw [black,fill=berkeley_lab_blue!20!white] (4,-0.85) rectangle (4.5,-0.45);
  \draw [black,fill=berkeley_lab_blue!20!white] (4,-0.35) rectangle (4.5,0.05);
  \draw [black,fill=berkeley_lab_blue!20!white] (4, 0.15) rectangle (4.5,0.55);
  \draw [black,fill=red] (4,0.65) rectangle (4.5,0.75);
  \draw [black,fill=berkeley_lab_blue] (4,0.85) rectangle (4.5,3.85);
  \draw [black,fill=green] (4,3.95) rectangle (4.5,4.05);
  \draw [black,fill=berkeley_lab_blue!20!white] (4,4.15) rectangle (4.5,4.55);
  \draw [black,fill=red] (4,4.65) rectangle (4.5,4.75);
  \draw [black,fill=berkeley_lab_blue] (4,4.85) rectangle (4.5,7.85); 
  \draw [black,fill=green] (4,7.95) rectangle (4.5,8.05);
  \draw [black,fill=berkeley_lab_blue!20!white] (4,8.15) rectangle (4.5,8.55);
  \draw [black,fill=red] (4,8.65) rectangle (4.5,8.75);


  \node (ib_1) at (-1.75,1.75) {};                
  \node (ib_2) at (-1.75,4.75) {};        
  \path [draw=berkeley_lab_blue, thick,->] (ib_1) -- (ib_2);      
  \node[rotate=90] (ib_0) at (-2,3.25) {\small Computational cost};

  \node (ib_1) at (0.25,-2.65) {};                
  \node (ib_2) at (3.25,-2.65) {};        
  \path [draw=berkeley_lab_blue, thick,->] (ib_1) -- (ib_2);      
  \node (ib_0) at (1.75,-2.9) {\small Simulated time};

  \node (ib_1) at (5, 0.9) {};                
  \node (ib_2) at (5,-1.2) {};        
  \path [draw=berkeley_lab_blue, thick] (ib_1) -- (ib_2);      
  \node[rotate=90] (ib_0) at (5.25,-0.1) {\small Prediction};

  \node (ib_1) at (5.135, 0.785) {};                
  \node (ib_2) at (4.6,0.785) {};        
  \path [draw=berkeley_lab_blue, thick] (ib_1) -- (ib_2);      

  \node (ib_1) at (5.135,-1.1) {};                
  \node (ib_2) at (4.6,-1.1) {};        
  \path [draw=berkeley_lab_blue, thick] (ib_1) -- (ib_2);

  \node (ib_1) at (5, 0.75) {};                
  \node (ib_2) at (5, 5.) {};        
  \path [draw=berkeley_lab_blue, thick] (ib_1) -- (ib_2);      
  \node[rotate=90] (ib_0) at (5.25,2.8) {\small Iteration 1};

  \node (ib_1) at (5.135, 0.86) {};                
  \node (ib_2) at (4.6,0.86) {};        
  \path [draw=berkeley_lab_blue, thick] (ib_1) -- (ib_2);      

  \node (ib_1) at (5.135,4.86) {};                
  \node (ib_2) at (4.6,4.86) {};        
  \path [draw=berkeley_lab_blue, thick] (ib_1) -- (ib_2);

  \node (ib_1) at (5, 4.8) {};                
  \node (ib_2) at (5, 9.05) {};        
  \path [draw=berkeley_lab_blue, thick] (ib_1) -- (ib_2);      
  \node[rotate=90] (ib_0) at (5.25,6.85) {\small Iteration 2};

  \node (ib_1) at (5.135,8.93) {};                
  \node (ib_2) at (4.6,8.93) {};        
  \path [draw=berkeley_lab_blue, thick] (ib_1) -- (ib_2);      

  \node (ib_1) at (5.135,4.93) {};                
  \node (ib_2) at (4.6,4.93) {};        
  \path [draw=berkeley_lab_blue, thick] (ib_1) -- (ib_2);

  \draw [black,fill=berkeley_lab_blue] (6,5) rectangle (6.5,5.5);
  \node (ib_2) at (7.4,5.425) {\small Fine sweep};
  \node (ib_2) at (9.05,5.025) {\small Full SH basis, fine temporal nodes};

  \draw [black,fill=berkeley_lab_blue!20!white] (6,4) rectangle (6.5,4.5);
  \node (ib_2) at (7.61,4.425) {\small Coarse sweep};
  \node (ib_2) at (9.75,4.025) {\small Truncated SH basis, coarse temporal nodes};

  \draw [black,fill=red] (6,3.5) rectangle (6.5,3);
  \node (ib_2) at (9.2,3.235) {\small Interpolation in space, then in time};

  \draw [black,fill=green] (6,2) rectangle (6.5,2.5);
  \node (ib_2) at (9,2.235) {\small Restriction in time, then in space};

  \node (ib_1) at (5.75,1.25) {};                
  \node (ib_2) at (6.75,1.25) {};        
  \path [draw=berkeley_lab_blue,dashed,  thick,->] (ib_1) -- (ib_2);
  \node (ib_2) at (8.7,1.235) {\small Coarse-level communication};

  \node (ib_1) at (5.75,0.5) {};                
  \node (ib_2) at (6.75,0.5) {};        
  \path [draw=black, very thick,->] (ib_1) -- (ib_2);      
  \node (ib_2) at (8.55,0.5) {\small Fine-level communication};

  \node (ib_2) at (8.4,-0.2) {\small Known initial condition};
  \node (ib_1) at (6.25,-0.2) {\Large \textcolor{red}{$\bullet$}};

  \node (ib_1) at (-1.015,-1.0825) {\Large \textcolor{red}{$\bullet$}};

\end{tikzpicture}
  }
\vspace{-0.1cm}
\caption{\label{fig:sketch_pfasst_iteration}
Sketch of the PFASST algorithm described in Section
\oldref{subsection_pfasst_algorithm}, including the prediction
step and the first two iterations for a block of three time steps.
Processor $\mathcal{P}_n$ owns time step $\Delta t^n$
($0 \leq n \leq 2$). The PFASST iteration starts with a sweep at
the fine level (dark blue). The solution is then restricted to the
coarse level (green), a coarse sweep is performed (light blue),
and the new coarse solution is interpolated back to the fine level. 
}
\end{figure}

We present the steps of the PFASST algorithm presented in
\cite{emmett2012toward}, with the improved communication pattern
of \cite{emmett2014efficient}. The presentation of the algorithm
is done for the two-level case that we exclusively
consider in this work. We consider that $n_{\textit{ts}}$
processors are available to solve a block of $n_{\textit{ts}}$ time
steps in parallel and we assume that the $n^{\text{th}}$ time step
is assigned to processor, or set of computing resources,
$\mathcal{P}_n$. We denote by $\boldsymbol{\Theta}^{n,m,(k)}$ the
approximate solution at time step $n$, at SDC node index $m$, and at
PFASST iteration $(k)$. In compact form,
$\vec{\boldsymbol{\Theta}}^{n,(k)}$ is the space-time vector
containing the approximate solution at time step $n$, at all
temporal nodes, and PFASST iteration $(k)$. This
algorithm, illustrated in Fig.~\oldref{fig:sketch_pfasst_iteration}
on a block of three time steps solved in parallel, starts with a
prediction step, followed by a sequence of iterations to correct
the prediction.

\subsubsection{Prediction step}

When the algorithm starts, the only known initial
condition -- that is, the solution at SDC node $m=0$ on a given
time step -- is at the first time step of the block, denoted by
$\Delta t^0$, and owned by processor $\mathcal{P}_0$. The 
prediction procedure aims at using serial coarse SDC sweeps
to generate an approximate fine solution for all the SDC nodes of
the time steps of the block.

At the beginning of the prediction step, processor
$\mathcal{P}_0$ sends its initial condition to the other processors.
Then, each processor restricts in space the fine initial
condition received from $\mathcal{P}_0$ using the methodology
described in Section~\oldref{subsection_coarsening_strategy}. After
that, processor $\mathcal{P}_n$ ($0 \leq n \leq n_{\textit{ts}}-1$)
performs $n+1$ relatively inexpensive coarse sweeps with
an updated initial condition (SDC node $m = 0$) received from
processor $\mathcal{P}_{n-1}$ before each sweep. A coarse sweep
consists in
applying \ref{modified_update_equation_sdcq_discrete_form}
at each coarse node of its time
interval. Finally, each processor interpolates the solution
resulting from the coarse sweeps to the fine level to
generate the fine prediction for the PFASST iterations. The
fine implicit and explicit right-hand sides are then reevaluated
using this fine prediction. We refer the reader to
\cite{emmett2012toward} for a detailed
description of the prediction step.

This prediction procedure does not require fine
sweeps and is therefore inexpensive compared to a
PFASST iteration. We note that the prediction step introduces a
load imbalance in the algorithm since the processors do not perform
the same number of coarse sweeps. This imbalance, illustrated in
Fig.~\oldref{fig:sketch_pfasst_iteration}, remains small
compared to  the total cost of the algorithm.
Without any predictor, the load imbalance would
still appear in the first iteration after the coarse sweep due
to its serial nature (see below). The idea of the predictor is
to use this fact to improve the initial guess for each time step.

\subsubsection{\label{subsubsection_pfasst_iteration}PFASST iteration}
  
After the prediction step, a sequence of iterations
is used to improve the quality of the predicted solution and achieve
high-order accuracy on the time intervals with expensive,
but parallel, fine SDC sweeps. The PFASST iteration is illustrated in
Fig.~\oldref{fig:sketch_pfasst_iteration} and written in pseudocode
in Algorithm~\oldref{alg:pfasst_iteration}. We use the PFASST
algorithm with a fixed number of iterations. The behavior of the
algorithm with a variable number of iterations and a convergence
check to stop the iterations will be studied in future work.

We describe below the steps
taken by processor $\mathcal{P}_n$
($0 \leq n \leq n_{\textit{ts}}-1$) during one PFASST iteration.
The implementation of the algorithm on a
high-performance computing platform as well as the optimization
of the communication pattern is discussed in
\cite{emmett2014efficient}.

\resizebox{0.9175\linewidth}{!}{
\begin{algorithm}[H]
  \SetAlgoLined
  \caption{\label{alg:pfasst_iteration}PFASST iteration on a fine and a coarse level} for processor $\mathcal{P}_{n}$.
  \BlankLine

 \vspace{-0.1cm}
 \KwData{Initial data $\boldsymbol{\Theta}^{n,0,(k)}_{f}$ and function evaluations $\vec{\boldsymbol{F}}^{n,(k)}_{I, f}$, $\vec{\boldsymbol{F}}^{n,(k)}_{E, f}$ from the previous PFASST 
 iteration $(k)$ on the fine level.}
 \KwResult{Approximate solution $\vec{\boldsymbol{\Theta}}^{n,(k+1)}_{\ell}$ and function evaluations $\vec{\boldsymbol{F}}^{n,(k+1)}_{I, \ell}$, $\vec{\boldsymbol{F}}^{n,(k+1)}_{E, \ell}$ on the fine and coarse levels}.
 
 \vspace{0.4cm}

 \textit{\textbf{A)} Perform a fine sweep and send to next processor} \\[2pt]
 $\vec{\boldsymbol{\Theta}}^{n,(k+1)}_{f}, \, \vec{\boldsymbol{F}}^{n,(k+1)}_{I, f}, \, \vec{\boldsymbol{F}}^{n,(k+1)}_{E, f} \longleftarrow \textbf{SweepFine}\big( \vec{\boldsymbol{\Theta}}^{n,(k)}_{f}, \, \vec{\boldsymbol{F}}^{n,(k)}_{I, f}, \, \vec{\boldsymbol{F}}^{n,(k)}_{E, f} \big) $ \\[2pt]

 \If{$(k+1)$ is last iteration}{
   return
 }
 
 $\textbf{SendLastNodeValueToNextProc}\big( \boldsymbol{\Theta}^{n, M_f ,(k+1)}_{f} \big)$ \\
 
 \vspace{0.5cm}

   \textit{\textbf{B)} Restrict, re-evaluate, save restriction, and compute FAS correction} \\ 
   \For{$m = 0, \dots, M_c$}{
     $\boldsymbol{\Theta}^{n,m,(k)}_{c} \longleftarrow \textbf{Restrict} \big( \boldsymbol{\vec{\Theta}}^{n,(k+1)}_{f} \big)$ \\
     $\boldsymbol{F}^{n,m,(k)}_{I, c}, \, \boldsymbol{F}^{n,m,(k)}_{E, c} \longleftarrow \textbf{Evaluate\_F} \big( \boldsymbol{\Theta}^{n,m,(k)}_{c} \big)$ \\
     $\boldsymbol{\tilde{\Theta}}^{n,m,(k)}_{c} \longleftarrow \boldsymbol{\Theta}^{n,m,(k)}_{c}$ \\
     $\boldsymbol{\tilde{F}}^{n,m,(k)}_{I,c}, \, \boldsymbol{\tilde{F}}^{n,m,(k)}_{E,c}  \longleftarrow \boldsymbol{F}^{n,m,(k)}_{I,c}, \, \boldsymbol{F}^{n,m,(k)}_{E,c}$
   }
   
   \vspace{-0.2cm}

   $\boldsymbol{\tau}_{c} \longleftarrow \text{FAS} \big( \vec{\boldsymbol{F}}^{n,(k)}_{I, f}, \, \vec{\boldsymbol{F}}^{n,(k)}_{E, f}, \, \vec{\boldsymbol{F}}^{n,(k)}_{I, c}, \, \vec{\boldsymbol{F}}^{n,(k)}_{E, c}, \,  \boldsymbol{\tau}_{f} \big)$

   \vspace{0.2cm}

   \textit{\textbf{C)} Receive new initial condition, sweep, and send to next processor} \\[2pt]
   $\textbf{ReceiveInitialConditionFromPreviousProc}\big( \boldsymbol{\Theta}^{n-1, M_c, (k+1)}_{c} \big)$ \\[4pt]
   $\vec{\boldsymbol{\Theta}}^{n,(k+1)}_{c}, \, \vec{\boldsymbol{F}}^{n,(k+1)}_{I, c}, \, \vec{\boldsymbol{F}}^{n,(k+1)}_{E, c} \longleftarrow \textbf{SweepCoarse} \big( \vec{\boldsymbol{\Theta}}^{n,(k)}_{c}, \, \vec{\boldsymbol{F}}^{n,(k)}_{I, c}, \, \vec{\boldsymbol{F}}^{n,(k)}_{E, c}, \, \boldsymbol{\tau}_{c} \big)$ \\[2pt]
   $\textbf{SendLastNodeValueToNextProc} \big( \boldsymbol{\Theta}^{n, M_c, (k+1)}_{c} \big)$ 

 \vspace{0.6cm}
 
 \textit{\textbf{D)} Return to finest level and receive new fine initial condition before next iteration} \\
 \For{$m = 0, \dots, M_f$}{
   $\boldsymbol{\Theta}^{n, m, (k+1)}_{f} \longleftarrow \boldsymbol{\Theta}^{n, m, (k+1)}_{f} + \textbf{Interpolate} \big( \boldsymbol{\vec{\Theta}}^{n, (k+1)}_{c} - \boldsymbol{\tilde{\vec{\Theta}}}^{n, (k)}_{c} \big)$ \\
   $\boldsymbol{F}^{n, m, (k+1)}_{I, f} \longleftarrow \boldsymbol{F}^{n, m, (k+1)}_{I, f} + \textbf{Interpolate} \big( \boldsymbol{\vec{F}}^{n, (k+1)}_{I, c} - \boldsymbol{\tilde{\vec{F}}}^{n, (k)}_{I, c}\big)$ \\ 
   $\boldsymbol{F}^{n, m, (k+1)}_{E, f} \longleftarrow\boldsymbol{F}^{n, m, (k+1)}_{E, f} +  \textbf{Interpolate} \big( \boldsymbol{\vec{F}}^{n, (k+1)}_{E, c} - \boldsymbol{\tilde{\vec{F}}}^{n, (k)}_{E, c}\big)$ 
 }
 $\textbf{ReceiveInitialConditionFromPreviousProc}\big( \boldsymbol{\Theta}^{n-1, M_f, (k+1)}_f \big)$ \\ 
 $\boldsymbol{\Theta}^{n, 0, (k+1)}_{f} \longleftarrow \boldsymbol{\Theta}^{n-1, M_f, (k+1)}_{f} + \textbf{Interpolate} \big( \boldsymbol{\Theta}^{n, m, (k+1)}_{c} - \boldsymbol{\tilde{\Theta}}^{n, m, (k)}_{c} \big)$ 

\end{algorithm}
}

\begin{enumerate}
\item[\textbf{\textit{A})}] The iteration starts with a sweep on
the fine level. This procedure consists in applying the discrete
correction \ref{update_equation_sdcq_discrete_form} at each
fine SDC node and is the most computationally
expensive step in the iteration. For the last
iteration of a time step, we skip Steps \textbf{\textit{B}},
\textbf{\textit{C}}, and \textbf{\textit{D}} and we stop
Algorithm~\oldref{alg:pfasst_iteration} here. This approach
is more robust than stopping the iterations after the coarse
calculations and interpolations (Step \textbf{\textit{D}})
when the coarsening strategy is very aggressive. If
this is not the last iteration, the updated fine value resulting
from this sweep at the last SDC node ($m = M_f$) is sent to
the next processor, $\mathcal{P}_{n+1}$.

\item[\textbf{\textit{B})}] The fine approximate
solution is restricted in time and in space to the coarse level.
Then, the coarse
implicit and explicit right-hand sides are reevaluated at all the
coarse SDC nodes using this restricted approximated
solution. We also compute the FAS correction. These quantities
will be used to compute the coarse correction during the sweep
of Step \textbf{\textit{C}}. Importantly, we also save the
solution and the right-hand sides at all the SDC nodes to
interpolate in space and in time the coarse change in these
quantities to the fine level at Step \textbf{\textit{D}}.

\item[\textbf{\textit{C})}] Processor
$\mathcal{P}_n$ receives the updated coarse initial condition from
$\mathcal{P}_{n-1}$. To reflect this change in initial condition, the
coarse right-hand sides at the first SDC node of this time interval
are reevaluated. At this point, we are ready to perform a coarse
sweep. This sweep yields a new coarse value at the last SDC node of
time step $n$ that is sent to $\mathcal{P}_{n+1}$ to conclude Step
$\textbf{\textit{C}}$.

\item[\textbf{\textit{D})}] We interpolate the
coarse solution update in space and in time to the fine level.
Then, instead of reevaluating the right-hand sides on the fine
level, we interpolate the coarse change in right-hand sides to
the fine level using the data saved at Step \textbf{\textit{B}}.
The same choice was made in \cite{hamon2018multi}. Our numerical
tests show that avoiding these fine function
evaluations reduces the computational cost without
significantly undermining the accuracy of the scheme. After the
interpolation, the new fine initial condition is
received from $\mathcal{P}_{n-1}$. For consistency, we
conclude the PFASST iteration by the application
of the interpolated coarse increment computed for SDC node $m = 0$
to this new initial condition. This initial condition will be used
at the next iteration performed by $\mathcal{P}_n$.

\end{enumerate}

For the special case of $n=0$, the processor in charge of the first time
step, $\mathcal{P}_0$, follows Algorithm~\oldref{alg:pfasst_iteration}
but skips the \textit{ReceiveInitialCondition} steps.
Next, we describe the restriction and interpolation
operators employed in Steps \textbf{\textit{B}} and
\textbf{\textit{D}}, respectively.

\subsection{\label{subsection_coarsening_strategy}Construction of the coarse level}

The coarsening strategy used to construct the space-time level
$\ell = c$ has a strong impact on the performance and accuracy
of the PFASST algorithm. The goal is to obtain a coarse problem
on which the SDC sweeps are inexpensive but can still capture
most of the features of the fine solution. This is a challenge
for nonlinear shallow-water problems whose temporal evolution
is characterized by the progressive amplification of
the modes corresponding to 
high-frequency wavenumbers that cannot be easily resolved by the
coarse sweeps. Therefore, the presence of these high-frequency
modes in the fine solution will inherently limit the spatial
coarsening that can be applied to this class of problems.

Here, we adopt the strategy proposed in \cite{hamon2018multi}
to construct the coarse level and transfer data between
levels. The spatial coarsening strategy is performed entirely
in spectral space and acts directly on the spectral basis.
This approach is based on mode truncation and
zero-padding and  does not introduce spurious modes in the
approximated solution. The restriction of the approximate
solution and right-hand sides from the fine level to the
coarse level is performed in two steps. This operation,
$\boldsymbol{R}^{c}_{f}$, can be decomposed into a restriction 
in time followed by a restriction in space as follows
\begin{equation}
\vec{\boldsymbol{\Theta}}_{c} = \boldsymbol{R}^{c}_{f} \vec{\boldsymbol{\Theta}}_{f} 
                                 = ( \boldsymbol{R}^{\textit{s}} )^{c}_{f} ( \boldsymbol{R}^{\textit{t}} )^{c}_{f} \vec{\boldsymbol{\Theta}}_{f},
\label{two_step_restriction}
\end{equation}
where $( \boldsymbol{R}^{\textit{t}} )^{c}_{f}$ and $( \boldsymbol{R}^{\textit{s}} )^{c}_{f}$ represents the
temporal and spatial operator, respectively. As explained above, $\vec{\boldsymbol{\Theta}}_{\ell} \in \mathbb{C}^{(M_{\ell}+1)K_{\ell}}$
is the space-time vector storing the state of the system at level $\ell$, where $K_{\ell}$ denotes the number of spectral
coefficients used in \ref{spectral_space} at level $\ell$ and $M_{\ell}+1$ is the number of SDC nodes at level $\ell$.

In \ref{two_step_restriction}, the restriction operator in time is defined using the Kronecker product as
\begin{equation}
( \boldsymbol{R}^{\textit{t}} )^{c}_{f} \equiv \boldsymbol{\Pi}^{c}_{f} \otimes \boldsymbol{I}_{K_{f}} \in \mathbb{R}^{(M_{c}+1)K_{f} \times (M_{f}+1)K_{f}},
\label{time_restriction}
\end{equation}
where $\boldsymbol{I}_{K_{f}} \in \mathbb{R}^{K_{f} \times K_{f}} $ is the identity matrix, and 
$\boldsymbol{\Pi}^{c}_{f} \in \mathbb{R}^{(M_{c}+1) \times (M_{f}+1)}$ is the rectangle matrix employed to interpolate a scalar function 
from the fine temporal discretization to the coarse temporal discretization. This matrix is defined using the Lagrange polynomials 
$L^{j}_{f}$ on the fine temporal discretization as follows
\begin{equation}
(\boldsymbol{\Pi}^{c}_{f})_{ij} = L^{j-1}_{f}(t^{i-1}_{c}),
\end{equation}
using the SDC node $i-1$ at the coarse level, denoted by $t^{i-1}_{c}$. In the numerical examples, we will limit the analysis to
the special cases of two, three, and five Gauss-Lobatto nodes. For this choice, applying this restriction operator in time is
equivalent to performing pointwise injection. The restriction in space requires choosing the number of spectral coefficients
that can be represented on the coarse level, denoted by $K_c$. Then, we truncate the spectral representation of the primary
variables \ref{spectral_space} in the SH transform to remove $K_f - K_c$ spectral coefficients corresponding to the high-frequency
spatial features of the approximate solution. This is achieved by applying the matrix
\begin{equation}
( \boldsymbol{R}^s )^{c}_{f} \equiv \boldsymbol{I}_{M_{c}+1} \otimes \boldsymbol{D}^{c}_{f} \in \mathbb{R}^{(M_{c}+1)K_{c} \times (M_{c} + 1)K_{f} }.
\label{space_restriction}
\end{equation}
In \ref{space_restriction}, $\boldsymbol{D}^{c}_{f} \in \mathbb{R}^{K_{c} \times K_{f}}$ is a rectangle truncation matrix defined
as such that $(\boldsymbol{D}^{c}_{f})_{ii} = 1$, and $(\boldsymbol{D}^{c}_{f})_{ij} = 0$ whenever $i \neq j$.

The interpolation procedure used to transfer the approximate solution from the coarse level to the fine level is also defined
in two steps. It begins with the application of the interpolation operator in space, $(\boldsymbol{P}^s)^{f}_{c}$ and
continues with the application of the interpolation operator in time, $(\boldsymbol{P}^t)^{f}_{c}$,
\begin{equation}
\vec{\boldsymbol{\Theta}}_{f} \equiv \boldsymbol{P}^{f}_{c} \vec{\boldsymbol{\Theta}}_{c} = (\boldsymbol{P}^t)^{f}_{c} (\boldsymbol{P}^{s})^{f}_{c} \vec{\boldsymbol{\Theta}}_{c}.
\end{equation}
The interpolation operator in space is defined as the transpose of the restriction operator in space:
\begin{equation}
( \boldsymbol{P}^s )^{f}_{c} \equiv \big( ( \boldsymbol{R}^s )^{c}_{f} \big)^T \in \mathbb{R}^{(M_{c}+1)K_{f} \times (M_{c}+1)K_{c}}.
\end{equation}
These operations are realized in an efficient way by padding the spectral representation of the primary variables at the coarse level with $K_{f} - K_{c}$  zeros.
Finally, the interpolation operator in time is analogous to \ref{time_restriction} and reads 
\begin{equation}
( \boldsymbol{P}^t )^{f}_{c} \equiv \boldsymbol{\Pi}^{f}_{c} \otimes \boldsymbol{I}_{K_{f}} \in \mathbb{R}^{(M_{f}+1)K_{f} \times (M_{c}+1)K_{f}},
\label{time_interpolation}
\end{equation}
where the rectangle interpolation matrix $\boldsymbol{\Pi}^{f}_{c}$ is constructed with the Lagrange polynomials 
$L^j_{c}$ on the coarse temporal discretization. For two, three, and five Gauss-Lobatto nodes,
\ref{time_interpolation} amounts to performing pointwise injection at the fine nodes that correspond to the 
coarse nodes, and then polynomial interpolation to compute the solution at the remaining fine nodes. In Section
\oldref{section_numerical_examples}, we will show that the accuracy of PFASST-SH is heavily dependent on the choice of
the spatial coarsening ratio in the coarsening step, since this parameter determines the range of spectral modes that can be captured on the coarse level.

\subsection{\label{subsection_solver_for_the_implicit_systems}Solver for the implicit systems}

The parallel-in-time integration scheme entails solving implicit systems at each SDC node such that $m > 0$. They
are in the form
\begin{equation}
\boldsymbol{\Theta}^{m+1,(k+1)} - \Delta t \tilde{q}^I_{m+1,m+1} \underbrace{\boldsymbol{F}_I( \boldsymbol{\Theta}^{m+1,(k+1)} )}_{\text{linear}} = \boldsymbol{b}.
\label{implicit_systems}
\end{equation}
In \ref{implicit_systems}, $\boldsymbol{b}$ is obtained from
\ref{update_equation_sdcq_discrete_form} or
\ref{modified_update_equation_sdcq_discrete_form}. For simplicity,
we have dropped the subscripts denoting the space-time 
levels. 
The structure of the linear right-hand side in
\ref{system_of_odes_implicit_part}
is determined by the spatial discretization and
temporal splitting described above. As in \cite{hamon2018multi},
this structure is exploited to circumvent the need for a
linear solver and efficiently get the updated solution in
spectral space via local updates only. We refer
the reader to \cite{hamon2018multi} for a presentation of this
solution strategy.

\subsection{\label{subsection_computational_cost}Computational cost of PFASST-SH}

Next we compare the theoretical computational cost of PFASST-SH
to those of the serial SDC and MLSDC-SH schemes
\citep{hamon2018multi}. The PFASST-SH integration scheme with
a block of $n_{\textit{ts}}$ time steps solved in parallel --
which is assumed here to correspond to the number of processors --,
$M_f + 1$ temporal nodes on the fine level, $M_c + 1$ temporal
nodes on the coarse level, $N_{\textit{PF}}$ iterations, and a
spatial coarsening ratio between the
number of coarse zonal wavenumbers and the number of fine
zonal wavenumbers, $\alpha = R_c / R_f$, is denoted by
PFASST($n_{\textit{ts}}$, $M_f+1$, $M_c+1$, $N_{\textit{PF}}$,
$\alpha$). We refer to the single-level SDC scheme with 
$M_f+1$ temporal nodes and $N_S$ fine sweeps as SDC($M_f+1$, $N_S$).
We denote by MLSDC($M_f+1$, $M_c+1$, $N_{ML}$, $\alpha$) the
MLSDC-SH scheme with $M_f+1$ nodes on the fine level, $M_c + 1$
nodes of the coarse level, $N_{\textit{ML}}$ iterations, and a
spatial coarsening ratio of $\alpha$. This notation for the
PFASST-SH, SDC, and MLSDC-SH schemes is summarized in
Tables~\oldref{tbl:sdc_mlsdc_overview} and
\oldref{tbl:pfasst_overview}. 

\begin{table}[!ht]
\begin{tabular}{|c|c|}
  \hline
  \multicolumn{2}{|c|}{SDC($M_f+1$, $N_{\textit{S}}$)} \\  
  \hline
  \textbf{Parameter} & \multicolumn{1}{c|}{\textbf{Description}} \\
  \hline
  $M_f+1$	     & SDC nodes on fine level\\ 
  \hline
  $N_{S}$            & SDC iterations  \\ 
  \hline
\end{tabular}
\hfill
\begin{tabular}{|c|c|}
  \hline
  \multicolumn{2}{|c|}{MLSDC($M_f+1$, $M_c+1$, $N_{\textit{ML}}$, $\alpha$)} \\  
  \hline
  \textbf{Parameter} & \multicolumn{1}{c|}{\textbf{Description}} \\
  \hline
  $M_f+1$	     & SDC nodes on fine level\\ 
  \hline
  $M_c+1$   	     & SDC nodes on coarse level\\ 
  \hline
  $N_{ML}$           & MLSDC iterations  \\ 
  \hline
  $\alpha$           & Spatial coarsening ratio \\ 
  \hline
\end{tabular}
\caption{\label{tbl:sdc_mlsdc_overview}Parameters for the serial SDC and MLSDC-SH schemes. The SDC iteration only involves one
  sweep on the fine level. The MLSDC-SH iteration is described in \cite{hamon2018multi}, and involves one sweep on the fine
  level and one sweep on the coarse level.} 
\end{table}

\begin{table}[!ht]
\begin{center}
\begin{tabular}{|c|c|}
  \hline
  \multicolumn{2}{|c|}{PFASST($n_{\textit{ts}}$, $M_f+1$, $M_c+1$, $N_{\textit{PF}}$, $\alpha$)} \\  
  \hline
  \textbf{Parameter} & \multicolumn{1}{c|}{\textbf{Description}} \\
  \hline
  $n_{\textit{ts}}$  & Time steps solved in parallel \\ 
  \hline
  $M_f+1$	     & SDC nodes on fine level\\ 
  \hline
  $M_c+1$   	     & SDC nodes on coarse level\\ 
  \hline
  $N_{\textit{PF}}$  & PFASST iterations  \\ 
  \hline
  $\alpha$           & Spatial coarsening ratio \\ 
  \hline
\end{tabular}
\end{center}
\vspace{-0.5cm}
\caption{\label{tbl:pfasst_overview}Parameters for the parallel-in-time PFASST-SH scheme.}
\end{table}

To approximate the computational cost of the integration
schemes, we consider a block of $n_{\textit{ts}}$ time steps of
individual size $\Delta t_{\textit{PF}}$ solved in parallel with
PFASST-SH. $C^s_{\ell}$ denotes
the cost of solving the implicit system \ref{implicit_systems} 
on level $\ell$. $C^{\textit{fi}}_{\ell}$ (respectively, $C^{\textit{fe}}_{\ell}$) is the cost of evaluating the implicit (respectively,
explicit) right-hand side on level $\ell$. We decompose the cost
of PFASST($n_{\textit{ts}}$, $M_f+1$, $M_c+1$, $N_{\textit{PF}}$,
$\alpha$) into the cost of the prediction, and the cost of the
PFASST-SH iterations. The cost of the prediction is
\begin{subequations}
\begin{align}
  C^{\textit{pred}} &= (M_c+1) (C^{\textit{fi}}_c + C^{\textit{fe}}_c) \label{cost_initialization_pfasst_-2} \\
                    &+ n_{\textit{ts}} (C^{\textit{fi}}_c + C^{\textit{fe}}_c) \label{cost_initialization_pfasst_-1} \\
                    &+ n_{\textit{ts}} M_c (C^{\textit{s}}_c + C^{\textit{fi}}_c + C^{\textit{fe}}_c),
\label{cost_initialization_pfasst}
\end{align}
\end{subequations}
where the term in \ref{cost_initialization_pfasst_-2} represents
the cost of re-evaluating the coarse right-hand sides at all the
temporal nodes after the restriction in Step \textbf{\textit{B}}.
The term in \ref{cost_initialization_pfasst_-1} corresponds to the
cost of re-evaluating the right-hand side at the first node once
the new initial condition has been received. Finally, the cost
of the $n_{\textit{ts}}$ serial coarse sweeps performed during
the prediction is in \ref{cost_initialization_pfasst}.

The cost of the PFASST iteration detailed in Algorithm
\oldref{alg:pfasst_iteration} reads
\begin{subequations}
\begin{align}
  C^{\textit{iter}} &= M_f (C^{\textit{s}}_f + C^{\textit{fi}}_f + C^{\textit{fe}}_f) \label{cost_iteration_pfasst_-4} \\
                    &+ M_c (C^{\textit{fi}}_c + C^{\textit{fe}}_c) \label{cost_iteration_pfasst_-3} \\
                    &+ C^{\textit{fi}}_c + C^{\textit{fe}}_c \label{cost_iteration_pfasst_-2} \\
                    &+ M_c (C^{\textit{s}}_c + C^{\textit{fi}}_c + C^{\textit{fe}}_c) \label{cost_iteration_pfasst_-1} \\
                    &+ C^{\textit{fi}}_f + C^{\textit{fe}}_f. 
\label{cost_iteration_pfasst}
\end{align}
\end{subequations}
The term in \ref{cost_iteration_pfasst_-4} represents the cost
of the fine sweep (Step \textbf{\textit{A}}). The term in
\ref{cost_iteration_pfasst_-3} is the cost of re-evaluating
the coarse right-hand sides at all nodes except the first one
after the restriction (Step \textbf{\textit{B}}). The term
in \ref{cost_iteration_pfasst_-2} represents the cost of
re-evaluating the coarse right-hand sides at the first node
after receiving the new initial condition in Step
\textbf{\textit{C}}. The cost of the coarse sweep is in
\ref{cost_iteration_pfasst_-1}, and the cost of re-evaluating
the fine right-hand side after receiving the new initial condition
is in \ref{cost_iteration_pfasst}.

Assuming that $N_{\textit{PF}}$ iterations are performed and
denoting the communication costs of the full algorithm --
i.e., prediction and iterations -- by $C^{\textit{comm}}$, we obtain
\begin{equation}
C^{\textit{PFASST}(n_{\textit{ts}}, M_f+1, M_c+1, N_{\textit{PF}}, \alpha)} = C^{\textit{pred}} + N_{\textit{PF}} C^{\textit{iter}} + C^{\textit{comm}}.
\label{total_computational_cost_parallel_pfasst}
\end{equation}
This derivation assumes that the cost of computing the FAS
correction at the end of Step \textbf{\textit{B}} (a linear
combination of function values already accounted for) is negligible.
The computational cost of the serial SDC and MLSDC-SH are considered
in \cite{hamon2018multi}. Briefly, using the same block 
of $n_{\textit{ts}}$ time steps of individual size
$\Delta t_{\textit{PF}} = \Delta t_{\textit{SDC}}$ the serial
SDC($M_f+1$, $N_S$) has a computational cost of
\begin{subequations}
\begin{align}
  C^{\textit{SDC}(M_f+1,N_{\textit{S}})} &= n_{\textit{ts}} M_f ( C^{\textit{fi}}_f + C^{\textit{fe}}_f )   \label{computational_cost_serial_sdc_-1} \\
  &+ n_{\textit{ts}} N_{\textit{S}} M_f ( C^s_{f} + C^{\textit{fi}}_{f} + C^{\textit{fe}}_{f} ),
  \label{computational_cost_serial_sdc}
\end{align}
\end{subequations}
where the term in \ref{computational_cost_serial_sdc_-1} represents
the cost of recomputing the right-hand sides at the beginning of
each time step -- not accounted for in \cite{hamon2018multi}
because it is the same for SDC and MLSDC-SH -- and the term in
\ref{computational_cost_serial_sdc} is the cost of the sweeps. 
To compare \ref{total_computational_cost_parallel_pfasst} and
\ref{computational_cost_serial_sdc}, we express the cost of
the coarse operators as functions of the fine grid quantities
and the spatial coarsening ratio. To this end, we assume that
the cost of the operators is proportional to the number of
spectral coefficients in \ref{spectral_space}, denoted by
$K_{\ell}$. We know that $K_{\ell} = 3 R_{\ell}(R_{\ell}+1)/2$,
where $R_{\ell}$ is the number of zonal wavenumbers
on level $\ell$, since there are three primary variables
individually represented with $R_{\ell}(R_{\ell}+1)/2$ spectral
coefficients. This yields
\begin{equation}
C^s_c = \frac{K_c}{K_f} C^s_f = \frac{R_c(R_c+1)}{R_f(R_f+1)} C^s_f = \alpha^2 \frac{R_f + 1/\alpha}{R_f + 1} C^s_f \approx \alpha^2 C^s_f,
\label{assumption_coarsening_ratio}  
\end{equation}
where we have used the definition of the spatial coarsening ratio
$\alpha = R_c / R_f$ and the fact that
$(R_f+1/\alpha)/(R_f+1) \approx 1$ for the parameters considered
in the numerical examples of Section
\oldref{section_numerical_examples} -- i.e., $R_f \geq 256$ and 
$\alpha \geq 1/5$. We obtain from an analogous derivation that
$C^{fi}_c \approx \alpha^2 C^{fi}_f$ and $C^{fe}_c \approx \alpha^2 C^{fe}_f$.
For simplicity, we further assume that the cost of
evaluating each right-hand side is equal to that of solving
the implicit system -- even though, in practice, evaluating
the nonlinear term requires applying the SH transform and is
therefore more expensive than the other operations. This assumption
reads
\begin{equation}
C^s_{f} = C^{\textit{fi}}_{f} = C^{\textit{fe}}_{f} = 1.
\label{unusual_assumption}
\end{equation}
Then, we use \ref{assumption_coarsening_ratio}
and \ref{unusual_assumption} to eliminate $C^s_{f}$ and $C^s_{c}$
from $C^{\textit{PFASST}(n_{\textit{ts}}, M_f+1, M_c+1, N_{\textit{PF}}, \alpha)}$
and $C^{\textit{SDC}(M_f+1,N_{\textit{S}})}$. After
dividing both $C^{\textit{PFASST}(n_{\textit{ts}}, M_f+1, M_c+1, N_{\textit{PF}}, \alpha)}$
and $C^{\textit{SDC}(M_f+1,N_{\textit{S}})}$
by $3 n_{\textit{ts}} M_f$, we obtain a compact form 
of the theoretical speedup of PFASST-SH with respect to SDC
\begin{equation}
\mathcal{S}^{\textit{theo}}_{\textit{S}}
= \frac{C^{\textit{SDC}(M_f+1,N_{\textit{S}})}}{C^{\textit{PFASST}(n_{\textit{ts}}, M_f+1, M_c+1, N_{\textit{PF}})}}
= \frac{N_S + a}{ b N_{\textit{PF}} + c \alpha^2 N_{\textit{PF}} + d \alpha^2},
\label{speedup_pfasst_sdc}
\end{equation}
where the coefficients are defined as
$a = 2/3$,
$b = (3M_f+2)/(3 n_{\textit{ts}} M_f)$, 
$c = (5M_c+2)/(3 n_{\textit{ts}} M_f)$, and
$d = \big(2M_c+(3M_c+2)n_{\textit{ts}}+2\big)/(3 n_{\textit{ts}} M_f)$.

The above derivation assumes that PFASST-SH and SDC can achieve
the same accuracy for the same time step size, which is only
true if and when the residuals in the PFASST-SH iterations
converge to the same level as as the SDC scheme. Whenever this
is not the case, we rescale the expression of
\ref{speedup_pfasst_sdc} by the ratio
$\Delta t_{\textit{PF}} / \Delta t_{\textit{SDC}}$.
The theoretical speedup of PFASST-SH with respect to MLSDC-SH, whose
derivation is omitted, reads
\begin{equation}
\mathcal{S}^{\textit{theo}}_{\textit{ML}} = \frac{N_{\textit{ML}} + a}{ b N_{\textit{PF}} + c \alpha^2 N_{\textit{PF}} + d \alpha^2} 
                                      +  \frac{N_{\textit{ML}}}{ e N_{\textit{PF}}/\alpha^2 + f N_{\textit{PF}} + g}.
\label{speedup_pfasst_mlsdc}
\end{equation}
In \ref{speedup_pfasst_mlsdc}, we have introduced
$e = (3M_f+2)/(5n_{\textit{ts}}M_c)$,
$f = (5M_c+2)/(5n_{\textit{ts}}M_c)$, 
and
$g = \big(2M_c+(3M_c+2)n_{\textit{ts}}+2\big)/(5n_{\textit{ts}}M_c)$.
We will compare the theoretical speedups obtained in this section
with the observed speedups based on wall-clock time in the next
section.


\section{\label{section_numerical_examples}Numerical examples}

We are ready to analyze the accuracy, stability, and
computational cost of the proposed parallel-in-time integration method
to determine the conditions in which PFASST-SH reduces the
time-to-solution. We use numerical examples
of increasing complexity to successively evaluate the sensitivity
of these PFASST-SH properties to the time step size, the number
of iterations, and the number of processors.
The examples are performed using the finest spatial resolution
given by $R_f = 256$ in \ref{spectral_to_physical_space}, which
is sufficient to generate complex flow dynamics over the time
integration window. The error is computed with the following norm
on the spectral coefficients
\begin{equation} \label{error_norm_r_norm}
  E_{R_{\textit{norm}}} = \frac{
    || \phi - \phi_{\textit{ref}} ||_{\infty, \, R_{\textit{norm}}}
  }{
    || \phi_{\textit{ref}} ||_{\infty, \, R_{\textit{norm}}}
  } 
  \qquad
  \text{with}
  \qquad
  || \phi ||_{\infty, \, R_{\textit{norm}}} =   \max_{\substack{r \in \{ 0, \dots, R_{\textit{norm}} \} \\ s \in \{ r, \dots, R_{\textit{norm}} \} }} | \phi^r_s |,
\end{equation}
where $\phi$ is the approximate solution and $\phi_{\textit{ref}}$
is the reference solution. We note that in
\ref{error_norm_r_norm}, the norm of the difference between the
computed solution and the reference solution is normalized by the
norm of the reference solution. If we set $R_{\textit{norm}} = R_f$
in \ref{error_norm_r_norm}, $|| \cdot ||_{\infty, \, R_f}$ is the
standard $L_{\infty}$-norm in spectral space that measures the
ability of PFASST-SH to capture the fully detailed solution.
If $R_{\textit{norm}} < R_f$, $|| \cdot ||_{\infty, \, R_{\textit{norm}}}$ 
is a semi-norm that assesses the accuracy of PFASST-SH for the
large-scale features represented by the slow modes of the spectral
basis.  As mentioned in  Section~\oldref{section_introduction},
the use of damping methods to resolve the spectral
blocking problem -- here, second-order artificial diffusion --
implies that the accuracy on the
modes corresponding to high-frequency wavenumbers
is not of primary interest.  Hence, in the numerical examples,
we consider
$R_{\textit{norm}} = 32$, $R_{\textit{norm}} = 64$, and
$R_{\textit{norm}} = R_f = 256$.  

The studies presented next have been conducted on the
Cori supercomputer at Lawrence Berkeley National Laboratory.
Even though the following test cases are run without any
parallelization in space for simplicity, our parallel-in-time
method is meant to allow access to an additional axis of
parallelism whenever spatial parallelization is saturated. 

\subsection{\label{subsection_propagation_of_a_gaussian_dome}Propagation of a Gaussian dome}

\begin{figure}[ht!]
\centering
\subfigure[]{
  \begin{tikzpicture}
    \node[anchor=south west,inner sep=0] at (0,0){\includegraphics[scale=0.27]{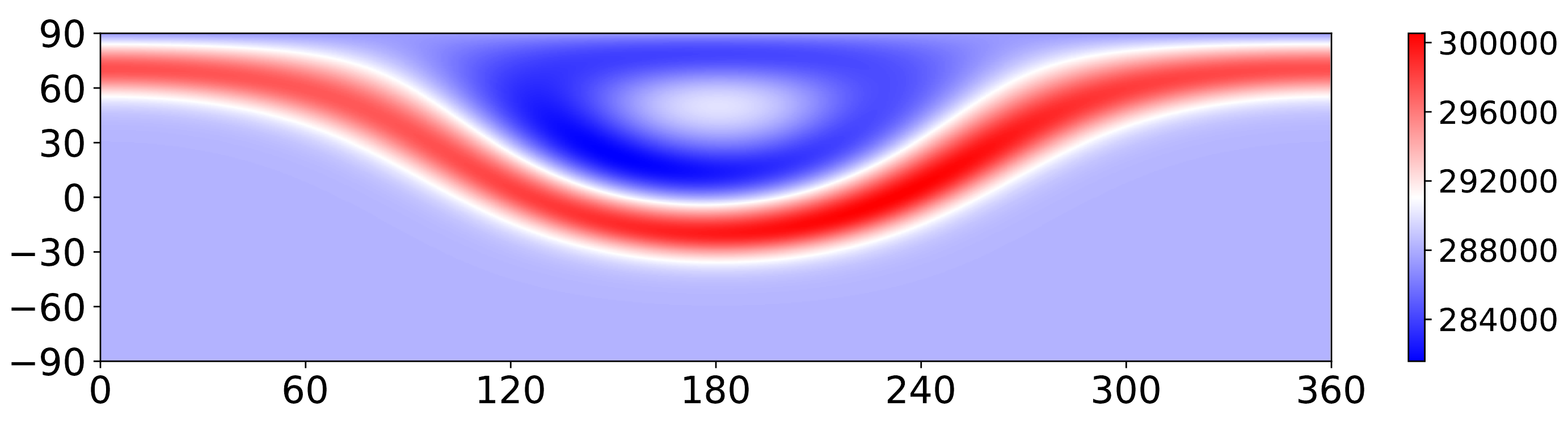}};
    \node (ib_1) at (3.4,-0.05) {\scriptsize Longitude (degrees)};
    \node[rotate=90] (ib_1) at (-0.12,1.2) {\scriptsize Latitude (degrees)};
\end{tikzpicture}
\label{fig:gaussian_bump_test_case_geopotential_field_1}
}
\hspace{-0.5cm}
\subfigure[]{
  \begin{tikzpicture}
    \node[anchor=south west,inner sep=0] at (0,0){\includegraphics[scale=0.27]{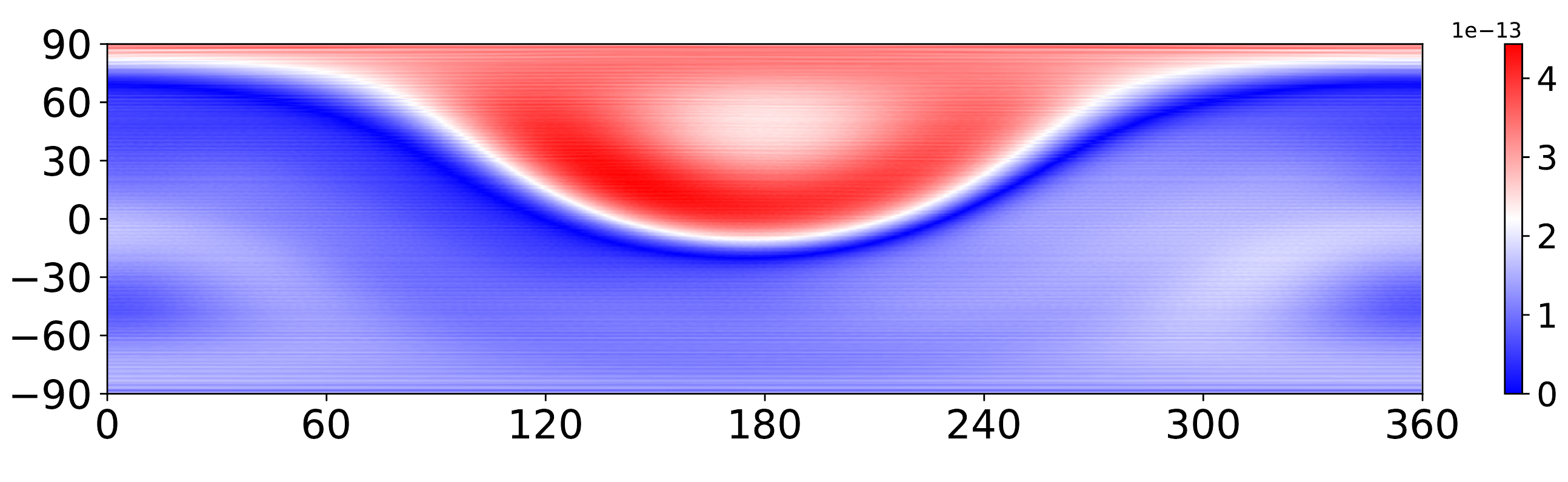}};
    \node (ib_1) at (3.4,-0.05) {\scriptsize Longitude (degrees)};
    \path (6.5,2.015) node (d) {};
    \path (6.9,2.115) node (e) {};
    \path [draw=white,fill=white] (d) rectangle (e); 
    \path (7.425,1.84) node (e) {\scalebox{0.55}{$\times 10^{-13}$}};
\end{tikzpicture}
\label{fig:gaussian_bump_test_case_geopotential_field_2}
}
\vspace{-0.5cm}
\caption{\label{fig:gaussian_bump_test_case_geopotential_field} 
  Gaussian bump: geopotential field after $2 \, \text{hours}$, $40 \, \text{min}$
  obtained using PFASST(16,5,3,8,4/5) with $R_f = 256$ and $\Delta t = 100 \, \text{s}$
  in \oldref{fig:gaussian_bump_test_case_geopotential_field_1}. Figure
  \oldref{fig:gaussian_bump_test_case_geopotential_field_2} shows
  the normalized difference (in the physical coefficients)
  between the PFASST-SH solution and the reference SDC(5,8)
  solution with the same
  spatial resolution and $\Delta t_{\textit{ref}} = 60 \, \text{s}$.
}
\end{figure}

\begin{figure}[ht!]
  \centering
  \subfigure[]{
    \begin{tikzpicture}
    \node[anchor=south west,inner sep=0] at (0,0.){\includegraphics[scale=0.385]{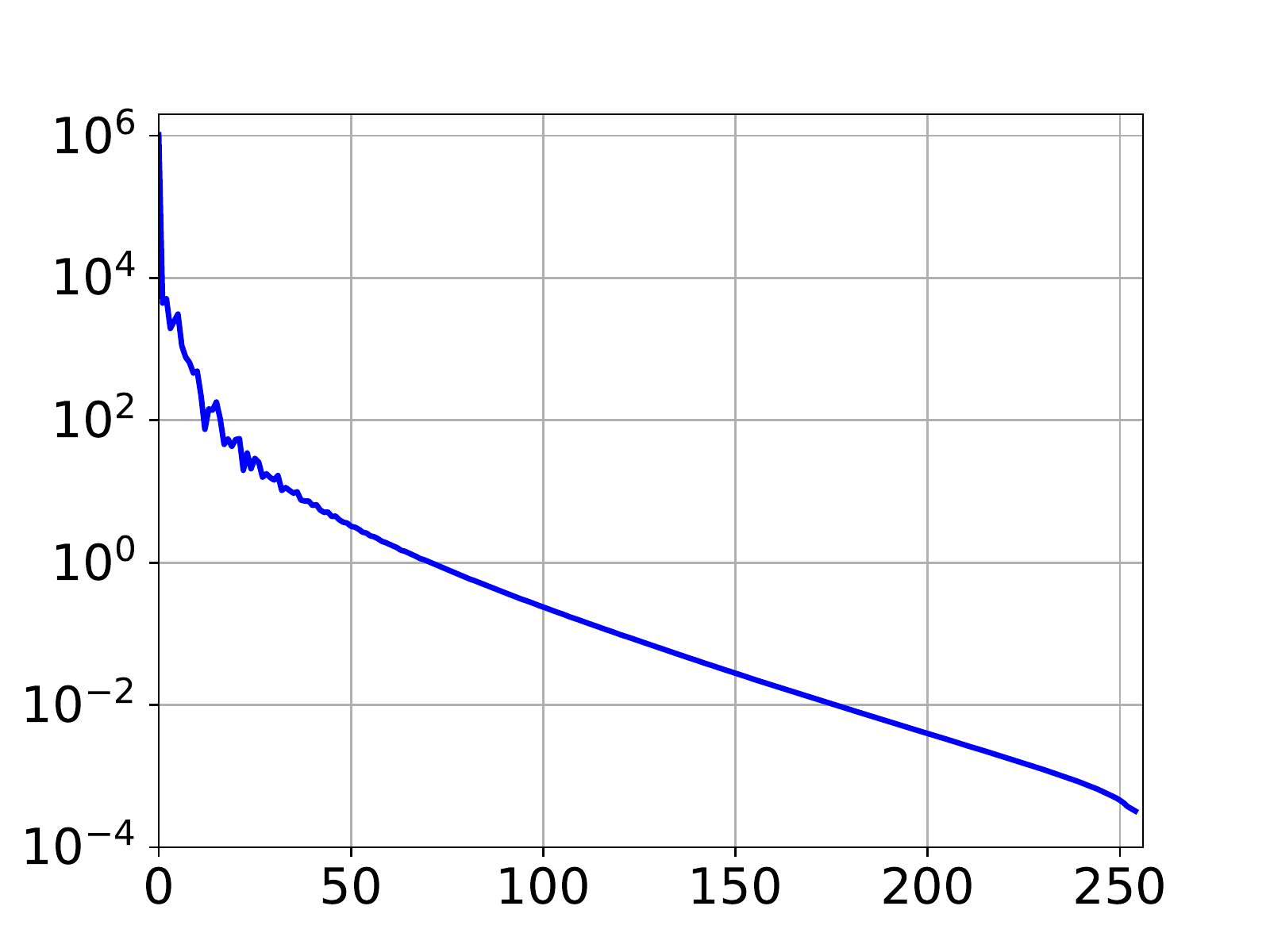}};
    \node (ib_1) at (3.1,-0.05) {$n_0$};
    \node[rotate=90] (ib_1) at (-0.05,2.15) {$|\Phi_{n_0}|$};

    
    \node (ib_1) at (1.75,4.25) {};
    \node (ib_2) at (1.75,0.4) {};
    \path [draw=black!100, thick, dashed] (ib_1) -- (ib_2);

    \node[rotate=-90] (ib_3) at (1.95,2.87) {\scriptsize $\alpha = 1/5$}; 

    \node (ib_1) at (3.21,4.25) {};
    \node (ib_2) at (3.21,0.4) {};
    \path [draw=black!100, thick, dashed] (ib_1) -- (ib_2);

    \node[rotate=-90] (ib_3) at (3.41,2.32) {\scriptsize $\alpha = 1/2$}; 

    \node (ib_1) at (4.68,4.25) {};
    \node (ib_2) at (4.68,0.4) {};
    \path [draw=black!100, thick, dashed] (ib_1) -- (ib_2);

    \node[rotate=-90] (ib_3) at (4.88,1.87) {\scriptsize $\alpha = 4/5$};

    \node (ib_3) at (3.2,4.35) {Geopotential}; 

    \end{tikzpicture}
    \label{fig:gaussian_bump_test_case_geopotential_spectrum}
}
  \subfigure[]{
    \begin{tikzpicture}
    \node[anchor=south west,inner sep=0] at (0,0){\includegraphics[scale=0.385]{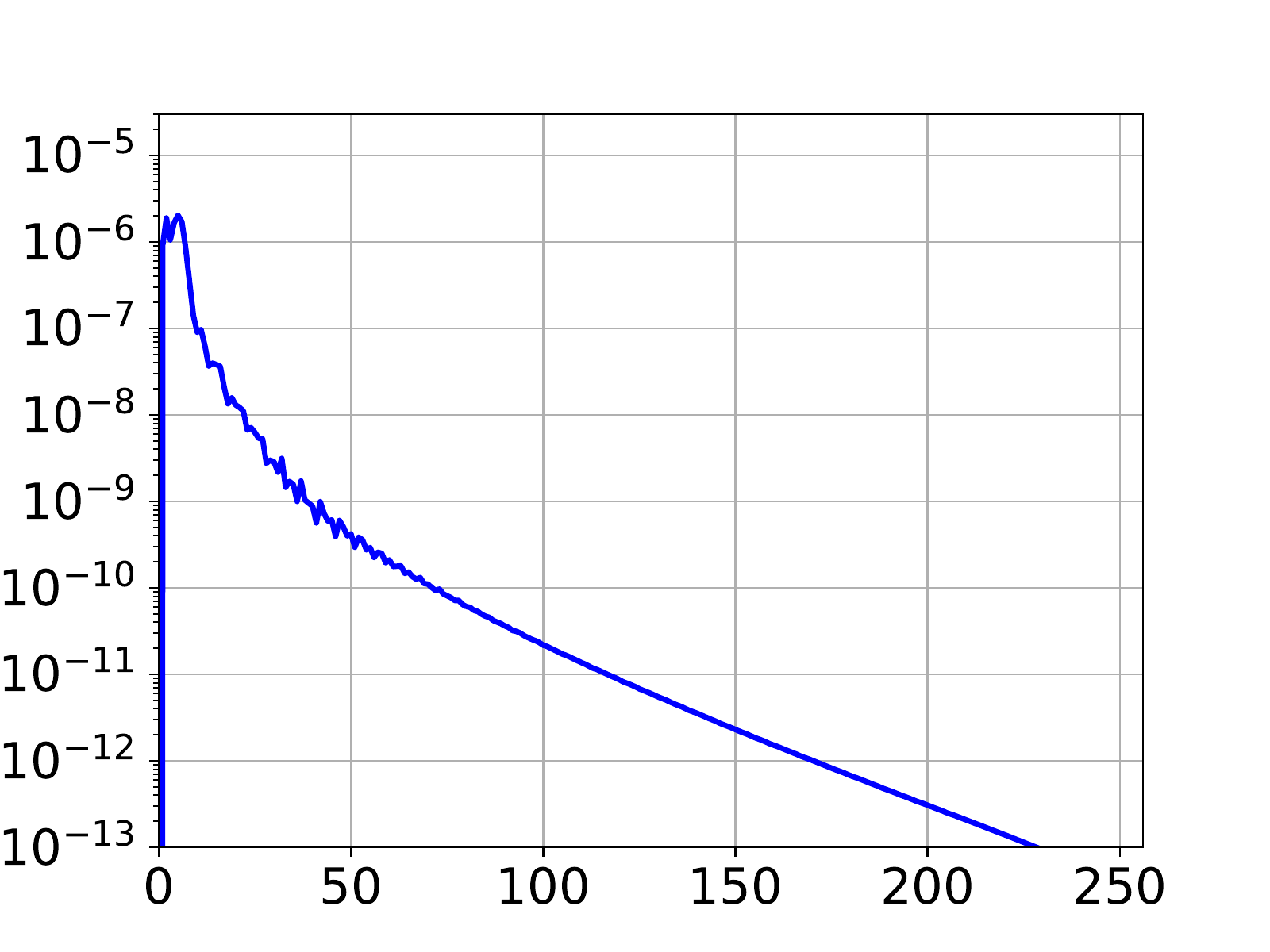}};
    \node (ib_1) at (3.1,-0.05) {$n_0$};
    \node[rotate=90] (ib_1) at (-0.25,2.15) {$|\zeta_{n_0}|$};

    
    \node (ib_1) at (1.75,4.25) {};
    \node (ib_2) at (1.75,0.4) {};
    \path [draw=black!100, thick, dashed] (ib_1) -- (ib_2);

    \node[rotate=-90] (ib_3) at (1.95,2.87) {\scriptsize $\alpha = 1/5$}; 

    \node (ib_1) at (3.21,4.25) {};
    \node (ib_2) at (3.21,0.4) {};
    \path [draw=black!100, thick, dashed] (ib_1) -- (ib_2);

    \node[rotate=-90] (ib_3) at (3.41,2.32) {\scriptsize $\alpha = 1/2$}; 

    \node (ib_1) at (4.68,4.25) {};
    \node (ib_2) at (4.68,0.4) {};
    \path [draw=black!100, thick, dashed] (ib_1) -- (ib_2);

    \node[rotate=-90] (ib_3) at (4.88,1.87) {\scriptsize $\alpha = 4/5$};

    \node (ib_3) at (3.2,4.35) {Vorticity}; 
    \end{tikzpicture}
    \label{fig:gaussian_bump_test_case_vorticity_spectrum}
}

\vspace{-0.4cm}
\caption{\label{fig:gaussian_bump_test_case_spectrum} 
  Gaussian bump: max-spectrum of the geopotential field in
  \oldref{fig:gaussian_bump_test_case_geopotential_spectrum} and of
  the vorticity field in \oldref{fig:gaussian_bump_test_case_vorticity_spectrum}
  obtained after $28 \, \text{hours}$, $26 \, \text{min}$, $40 \, \text{s}$.
  These figures are obtained using SDC(5,8) with $R_f = 256$
  and $\Delta t_{\textit{ref}} = 60 \, \text{s}$. The quantities on the $y$-axis
  are defined as $|\Phi_{n_0}| = \max_{r} |\Phi^r_{n_0}|$ and
  $|\zeta_{n_0}| = \max_{r} |\zeta^r_{n_0}|$. We observe rapidly decaying
  modes across the full spectrum. For each coarsening ratio $\alpha$, a
  vertical dashed line shows the fraction of the spectrum that is truncated during
  the coarsening step in MLSDC-SH and PFASST-SH. 
}
\end{figure}

We first consider  the propagation of a Gaussian bump on the sphere. This example is adapted
from \cite{swarztrauber2004shallow} and was used in \cite{hamon2018multi} to assess
the performance of MLSDC-SH. The initial velocities are set to zero ($u = v = 0$),
and a Gaussian dome centered at $\lambda_c = \pi$ and $\phi_c = \pi/4$ is introduced
in the geopotential field,
such that the height is given by
\begin{equation}
h( \lambda, \phi ) = \bar{h} + A \text{e}^{-\alpha (d / a )^2}, \label{gaussian_bump}
\end{equation}
where $a$ denotes the earth radius. The distance $d$ is defined as
\begin{equation}
d = \sqrt{x^2 + y^2 + z^2}, 
\end{equation}
with
\begin{align}
x &= a \big( \cos(\lambda) \cos(\phi) - \cos(\lambda_c) \cos(\phi_c) \big), \\
y &= a \big( \sin(\lambda) \cos(\phi) - \sin(\lambda_c) \cos(\phi_c) \big), \\
z &= a \big( \sin(\phi) - \sin(\phi_c) \big).
\end{align} 
We use the same values as in \cite{swarztrauber2004shallow} for the Earth radius,
the gravitational acceleration, and the angular rate of rotation $\Omega$
involved in the Coriolis force. In \ref{gaussian_bump}, we set
$\bar{h} = 29,400 \, \text{m}$ and $A = 6,000 \, \text{m}$. We simulate
$102,400 \, \text{s}$ ($28 \, \text{hours}$, $26 \, \text{min}$, $40 \, \text{s}$)
of propagation with a spatial resolution $R_f = 256$ and a diffusion coefficient
$\nu = 10^5 \, \text{m}^2.\text{s}^{-1}$. PFASST-SH is run with 16 processors, using
one time step per compute node. The geopotential field obtained using
PFASST($n_{\textit{ts}} = 16$, $M_f + 1= 5$, $M_c+1 = 3$, $N_{\textit{PF}} =8$, $\alpha = 4/5$)
-- i.e., with five SDC nodes on the fine level, three SDC nodes on the coarse level,
eight iterations, and a spatial coarsening ratio $\alpha = 4/5$ -- with a time step
size $\Delta t = 100 \, \text{s}$ is  shown in Fig.~\oldref{fig:gaussian_bump_test_case_geopotential_field}.
This figure also shows the difference between this PFASST-SH solution and the reference
solution obtained using SDC($M_f + 1 = 5$, $N_{S} = 8$) based on the same spatial
resolution and $\Delta t_{\textit{ref}} = 60 \, \text{s}$, illustrating the fact
that PFASST-SH can accurately reproduce the reference geopotential map without numerical
artifacts. The spectrum of the reference solution is in
Fig.~\oldref{fig:gaussian_bump_test_case_spectrum}.

\begin{figure}[ht!]
\centering
\subfigure[]{
\begin{tikzpicture}
\node[anchor=south west,inner sep=0] at (0,0.045){\includegraphics[scale=0.385]{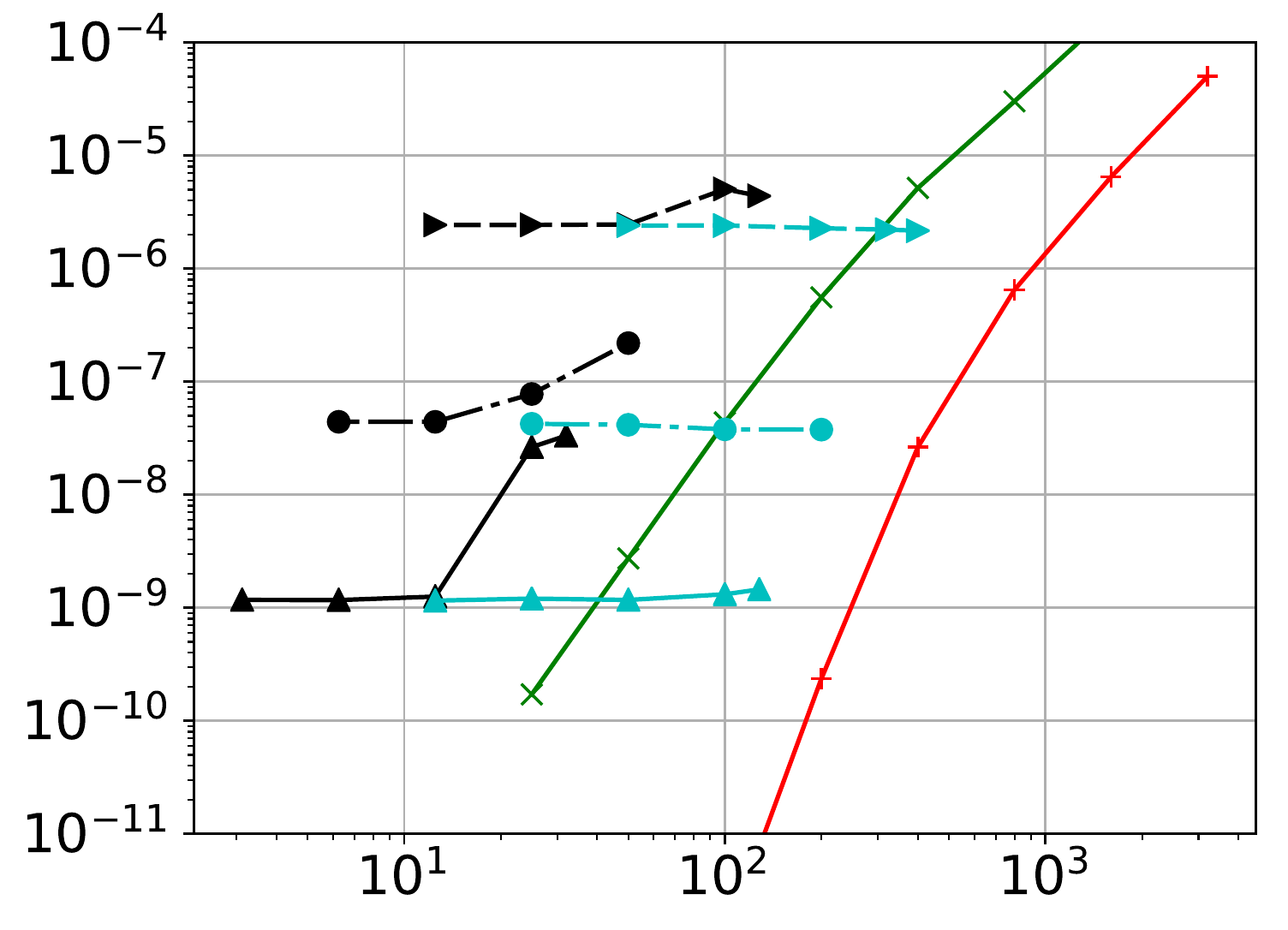}};
\node (ib_1) at (3.2,-0.05) {$\Delta t$ $(s)$};
\node[rotate=90] (ib_1) at (-0.1,2.25) {Normalized $L_{\infty}$-error};
\node (ib_1) at (3.3,4.3) {$R_{\textit{norm}} = 256$};
\node (ib_1) at (2.4,3.45) {\scriptsize $\alpha = 1/5$};
\node (ib_1) at (1.65,2.55) {\scriptsize $\alpha = 1/2$};
\node (ib_1) at (1.5,1.7) {\scriptsize $\alpha = 4/5$};
\end{tikzpicture}
\label{fig:gaussian_bump_test_case_accuracy_geopotential_256}
}
\hspace{-0.5cm}
\subfigure[]{
\begin{tikzpicture}
\node[anchor=south west,inner sep=0] at (0,0.05){\includegraphics[scale=0.385]{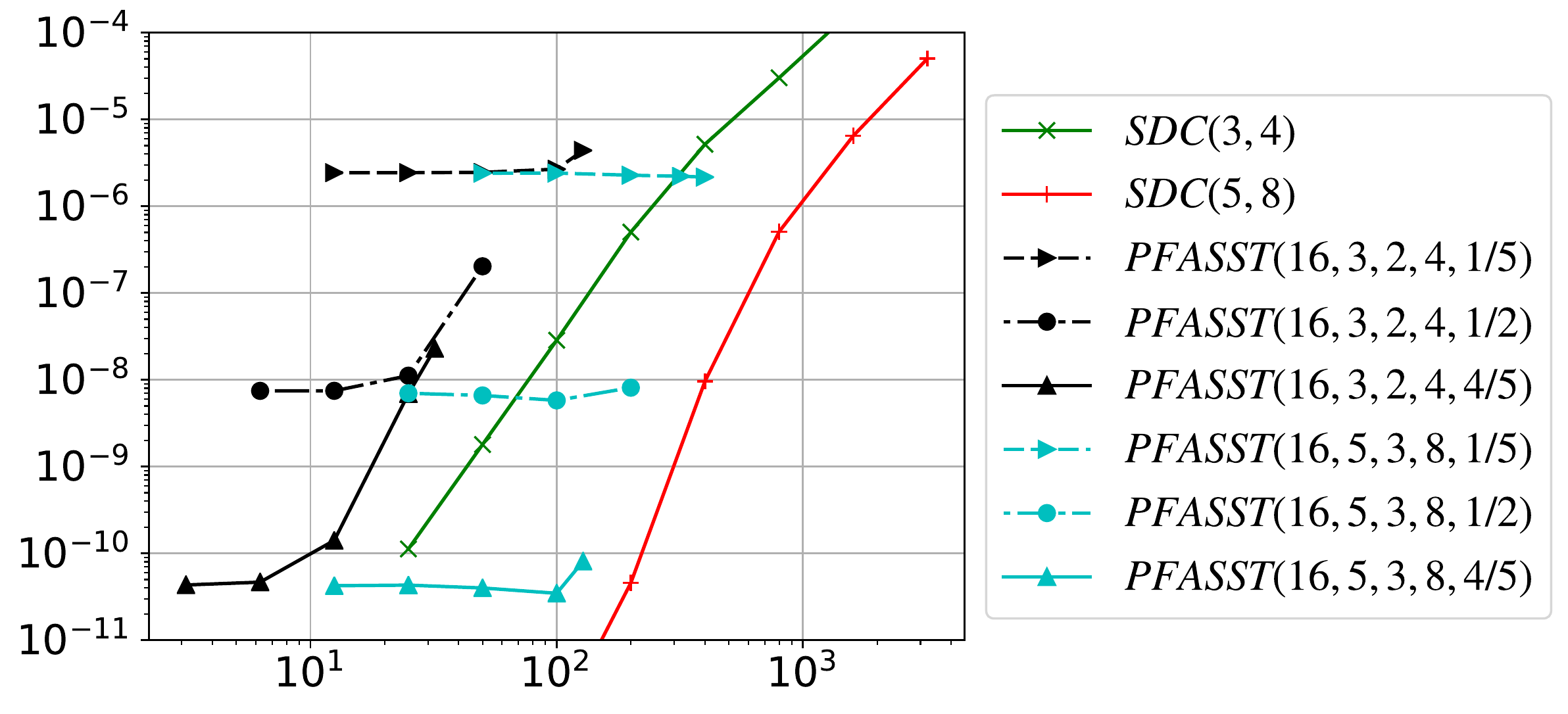}};
\node (ib_1) at (3.2,-0.05) {$\Delta t$ $(s)$};
\node (ib_1) at (3.1,4.325) {$R_{\textit{norm}} = 32$};
\end{tikzpicture}
\label{fig:gaussian_bump_test_case_accuracy_geopotential_32}
}
\vspace{-0.45cm}
\caption{\label{fig:gaussian_bump_test_case_accuracy_geopotential_field} 
  Gaussian bump: $L_{\infty}$-norm of the error in the geopotential field after
  $28 \, \text{hours}$, $26 \, \text{min}$, $40 \, \text{s}$
  as a function of time step size. We apply the same norm to the SDC and PFASST-SH
  schemes, with $R_{\text{norm}} = 256$ in \oldref{fig:gaussian_bump_test_case_accuracy_geopotential_256},
  and $R_{\text{norm}} = 32$
  in \oldref{fig:gaussian_bump_test_case_accuracy_geopotential_32}. In the legend, the fourth (respectively,
  fifth) parameter of the PFASST-SH schemes denotes the number of iterations (respectively, the coarsening
  ratio in space). For each scheme,
  the rightmost point corresponds to the largest stable time step.
}
\end{figure}

To gain an understanding of the accuracy and stability of PFASST-SH,
we perform a study of the error norm as a function of time step
size for a fixed number of iterations. The geopotential error with
respect to the reference solution is shown as a function of time
step size in
Fig.~\oldref{fig:gaussian_bump_test_case_accuracy_geopotential_field}.
In this figure and in the numerical results shown in the remainder
of this paper, we consider a wide range of normalized errors
-- from $10^{-4}$ to $10^{-11}$ -- to compare PFASST-SH and SDC
in a regime where the SDC schemes achieve their formal order of
accuracy -- up to eighth order for SDC(5,8).
We focus on PFASST($n_{\textit{ts}} = 16$,
$M_f + 1 = 3$, $M_c + 1 = 2$, $N_{\textit{PF}} = 4$, $\alpha$) and
PFASST($n_{\textit{ts}} = 16$, $M_f + 1 = 5$, $M_c + 1 = 3$,
$N_{\textit{PF}} = 8$, $\alpha$) with four and eight iterations,
respectively. We explain at the end of this
section how this fixed number of iterations, $N_{\textit{PF}}$,
is determined. We consider the spatial coarsening ratios
$\alpha = 1/5, \, 1/2, \, \text{and} \, 4/5$.
The range of error magnitudes observed here is
consistent with the results in the numerical examples of
\cite{jia2013spectral}.

Figure~\oldref{fig:gaussian_bump_test_case_accuracy_geopotential_field} 
illustrates two key properties of PFASST-SH. First, considering
stability, we observe that the PFASST-SH schemes can take
relatively large time steps, on the order of $10^2 \, \text{s}$,
but we also note that their largest time step size decreases when
the spatial coarsening ratio, $\alpha$, is increased.  The
largest time steps taken by PFASST-SH are about ten times smaller
than with the corresponding serial SDC schemes. Second, we see that
the PFASST-SH accuracy heavily depends on the choice of a spatial
coarsening ratio, as this parameter determines the portion of the
spectrum of the solution that can be resolved on the coarse level
(see Fig.~\oldref{fig:gaussian_bump_test_case_spectrum}).

This is because for an aggressive spatial coarsening -- i.e.,
$\alpha$ close to zero --, the main features of the solution
cannot be captured on the coarse level and as a result, 
the PFASST-SH error stagnates at a relatively large magnitude. 
Conversely, for a milder spatial coarsening -- with $\alpha$
close to one --, higher-frequency modes are included on the coarse
level. This significantly improves the accuracy of PFASST-SH, but
also imposes a more severe stability limit on the time
step size.

Finally, for $\alpha = 1/2$ and $\alpha = 4/5$, we observe that
the PFASST-SH accuracy better matches that of the SDC schemes 
when we only consider the low-frequency modes in the computation
of the error ($R_{\textit{norm}} = 32$). This means that
PFASST-SH can efficiently 
resolve the spectral coefficients corresponding to the
large-scale features of the solution, at the expense of the accuracy
for the small-scale features. 

\begin{table}[!ht]
\centering
 \scalebox {0.9}{
         \begin{tabular}{lccccc}
           \\ \toprule 
           Integration               &        $L_{\infty}$-error        & Time step & Wall-clock & Observed & Observed  \\
           scheme                    &        $R_{\textit{norm}} = 32$  & size (s)  & time (s)   & speedup $S^{\textit{obs}}_{S}$ & speedup $S^{\textit{obs}}_{\textit{ML}}$ \\\toprule
           SDC(3,4)                  &       $5.2 \times 10^{-6}$       &        400             &        528      & -     &  -          \\
           MLSDC(3,2,2,1/2)          &       $5.2 \times 10^{-6}$       &        400             &        354      & 1.5    &  -          \\
           PFASST(16,3,2,4,1/5)      &       $4.4 \times 10^{-6}$       &        128             &        120      & 4.1    &  3.0        \\ \midrule
           SDC(5,8)                  &       $4.6 \times 10^{-11}$      &        200             &        3,389   & -     &  -          \\
           MLSDC(5,3,4,1/2)          &       $7.6 \times 10^{-10}$      &       200              &        2,269   & 1.5   &  -          \\
           PFASST(16,5,3,8,4/5)      &       $3.5 \times 10^{-11}$      &       100              &        917     & 3.7   &  2.5        \\
           \bottomrule 
         \end{tabular}}
\caption{\label{tab:gaussian_bump_timings}
  Gaussian bump: $L_{\infty}$-error for $R_{\textit{norm}} = 32$, time step size,
  and wall-clock time for serial and parallel SDC-based schemes.
  In the last two columns, $S^{\textit{obs}}_S$ denotes the observed speedup
  achieved with MLSDC-SH and PFASST-SH with respect to SDC, while $S^{\textit{obs}}_{\textit{ML}}$
  denotes the speedup achieved with PFASST-SH with respect  to MLSDC-SH.
  The speedups for $R_{\textit{norm}} = 256$ can be read from
  Fig.~\oldref{fig:gaussian_bump_test_case_computational_cost_geopotential_256}.
}
\end{table}

\begin{figure}[ht!]
\centering
\subfigure[]{
\begin{tikzpicture}
\node[anchor=south west,inner sep=0] at (0,0.045){\includegraphics[scale=0.385]{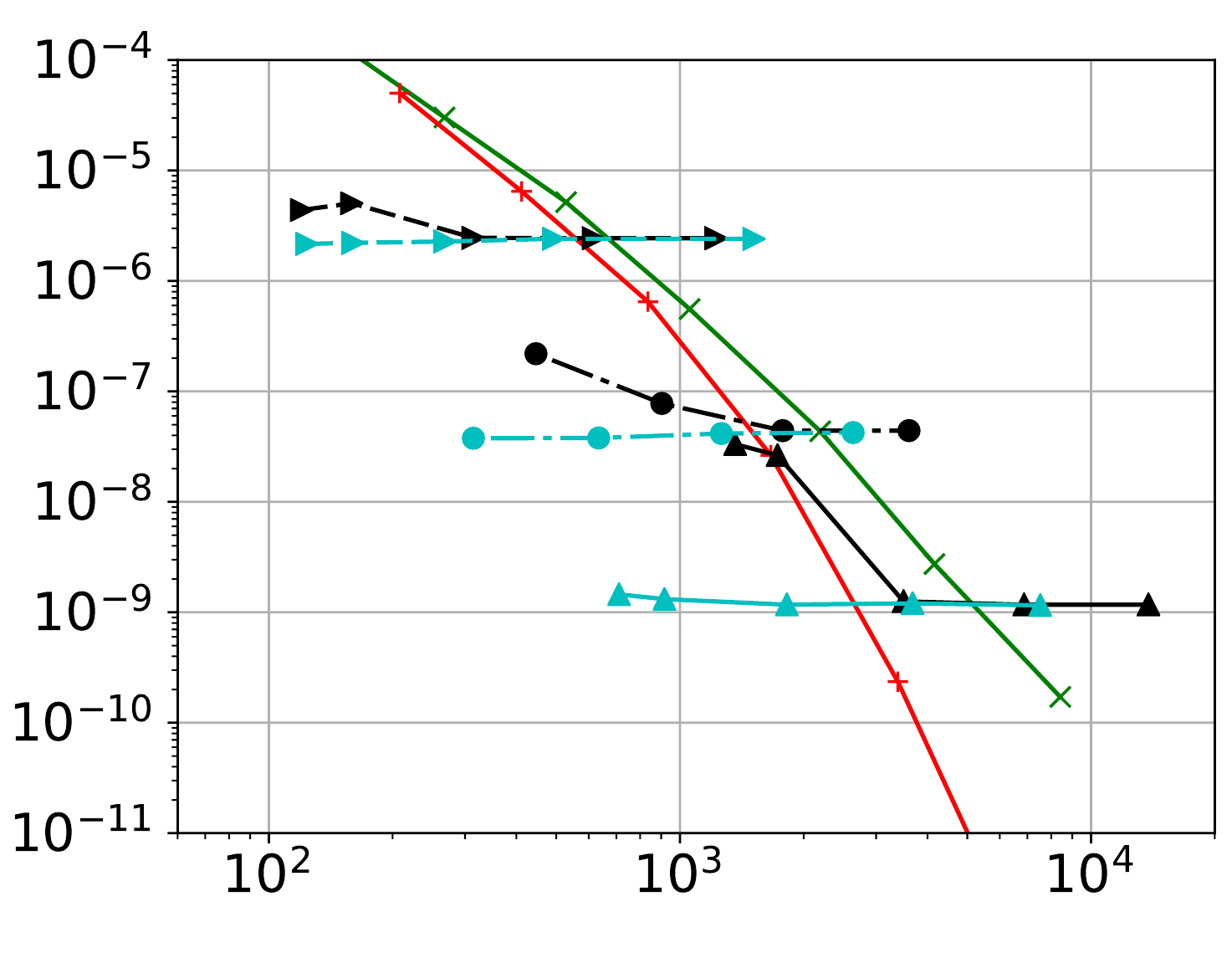}};
\node (ib_1) at (3.2,-0.05) {Wall-clock time $(s)$};
\node[rotate=90] (ib_1) at (-0.1,2.35) {Normalized $L_{\infty}$-error};
\node (ib_1) at (3.3,4.45) {$R_{\textit{norm}} = 256$};
\node (ib_1) at (3.65,3.56) {\scriptsize $\alpha = 1/5$};
\node (ib_1) at (4.55,2.65) {\scriptsize $\alpha = 1/2$};
\node (ib_1) at (5.2,1.825) {\scriptsize $\alpha = 4/5$};

\end{tikzpicture}
\label{fig:gaussian_bump_test_case_computational_cost_geopotential_256}
}
\hspace{-0.5cm}
\subfigure[]{
\begin{tikzpicture}
\node[anchor=south west,inner sep=0] at (0,0.05){\includegraphics[scale=0.385]{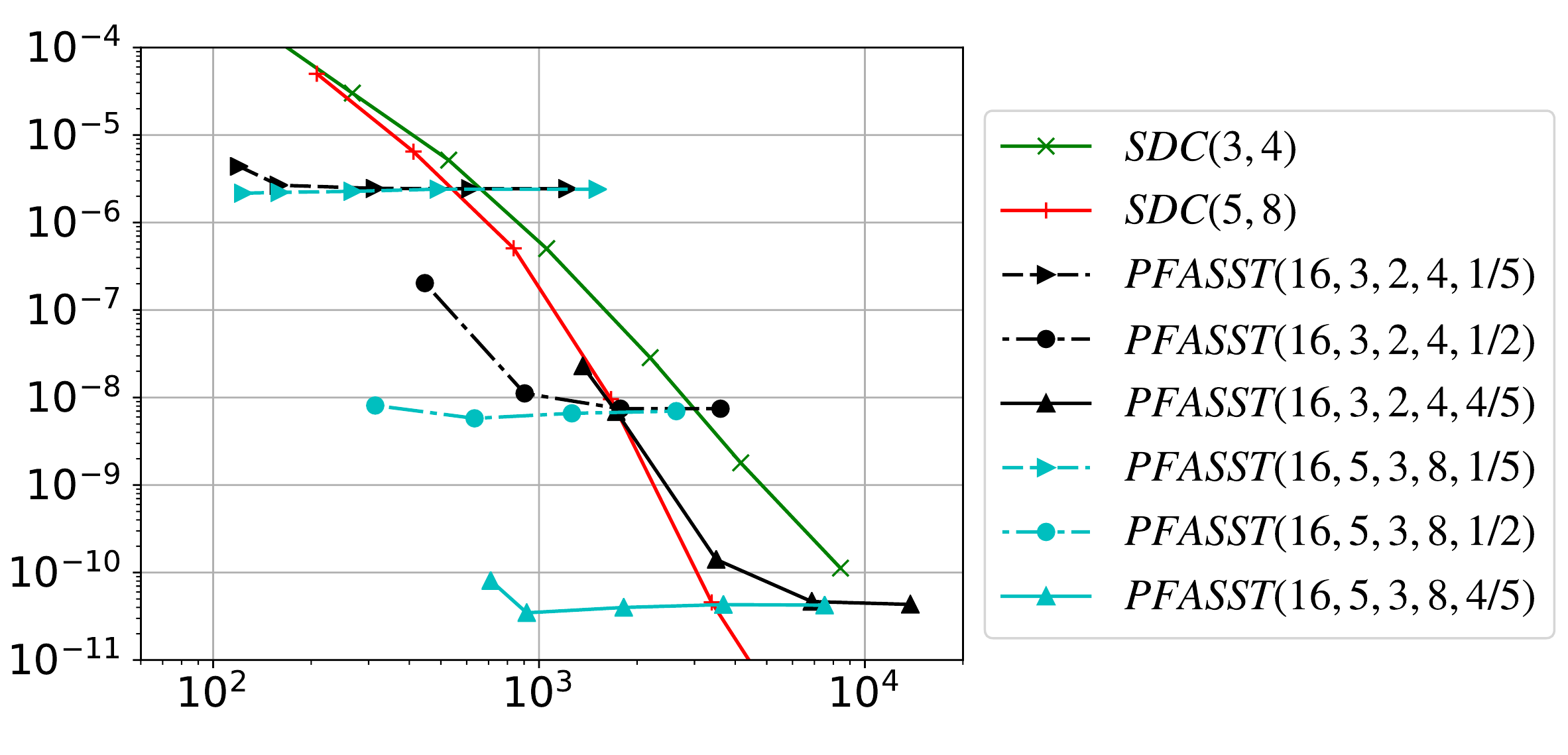}};
\node (ib_1) at (3.2,-0.05) {Wall-clock time $(s)$};
\node (ib_1) at (3.1,4.45) {$R_{\textit{norm}} = 32$};
\end{tikzpicture}
\label{fig:gaussian_bump_test_case_computational_cost_geopotential_32}
}
\vspace{-0.45cm}
\caption{\label{fig:gaussian_bump_test_case_computational_cost_geopotential_field} 
Gaussian bump: $L_{\infty}$-norm of the error in the geopotential field after
$28 \, \text{hours}$, $26 \, \text{min}$, $40 \, \text{s}$
as a function of wall-clock time. We apply the same norm to the SDC and PFASST-SH
schemes, with $R_{\text{norm}} = 256$ in
\oldref{fig:gaussian_bump_test_case_computational_cost_geopotential_256}, and $R_{\text{norm}} = 32$
in \oldref{fig:gaussian_bump_test_case_computational_cost_geopotential_32}.
In the legend, the fourth (respectively,
  fifth) parameter of the PFASST-SH schemes denotes the number of iterations (respectively, the coarsening
  ratio in space).
  For each scheme, the leftmost point corresponds to the largest stable time step. For their largest
  stable time step, we observe that the PFASST-SH schemes
  reduce the computational cost compared to the serial SDC schemes. 
}
\end{figure}

We now evaluate the computational cost of the parallel-in-time
scheme. In 
Fig.~\oldref{fig:gaussian_bump_test_case_computational_cost_geopotential_field},
the geopotential error norm is shown as a function of wall-clock
time for the same PFASST-SH schemes as in
Fig.~\oldref{fig:gaussian_bump_test_case_accuracy_geopotential_field}.
We observe in Table~\oldref{tab:gaussian_bump_timings} that for
their largest time step size, these PFASST-SH
schemes are consistently faster than the corresponding serial
SDC schemes for the range of spatial coarsening ratios considered
here.

To achieve a normalized
error norm of about $10^{-6}$ in the geopotential variable, 
PFASST($n_{\textit{ts}} = 16$, $M_f + 1 = 3$, $M_c + 1 = 2$,
$N_{\textit{PF}} = 4$, $\alpha = 1/5$) has to take a time step size
smaller or equal to $\Delta t =  128 \, \text{s}$, whereas
SDC($n_{\textit{ts}} = 3$, $N_S = 4$) and MLSDC($M_f +1 =3$, $M_c + 1=2$, $N_{\textit{ML}} = 2$,
$\alpha = 1/2$) can achieve the same accuracy with a time step
size $\Delta t = 400 \, \text{s}$. Still, for this relatively low
accuracy, PFASST-SH yields a speedup of $S^{\textit{obs}}_S = 4.1$
(respectively, $S^{\textit{obs}}_{\textit{ML}} = 3.0$)
compared to SDC($M_f + 1 = 3$, $N_S = 4$)
(respectively, MLSDC($M_f + 1 = 3$, $M_c + 1 = 2$,
$N_{\textit{ML}} = 2$, $\alpha = 1/2$)). 
These observed speedups are consistent with the values derived
in Section~\oldref{subsection_computational_cost}.
With a milder coarsening in space, the accuracy of the
parallel-in-time algorithm improves, but the speedup is slightly
reduced. Specifically, for a normalized error norm of
about $10^{-9}$, PFASST($n_{\textit{ts}} = 16$,
$M_f + 1 = 5$, $M_c + 1 = 3$, $N_{\textit{PF}} = 8$, $\alpha = 4/5$),
achieves a speedup of $S^{\textit{obs}}_S = 3.7$
(respectively, $S^{\textit{obs}}_{\textit{ML}} = 2.5$)
compared to SDC($M_f + 1 = 5$, $N_S = 8$) (respectively,
MLSDC($M_f + 1 = 5$, $M_c + 1 = 3$, $N_{\textit{ML}} = 4$,
$\alpha = 1/2$)). 
We note that the MLSDC(5,3,4,1/2) error is about 20 times
larger than the SDC(5,8) error for the same time step size
since the four coarse sweeps do not resolve well the high-frequency
modes of the solution. This is documented in the discussion
of MLSDC-SH found in \cite{hamon2018multi}.

\begin{figure}[ht!]
\centering
\subfigure[]{
\begin{tikzpicture}
\node[anchor=south west,inner sep=0] at (0,0.045){\includegraphics[scale=0.385]{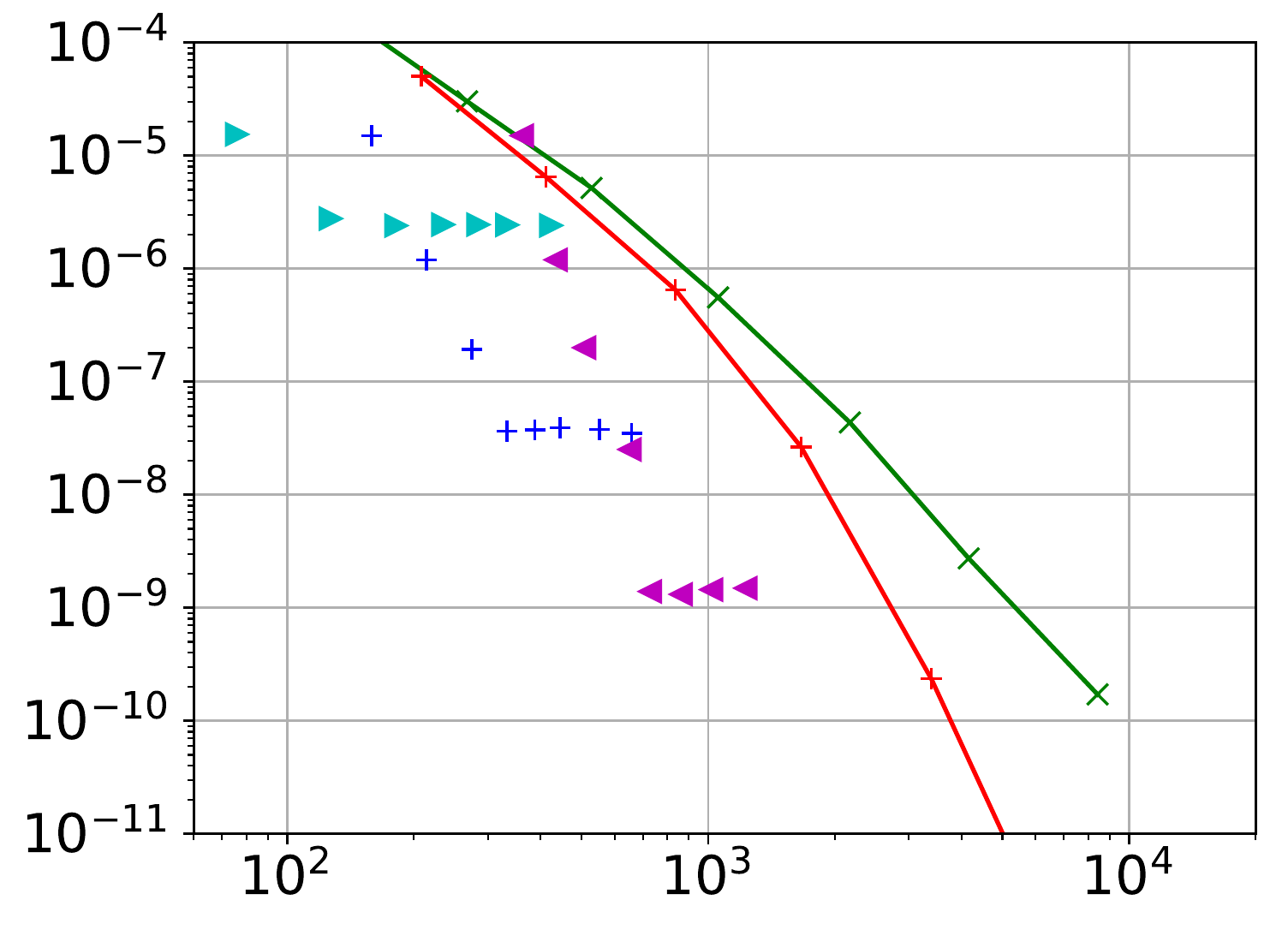}};
\node (ib_1) at (3.2,-0.05) {Wall-clock time $(s)$};
\node[rotate=90] (ib_1) at (-0.1,2.35) {Normalized $L_{\infty}$-error};
\node (ib_1) at (3.3,4.325) {$R_{\textit{norm}} = 256$};
\node (ib_1) at (1.075,3.825) {\tiny \textcolor{cyan}{1}};
\node (ib_1) at (1.5,3.44) {\tiny \textcolor{cyan}{2}};
\node (ib_1) at (1.8,3.4) {\tiny \textcolor{cyan}{3}};
\node (ib_1) at (2.017,3.4) {\tiny \textcolor{cyan}{4}};
\node (ib_1) at (2.175,3.4) {\tiny \textcolor{cyan}{5}};
\node (ib_1) at (2.325,3.4) {\tiny \textcolor{cyan}{6}};
\node (ib_1) at (1.7,3.825) {\tiny \textcolor{blue}{1}};
\node (ib_1) at (2.065,3.115) {\tiny \textcolor{blue}{2}};
\node (ib_1) at (2.275,2.7) {\tiny \textcolor{blue}{3}};
\node (ib_1) at (2.325,2.455) {\tiny \textcolor{blue}{4}};
\node (ib_1) at (2.455,2.46) {\tiny \textcolor{blue}{5}};
\node (ib_1) at (2.565,2.475) {\tiny \textcolor{blue}{6}};
\node (ib_1) at (2.75,2.47) {\tiny \textcolor{blue}{8}};
\node (ib_1) at (2.95,2.4525) {\tiny \textcolor{blue}{10}};
\node (ib_1) at (2.4,3.825) {\tiny \textcolor{magenta}{1}};
\node (ib_1) at (2.68,3.115) {\tiny \textcolor{magenta}{2}};
\node (ib_1) at (2.81,2.705) {\tiny \textcolor{magenta}{3}};
\node (ib_1) at (3.02,2.24) {\tiny \textcolor{magenta}{4}};
\node (ib_1) at (3.0,1.725) {\tiny \textcolor{magenta}{5}};
\node (ib_1) at (3.1285,1.725) {\tiny \textcolor{magenta}{6}};
\node (ib_1) at (3.275,1.75) {\tiny \textcolor{magenta}{8}};
\node (ib_1) at (3.5,1.75) {\tiny \textcolor{magenta}{10}};
\end{tikzpicture}
\label{fig:gaussian_bump_test_case_computational_cost_geopotential_iters_all_256}
}
\hspace{-0.5cm}
\subfigure[]{
\begin{tikzpicture}
\node[anchor=south west,inner sep=0] at (0,0.05){\includegraphics[scale=0.385]{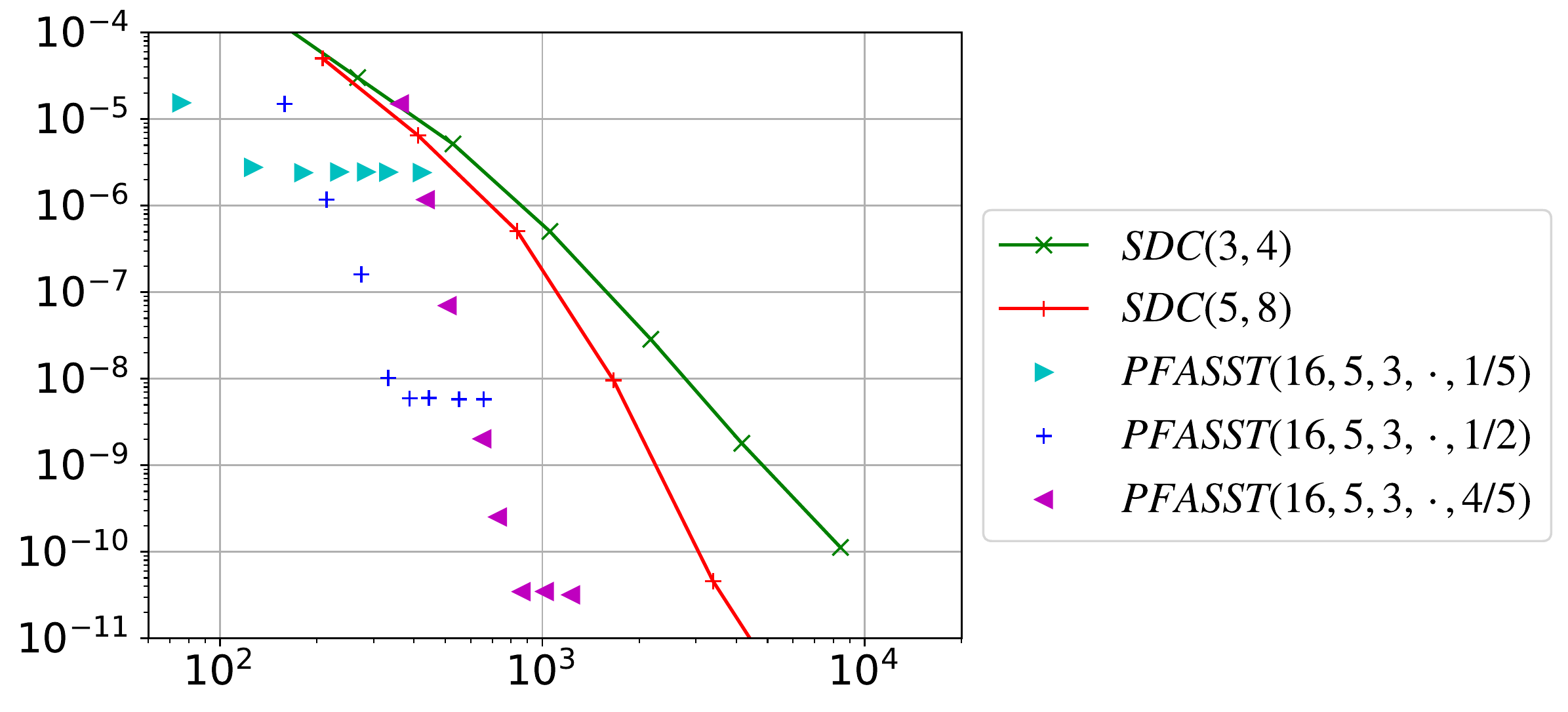}};
\node (ib_1) at (3.2,-0.05) {Wall-clock time $(s)$};
\node (ib_1) at (3.1,4.325) {$R_{\textit{norm}} = 32$};
\end{tikzpicture}
\label{fig:gaussian_bump_test_case_computational_cost_geopotential_iters_all_32}
}
\vspace{-0.45cm}
\caption{\label{fig:gaussian_bump_test_case_computational_cost_geopotential_field_iters_all} 
Gaussian bump: $L_{\infty}$-norm of the error in the geopotential field after
$28 \, \text{hours}$, $26 \, \text{min}$, $40 \, \text{s}$ as a function
of wall-clock time. We fix the time step size to $\Delta t = 100 \, \text{s}$
and increase the number of iterations.
In \oldref{fig:gaussian_bump_test_case_computational_cost_geopotential_iters_all_256},
we annotate the data points with the corresponding number of iterations
that was used to run the simulation. We see that there is a threshold
beyond which increasing the number of iterations does not reduce the error
anymore. This threshold is reached for $N_{\textit{PF}} = 2$ for
$\alpha = 1/5$, for $N_{\textit{PF}} = 4$ for $\alpha = 1/2$, and
for $N_{\textit{PF}} = 6$ for $\alpha = 4/5$.}
\end{figure}

To conclude this section, we motivate the choice of the number of
PFASST-SH iterations by studying the sensitivity of the algorithm
to $N_{\textit{PF}}$ for a fixed time step size ($\Delta t = 100 \, s$).
These results are shown in
Fig.~\oldref{fig:gaussian_bump_test_case_computational_cost_geopotential_field}.
We see that for a small $N_{\textit{PF}}$ (for instance,
$N_{\textit{PF}} < 6$ for PFASST($n_{\textit{ts}} = 16$, $M_f + 1 = 5$,
$M_c + 1 = 3$, $\cdot$, $\alpha = 4/5$)), performing an additional
iteration slightly increases the computational cost, but
significantly reduces the error. Instead, for a larger
$N_{\textit{PF}}$ (larger than the threshold $N_{\textit{PF}} = 6$
for PFASST($n_{\textit{ts}} = 16$, $M_f + 1 = 5$,
$M_c + 1 = 3$, $\cdot$, $\alpha = 4/5$)), increasing the number
of iterations still increases the computational cost, but does
not reduce the error anymore. This remaining error persists
because the high-frequency modes cannot be captured on the
coarse level.

We used this sensitivity analysis to choose the fixed
$N_{\textit{PF}}$ in the refinement studies on the time step size
(Figs.~\oldref{fig:gaussian_bump_test_case_accuracy_geopotential_field}
and \oldref{fig:gaussian_bump_test_case_computational_cost_geopotential_field}).
Specifically, we set $N_{PF}$ slightly above the threshold
mentioned in the previous paragraph to minimize the number
of iterations that increase the computational cost without
reducing the error. This methodology yielded $N_{\textit{PF}} = 4$
for PFASST($n_{\textit{ts}} = 16$, $M_f + 1 = 3$, $M_c + 1 = 2$,
$\cdot$, $\alpha = 4/5$) and $N_{\textit{PF}} = 8$ for
PFASST($n_{\textit{ts}} = 16$, $M_f + 1 = 5$, $M_c + 1 = 3$,
$\cdot$, $\alpha = 4/5$). A more flexible way to choose
$N_{\textit{PF}}$ would rely on a convergence check to stop
the iterations, which we will implement in future work.

\subsection{\label{subsection_rossby_haurwitz_wave}Rossby-Haurwitz wave}

\begin{figure}[ht!]
  \centering
  \subfigure[]{
    \begin{tikzpicture}
    \node[anchor=south west,inner sep=0] at (0,0){\includegraphics[scale=0.385]{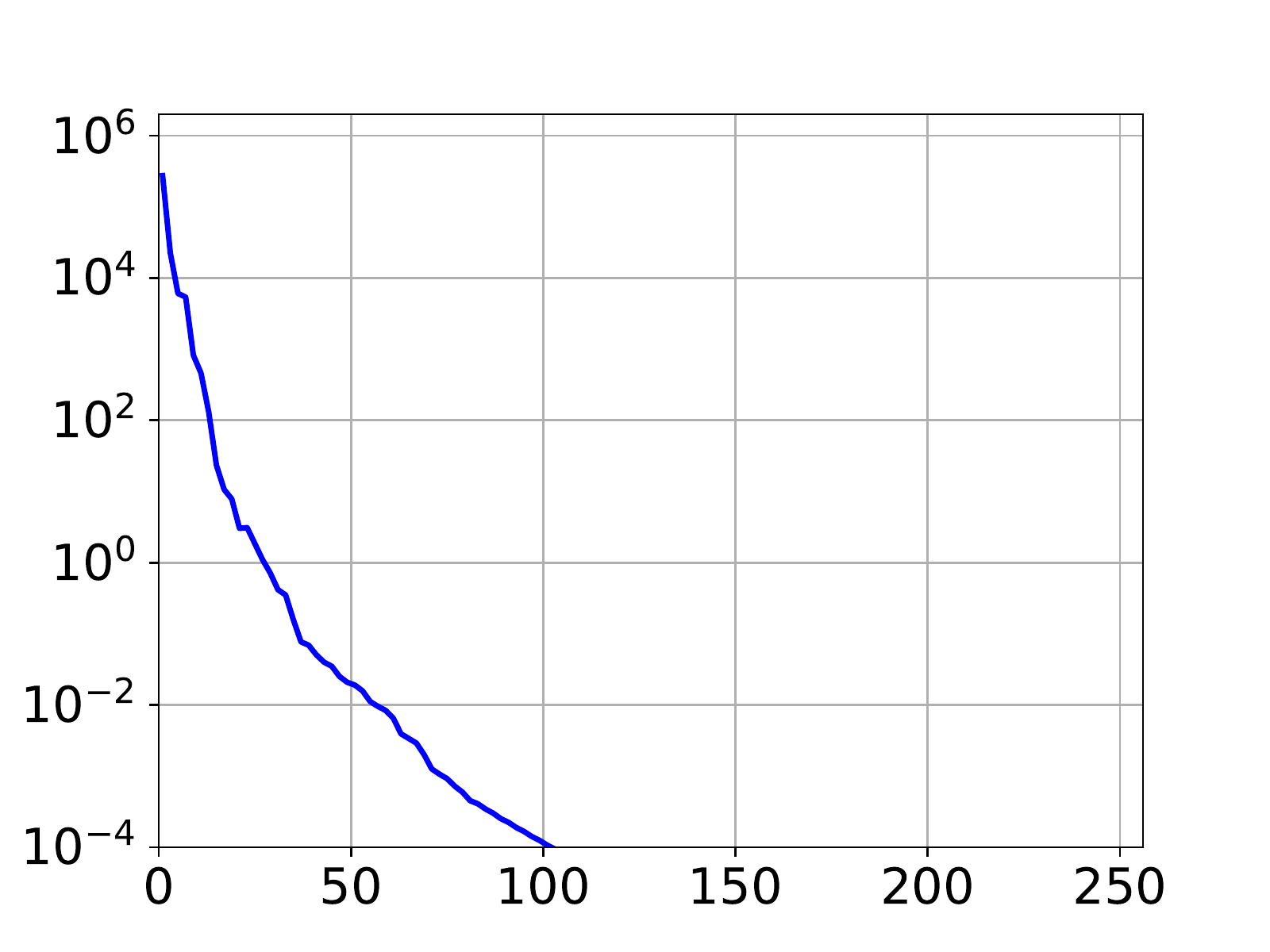}};
    \node (ib_1) at (3.1,-0.05) {$n_0$};
    \node[rotate=90] (ib_1) at (-0.05,2.15) {$|\Phi_{n_0}|$};

    
    \node (ib_1) at (1.75,4.25) {};
    \node (ib_2) at (1.75,0.4) {};
    \path [draw=black!100, thick, dashed] (ib_1) -- (ib_2);

    \node[rotate=-90] (ib_3) at (1.95,2.87) {\scriptsize $\alpha = 1/5$}; 

    \node (ib_1) at (3.21,4.25) {};
    \node (ib_2) at (3.21,0.4) {};
    \path [draw=black!100, thick, dashed] (ib_1) -- (ib_2);

    \node[rotate=-90] (ib_3) at (3.41,2.32) {\scriptsize $\alpha = 1/2$}; 

    \node (ib_1) at (4.68,4.25) {};
    \node (ib_2) at (4.68,0.4) {};
    \path [draw=black!100, thick, dashed] (ib_1) -- (ib_2);

    \node[rotate=-90] (ib_3) at (4.88,1.87) {\scriptsize $\alpha = 4/5$};

    \node (ib_3) at (3.2,4.35) {Geopotential}; 
    \end{tikzpicture}
    \label{fig:rossby_haurwitz_wave_test_case_geopotential_spectrum}
}
  \subfigure[]{
    \begin{tikzpicture}
    \node[anchor=south west,inner sep=0] at (0,0){\includegraphics[scale=0.385]{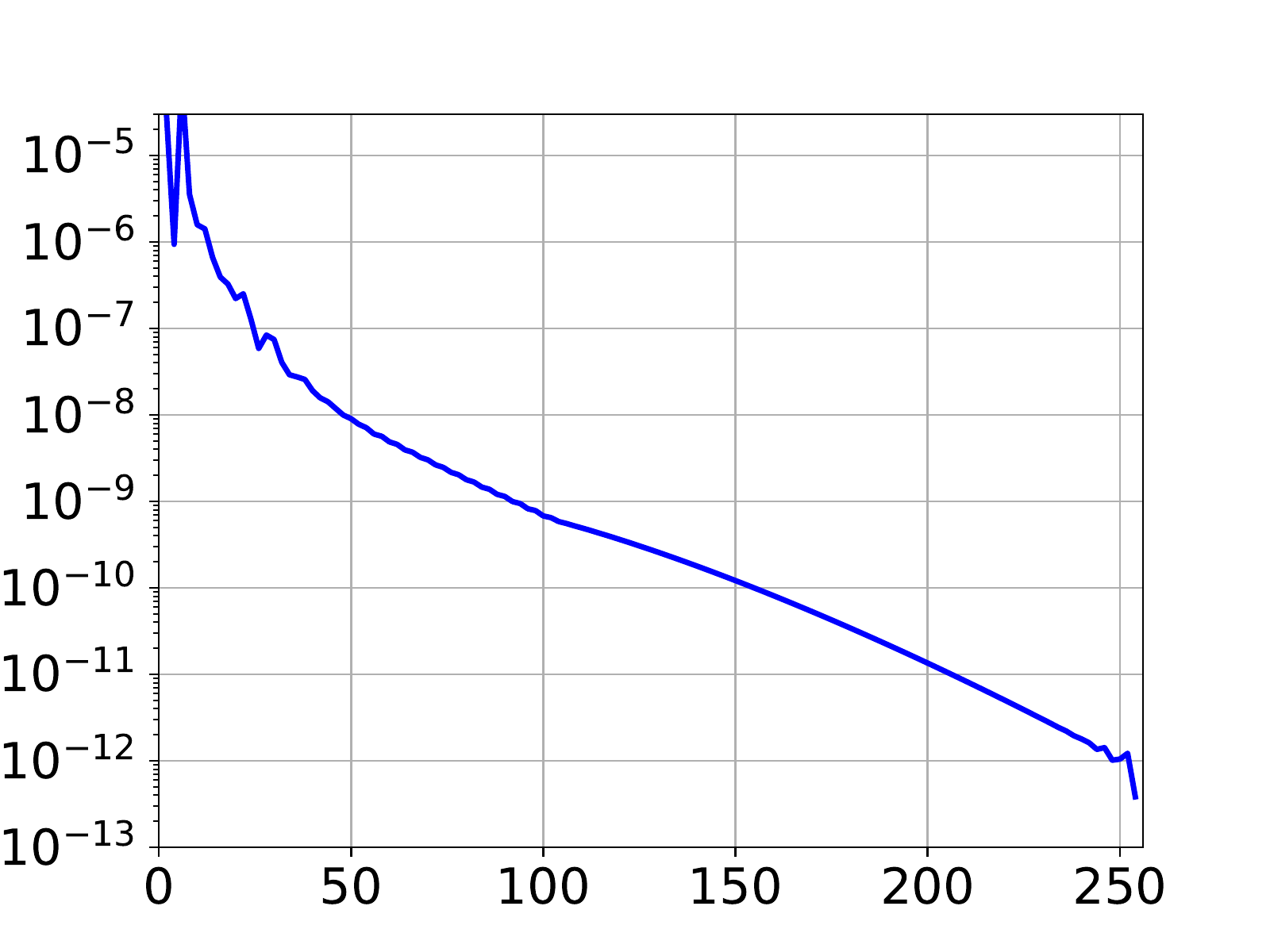}};
    \node (ib_1) at (3.1,-0.05) {$n_0$};
    \node[rotate=90] (ib_1) at (-0.25,2.15) {$|\zeta_{n_0}|$};

    
    \node (ib_1) at (1.75,4.25) {};
    \node (ib_2) at (1.75,0.4) {};
    \path [draw=black!100, thick, dashed] (ib_1) -- (ib_2);

    \node[rotate=-90] (ib_3) at (1.95,3.45) {\scriptsize $\alpha = 1/5$}; 

    \node (ib_1) at (3.21,4.25) {};
    \node (ib_2) at (3.21,0.4) {};
    \path [draw=black!100, thick, dashed] (ib_1) -- (ib_2);

    \node[rotate=-90] (ib_3) at (3.35,2.82) {\scriptsize $\alpha = 1/2$}; 

    \node (ib_1) at (4.68,4.25) {};
    \node (ib_2) at (4.68,0.4) {};
    \path [draw=black!100, thick, dashed] (ib_1) -- (ib_2);

    \node[rotate=-90] (ib_3) at (4.88,1.97) {\scriptsize $\alpha = 4/5$};

    \node (ib_3) at (3.2,4.35) {Vorticity}; 
    \end{tikzpicture}
    \label{fig:rossby_haurwitz_wave_test_case_vorticity_spectrum}
}

\vspace{-0.4cm}
\caption{\label{fig:rossby_haurwitz_wave_test_case_spectrum} 
  Rossby-Haurwitz wave: max-spectrum of the geopotential field in
  \oldref{fig:rossby_haurwitz_wave_test_case_geopotential_spectrum} and of
  the vorticity field in
  \oldref{fig:rossby_haurwitz_wave_test_case_vorticity_spectrum} after
  $28 \, \text{hours}$, $26 \, \text{min}$, $40 \, \text{s}$. These figures
  are obtained using SDC(5,8) with $R_f = 256$ and  a time step size
  $\Delta t_{\textit{ref}} = 100 \, \text{s}$. The quantities on the $y$-axis
  are defined as $|\Phi_{n_0}| = \max_{r} |\Phi^r_{n_0}|$ and
  $|\zeta_{n_0}| = \max_{r} |\zeta^r_{n_0}|$, respectively.
  For each coarsening ratio $\alpha$, a
  vertical dashed line shows the fraction of the spectrum that is truncated during
  the coarsening step in MLSDC-SH and PFASST-SH.
}
\end{figure}

We continue the analysis of the PFASST-SH algorithm
with the Rossby-Haurwitz wave test case. The initial
velocity field is non-divergent, and the initial
geopotential distribution, given in \cite{williamson1992standard},
is such that the temporal derivative of the divergence
is zero. The diffusion coefficient is set
to $10^5 \, \text{m}^2.\text{s}^{-1}$.
We simulate $102,400 \, \text{s}$ ($28 \, \text{hours}$,
$26 \, \text{min}$, $40 \, \text{s}$) of propagation with a spatial
resolution given by $R_{f} = 256$.
We note that when PFASST-SH is used with an
aggressive coarsening ratio and/or large time steps to simulate
the problem over a longer period of time, the wave becomes unstable,
as documented for standard time integration schemes in \cite{thuburn2000numerical,benard2019numerical}.
The spectrum of the reference
solution obtained using SDC($M_f + 1 = 5$, $N_S = 8$) with the same
spatial resolution and a time step size
$\Delta t_{\textit{ref}} = 100 \, \text{s}$ is shown in
Fig.~\oldref{fig:rossby_haurwitz_wave_test_case_spectrum}. 
As in Section~\oldref{subsection_propagation_of_a_gaussian_dome},
we see that even though relatively large spectral coefficients are
truncated on the coarse level for $\alpha = 1/5$, they are
much smaller for $\alpha = 4/5$.

We first study the error norm as a function of time step size
for a fixed number of iterations. We consider the same PFASST-SH
schemes with the same number of iterations as in the previous
test case, namely, PFASST($n_{\textit{ts}} = 16$, $M_f + 1=3$, 
$M_c + 1=2$, $N_{\textit{PF}}=4$, $\alpha$) and
PFASST($n_{\textit{ts}} = 16$, $M_f + 1 = 5$, $M_c + 1 = 3$,
$N_{\textit{PF}} = 8$, $\alpha$). Here, the fixed
number of iterations, $N_{\textit{PF}}$, has been chosen using
the approach described at the end of the previous section.

\begin{figure}[ht!]
\centering
\subfigure[]{
\begin{tikzpicture}
\node[anchor=south west,inner sep=0] at (0,0.045){\includegraphics[scale=0.385]{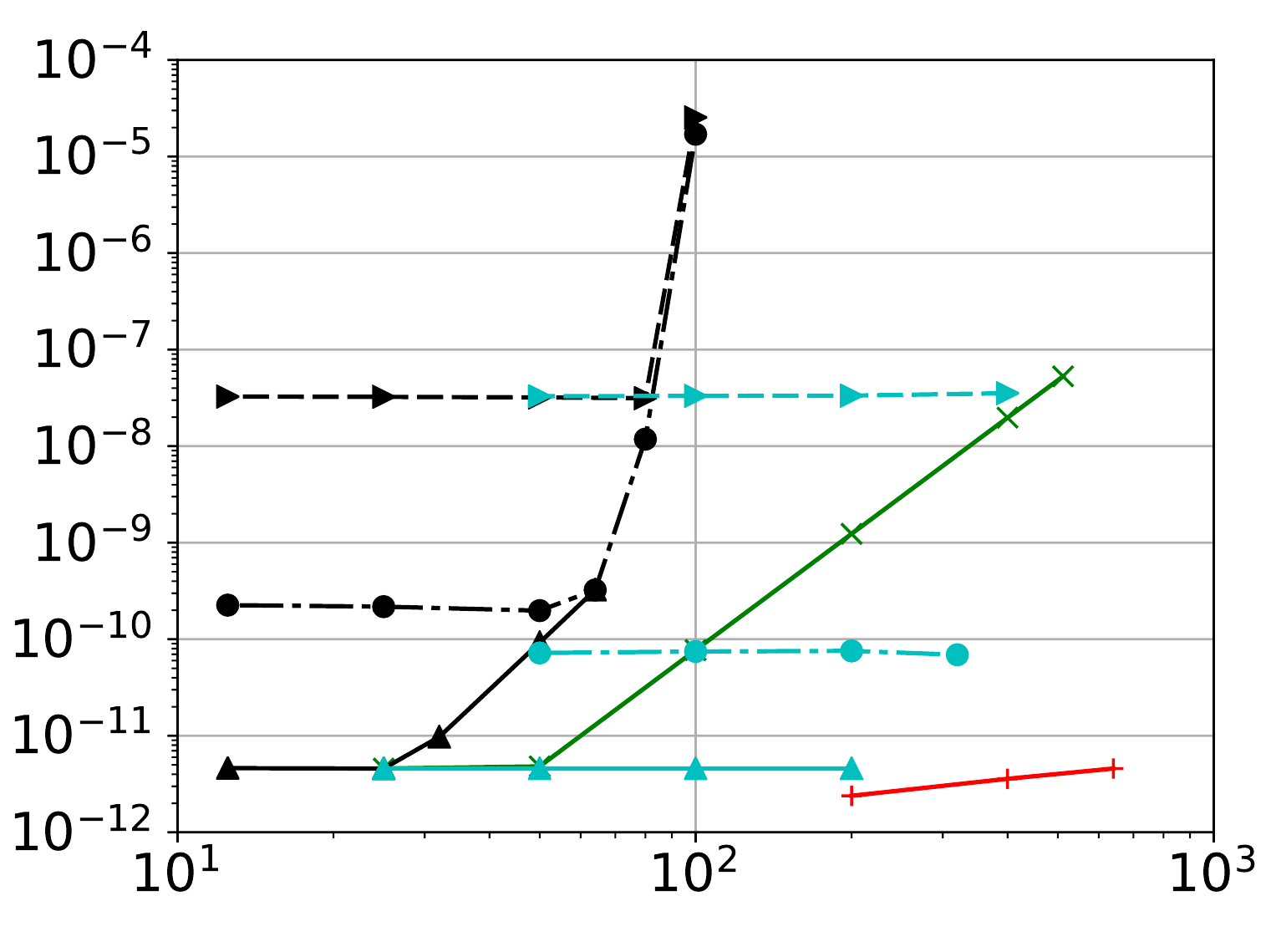}};
\node (ib_1) at (3.2,-0.05) {$\Delta t$ $(s)$};
\node[rotate=90] (ib_1) at (-0.1,2.35) {Normalized $L_{\infty}$-error};
\node (ib_1) at (3.3,4.4) {$R_{\textit{norm}} = 256$};

\node (ib_1) at (1.45,2.82) {\scriptsize $\alpha = 1/5$};
\node (ib_1) at (1.45,1.85) {\scriptsize $\alpha = 1/2$};
\node (ib_1) at (1.45,1.05) {\scriptsize $\alpha = 4/5$};

\end{tikzpicture}
\label{fig:rossby_haurwitz_wave_test_case_accuracy_geopotential_256}
}
\hspace{-0.5cm}
\subfigure[]{
\begin{tikzpicture}
\node[anchor=south west,inner sep=0] at (0,0.05){\includegraphics[scale=0.385]{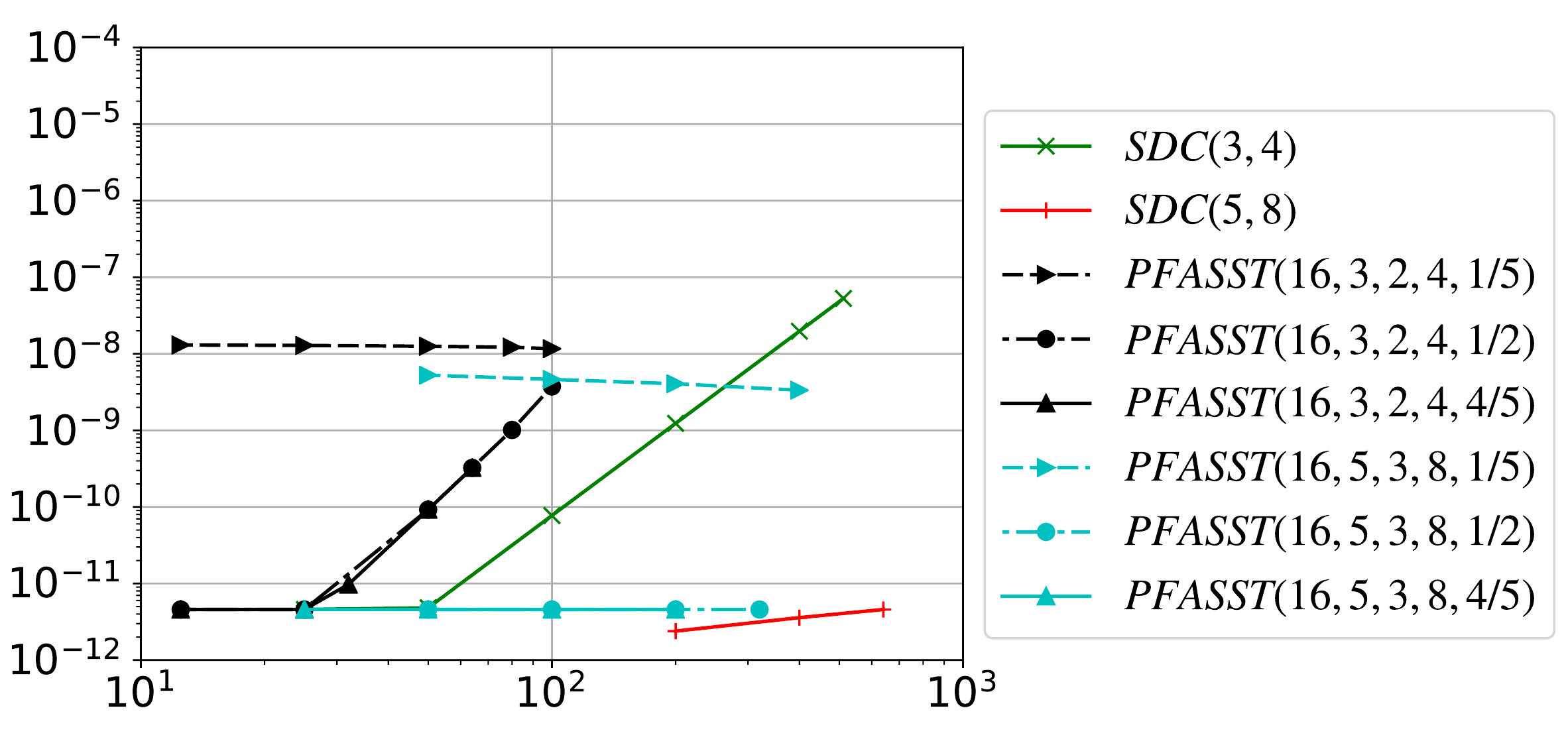}};
\node (ib_1) at (3.2,-0.05) {$\Delta t$ $(s)$};
\node (ib_1) at (3.1,4.45) {$R_{\textit{norm}} = 32$};
\end{tikzpicture}
\label{fig:rossby_haurwitz_wave_test_case_accuracy_geopotential_32}
}
\vspace{-0.45cm}
\caption{\label{fig:rossby_haurwitz_wave_test_case_accuracy_geopotential_field} 
  Rossby-Haurwitz wave: $L_{\infty}$-norm of the error in the geopotential field
  after $28 \, \text{hours}$, $26 \, \text{min}$, $40 \, \text{s}$
  as a function of time step size. We apply the same norm to the SDC and PFASST-SH
  schemes, with $R_{\text{norm}} = 256$ in
  \oldref{fig:rossby_haurwitz_wave_test_case_accuracy_geopotential_256}, and
  $R_{\text{norm}} = 32$
  in \oldref{fig:rossby_haurwitz_wave_test_case_accuracy_geopotential_32}.
 In the legend, the fourth (respectively, fifth) parameter of the PFASST-SH
  schemes denotes the number of iterations (respectively, the coarsening ratio in
    space).
}
\end{figure}

The results shown in
Fig.~\oldref{fig:rossby_haurwitz_wave_test_case_accuracy_geopotential_field}
are consistent with those of the previous test case.
Specifically, the largest stable time steps achieved with
PFASST-SH decrease when the spatial coarsening ratio is
close to one, and are still smaller than with the serial
SDC schemes. In addition, the error obtained with the
PFASST-SH schemes first decreases as the time step size
is reduced, and then stagnates for small time steps. As
discussed earlier, the error norm achieved for small
time steps depends on the spatial coarsening ratio as
this parameter determines the fraction of the spectrum
that can be represented on the coarse level.
Figure~\oldref{fig:rossby_haurwitz_wave_test_case_accuracy_geopotential_field}
shows that the error norm remains large with
aggressive spatial coarsening, but essentially vanishes
for mild coarsening in space. We observe that
PFASST($n_{\textit{ts}} = 16$, $M_f + 1 = 5$, 
$M_c + 1 = 3$, $N_{\textit{PF}} = 8$, $\alpha = 4/5$) is
as accurate as the serial SDC($M_f + 1 = 5$, $N_S = 8$)
for all the time step sizes considered here.
PFASST($n_{\textit{ts}} = 16$, $M_f + 1 = 5$, $M_c + 1 = 3$,
$N_{\textit{PF}} = 8$, $\alpha = 1/2$) is less accurate
than PFASST(16,5,3,8,4/5) for $R_{\textit{norm}} = 256$, 
but the two schemes have the same accuracy for
$R_{\textit{norm}} = 32$.

\begin{figure}[ht!]
\centering
\subfigure[]{
\begin{tikzpicture}
\node[anchor=south west,inner sep=0] at (0,0.045){\includegraphics[scale=0.385]{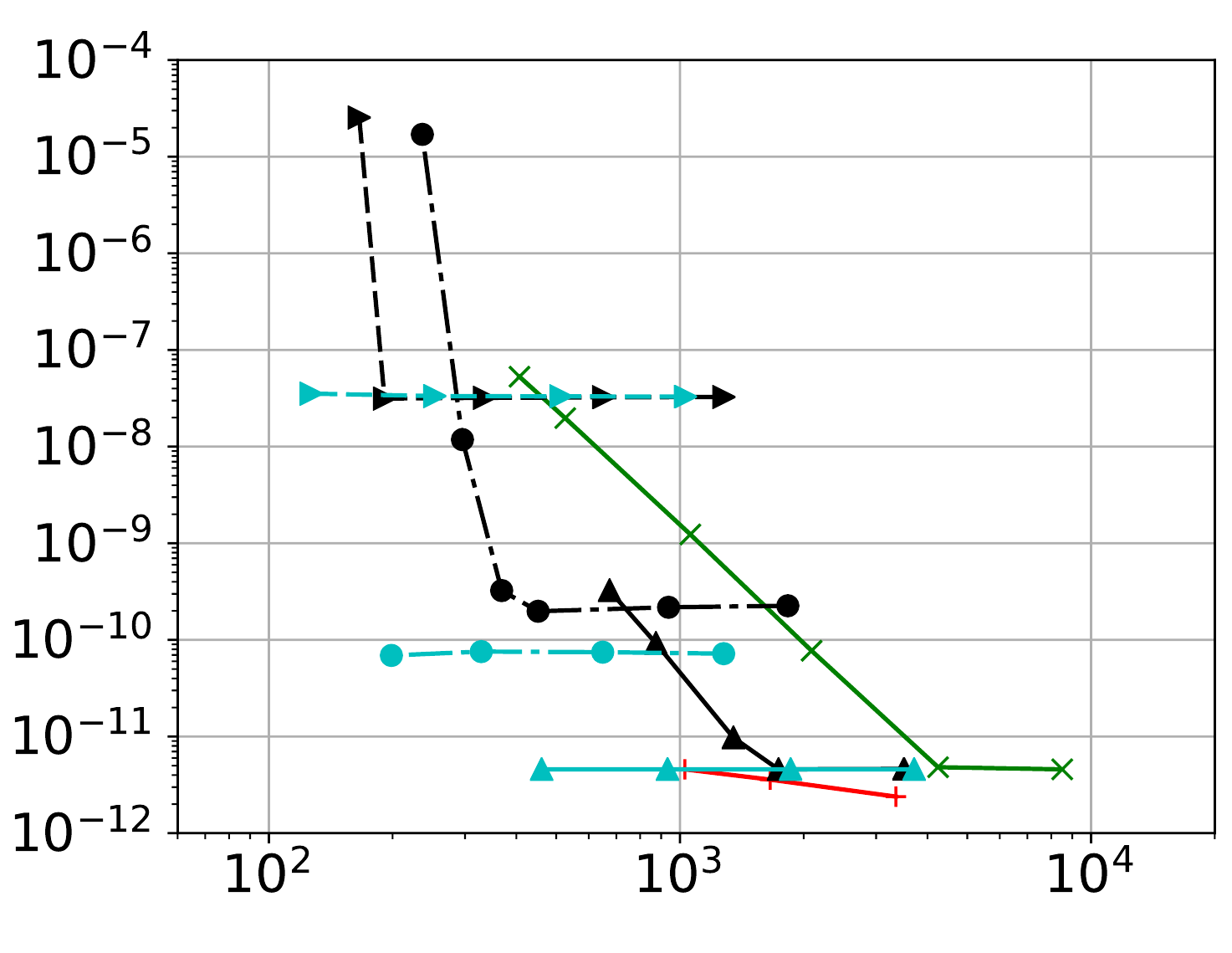}};
\node (ib_1) at (3.2,-0.05) {Wall-clock time $(s)$};
\node[rotate=90] (ib_1) at (-0.1,2.35) {Normalized $L_{\infty}$-error};
\node (ib_1) at (3.3,4.45) {$R_{\textit{norm}} = 256$};

\node (ib_1) at (3.6,2.82) {\scriptsize $\alpha = 1/5$};
\node (ib_1) at (1.6,1.65) {\scriptsize $\alpha = 1/2$};
\node (ib_1) at (2.4,1.05) {\scriptsize $\alpha = 4/5$};

\end{tikzpicture}
\label{fig:rossby_haurwitz_wave_test_case_computational_cost_geopotential_256}
}
\hspace{-0.5cm}
\subfigure[]{
\begin{tikzpicture}
\node[anchor=south west,inner sep=0] at (0,0.05){\includegraphics[scale=0.385]{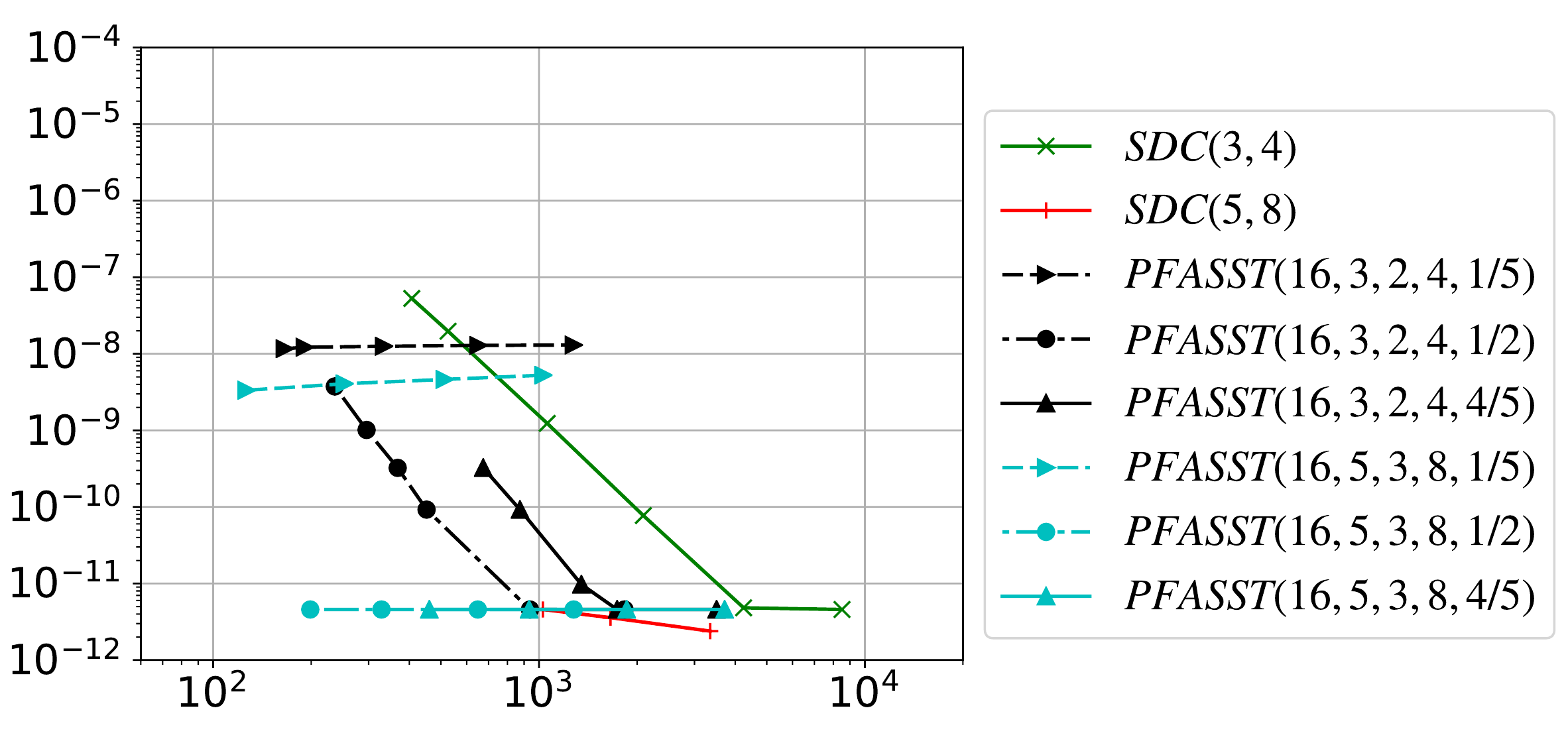}};
\node (ib_1) at (3.2,-0.05) {Wall-clock time $(s)$};
\node (ib_1) at (3.1,4.45) {$R_{\textit{norm}} = 32$};

\node (ib_1) at (3.8,2.65) {\scriptsize $\alpha = 1/5$};
\node (ib_1) at (1.9,1.06) {\scriptsize $\alpha = 1/2, \, 4/5$};

\end{tikzpicture}
\label{fig:rossby_haurwitz_wave_test_case_computational_cost_geopotential_32}
}
\vspace{-0.45cm}
\caption{\label{fig:rossby_haurwitz_wave_test_case_computational_cost_geopotential_field} 
  Rossby-Haurwitz wave: $L_{\infty}$-norm of the error in the geopotential field after
  $28 \, \text{hours}$, $26 \, \text{min}$, $40 \, \text{s}$ as a function of wall-clock time.
  We apply the same norm to the SDC and PFASST-SH schemes, with $R_{\text{norm}} = 256$ in 
  \oldref{fig:rossby_haurwitz_wave_test_case_computational_cost_geopotential_256}, and
  $R_{\text{norm}} = 32$ in \oldref{fig:rossby_haurwitz_wave_test_case_computational_cost_geopotential_32}.
  In the legend, the fourth (respectively, fifth) parameter of the PFASST-SH schemes
  denotes the number of iterations (respectively, the coarsening ratio in space). 
  For their largest
  time steps, we see that the PFASST-SH schemes have a smaller wall-clock time
  than the serial SDC schemes for a similar error.
}
\end{figure}  

Next, we assess the computational cost of these PFASST-SH schemes by considering
the error norm as a function of wall-clock time.
As in the previous test case, we note in
Fig.~\oldref{fig:rossby_haurwitz_wave_test_case_computational_cost_geopotential_field}
that for their largest stable time steps, the PFASST-SH schemes are more efficient
than the corresponding serial SDC schemes.
We also see that for a given accuracy, the difference
in wall-clock time is larger when the error is computed with $R_{\textit{norm}} = 32$,
which suggests that the PFASST-SH algorithm resolves the main features of the solution
significantly faster than
the serial schemes.

These results are summarized in Table~\oldref{tab:rossby_haurwitz_wave_timings}.
We observe that PFASST($n_{\textit{ts}} = 16$, $M_f + 1 = 3$, $M_c + 1 = 2$,
$N_{\textit{ML}} = 4$, $\alpha = 1/5$) can achieve a geopotential error norm of
$2 \times 10^{-8}$ for a time step size of $\Delta t  = 100 \, \text{s}$, which
is  a quarter of the time step size required by SDC($M_f + 1 = 3$, $N_S = 4$)
and MLSDC($M_f +1 = 3$, $M_c +1 = 2$, $N_{\textit{ML}} = 2$,
$\alpha = 1/2$) to reach this tolerance. The reduced stability of the parallel-in-time
scheme explains why the speedups obtained
with PFASST($n_{\textit{ts}} = 16$, $M_f + 1 = 3$, $M_c + 1 = 2$, $N_{\textit{PF}} = 4$,
$\alpha = 1/5$) are slightly smaller than for the previous test case. We observe
a speedup of  $S^{\textit{obs}}_S = 3.2$ compared to SDC($M_f + 1 = 3$, $N_S = 4$) and of
$S^{\textit{obs}}_{ML} = 2.2$ compared to MLSDC($M_f + 1 = 3$, $M_c + 1 = 2$, $N_{\textit{ML}} = 2$,
$\alpha = 1/2$), which is close to the theoretical speedup of Section~\oldref{subsection_computational_cost}.
This remark also applies to the speedup obtained with PFASST($16$, $M_f + 1 = 5$,
$M_c + 1 = 3$, $N_{\textit{PF}} = 8$, $\alpha = 1/2$). 

\begin{table}[!ht]
\centering
\scalebox {0.9}{
           \begin{tabular}{lccccc}
           \\ \toprule 
           Integration         &    $L_{\infty}$-error & Time step & Wall-clock & Observed & Observed  \\
           scheme              &    $R_{\textit{norm}} = 32$ & size (s)  & time (s)   & speedup $S^{\textit{obs}}_{S}$ & speedup $S^{\textit{obs}}_{\textit{ML}}$ \\\toprule
           SDC(3,4)                  &     $2.0 \times 10^{-8}$        &     400      &   526   & -    & -    \\
           MLSDC(3,2,2,1/2)          &     $1.7 \times 10^{-8}$        &     400      &   366   & 1.4  & -    \\
           PFASST(16,3,2,4,1/5)      &     $1.2 \times 10^{-8}$        &     100      &   166    & 3.2   & 2.2  \\ \midrule
           SDC(5,8)                  &     $4.6 \times 10^{-12}$       &     640      &   1,027   & -    & -    \\
           MLSDC(5,3,4,1/2)          &     $4.6 \times 10^{-12}$       &     640      &   710     & 1.4  & -    \\
           PFASST(16,5,3,8,1/2)      &     $4.6 \times 10^{-12}$       &     320       &  198     & 5.2   & 3.6  \\
           \bottomrule 
         \end{tabular}}
\caption{\label{tab:rossby_haurwitz_wave_timings}
  Rossby-Haurwitz wave: $L_{\infty}$-error computed with $R_{\textit{norm}} = 32$, time
  step size, and wall-clock time for the serial and parallel SDC-based schemes.
  $S^{\textit{obs}}_S$ denotes the observed speedup
  achieved with MLSDC-SH and PFASST-SH with respect to SDC, while $S^{\textit{obs}}_{\textit{ML}}$
  denotes the speedup achieved with PFASST-SH with respect  to MLSDC-SH.   The speedups for $R_{\textit{norm}} = 256$ can be read from
  Fig.~\oldref{fig:rossby_haurwitz_wave_test_case_computational_cost_geopotential_256}.
}
\end{table}
\begin{figure}[ht!]
\centering
\subfigure[]{
\begin{tikzpicture}
\node[anchor=south west,inner sep=0] at (0,0){\includegraphics[scale=0.27]{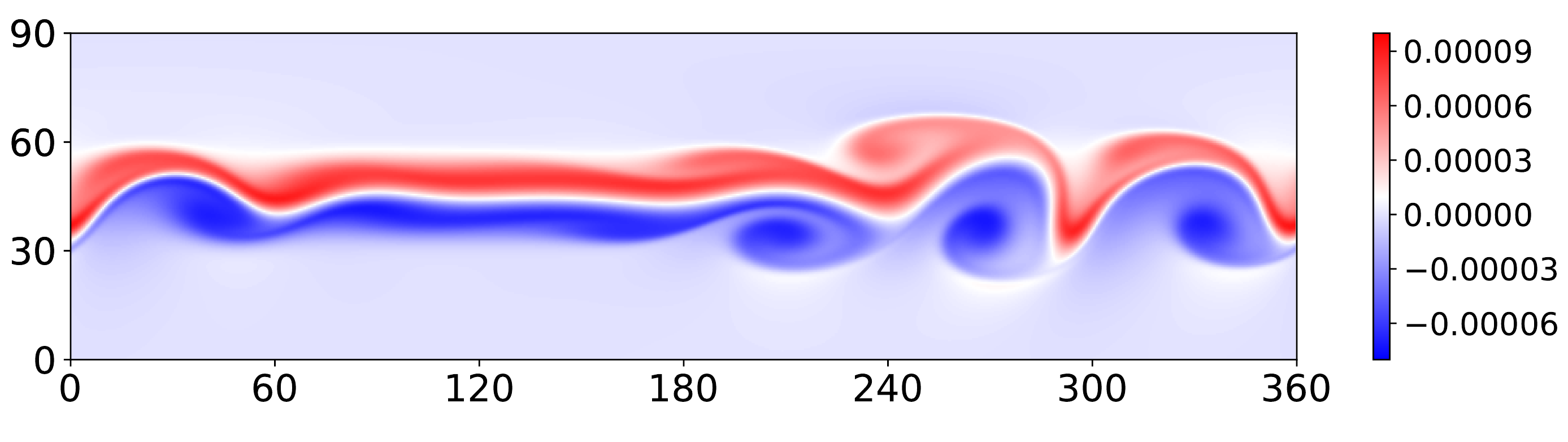}};
\node (ib_1) at (3.4,-0.05) {\scriptsize Longitude (degrees)};
\node[rotate=90] (ib_1) at (-0.12,1.2) {\scriptsize Latitude (degrees)};
\end{tikzpicture}
\label{fig:galewsky_with_bump_test_case_vorticity_field_256_visc_1e5}
}
\hspace{-0.4cm}
\subfigure[]{
\begin{tikzpicture}
\node[anchor=south west,inner sep=0] at (0,0){\includegraphics[scale=0.27]{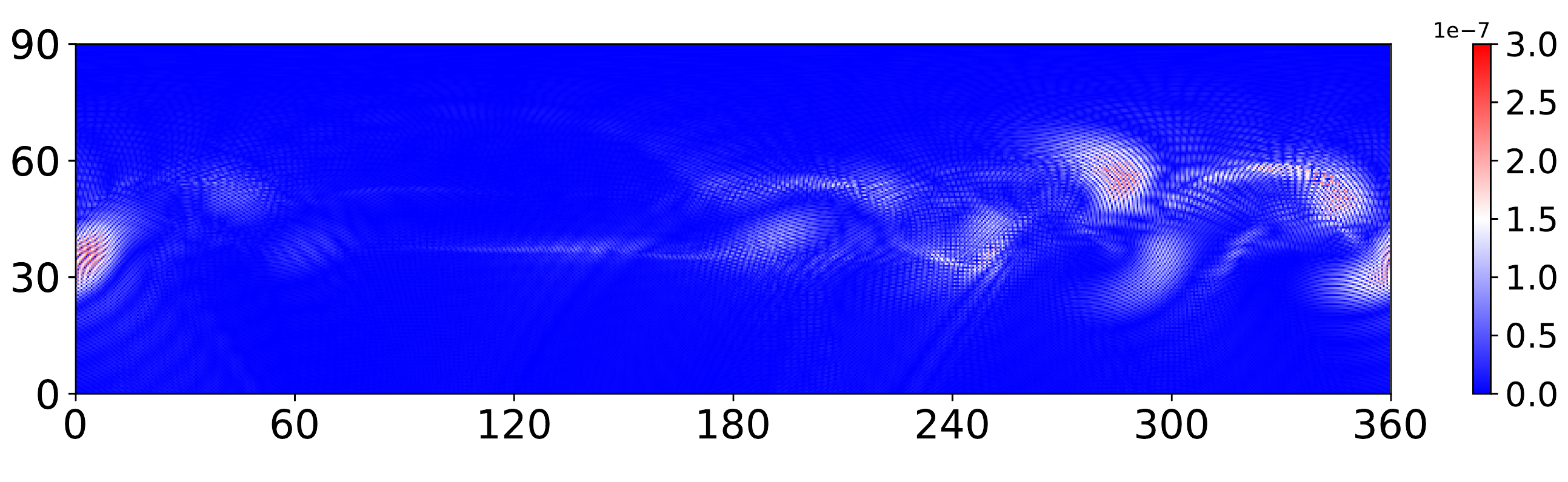}};
\node (ib_1) at (3.4,-0.05) {\scriptsize Longitude (degrees)};
\path (6.5,2.015) node (d) {};
\path (6.8,2.115) node (e) {};
\path [draw=white,fill=white] (d) rectangle (e); 
\path (7.3,1.995) node (e) {\scalebox{0.45}{$\times 10^{-7}$}};
\end{tikzpicture}
\label{fig:galewsky_with_bump_test_case_vorticity_field_error_256_visc_1e5}
}
\vspace{-0.5cm}
\caption{\label{fig:galewsky_with_bump_test_case_vorticity_field_and_error_visc_1e5} 
  Unstable barotropic wave: vorticity field obtained using PFASST(16,5,3,10,4/5) with
  a time step size $\Delta t = 240 \, \text{s}$ in  
  \oldref{fig:galewsky_with_bump_test_case_vorticity_field_256_visc_1e5}.
  This PFASST-SH scheme relies on ten iterations and a spatial coarsening ratio
  $\alpha = 4/5$. Figure
  \oldref{fig:galewsky_with_bump_test_case_vorticity_field_error_256_visc_1e5} shows the
  normalized difference (in the physical coefficients) in the vorticity field of this solution with respect to the reference solution
  obtained using SDC(5,8) with a time step size $\Delta t_{\textit{ref}} = 60 \, \text{s}$. 
  These maps are obtained after 144 hours with a resolution of $R_f = 256$.
}
\end{figure}


\subsection{\label{subsection_nonlinear_evolution_of_an_unstable_barotropic_wave}Nonlinear evolution of an unstable barotropic wave}

We conclude the analysis of PFASST-SH by considering the initial condition
proposed by \cite{galewsky2004initial}. This test case consists of a stationary zonal
jet perturbed by the introduction of a Gaussian bump in the geopotential field. This leads 
to the propagation of gravity waves followed by the development of complex
vortical dynamics. These processes operate on multiple time scales and are
representative of the horizontal aspects of atmospheric flows. 
We run the simulations for 144 hours using a diffusion coefficient of
$\nu = 10^5 \, \text{m}^2.\text{s}^{-1}$ as in \cite{galewsky2004initial}.

We have
shown in \cite{hamon2018multi} that this is a particularly challenging test 
for the MLSDC-SH scheme. This is due to the fast amplification of small-scale
features that significantly undermines the accuracy of the multi-level integration
scheme when the spatial coarsening strategy is too aggressive. Based on a result
from our previous work on MLSDC-SH, we only consider coarsening ratios closer
to one -- i.e., $\alpha = R_c / R_f \geq 3/5$ -- in this example. This test
case is therefore well suited to assess the robustness and efficiency of
PFASST-SH when the coarse corrections are relatively expensive due to the limited
spatial coarsening. 

\begin{figure}[ht!]
  \centering
  \subfigure[]{
    \begin{tikzpicture}
    \node[anchor=south west,inner sep=0] at (0,0){\includegraphics[scale=0.385]{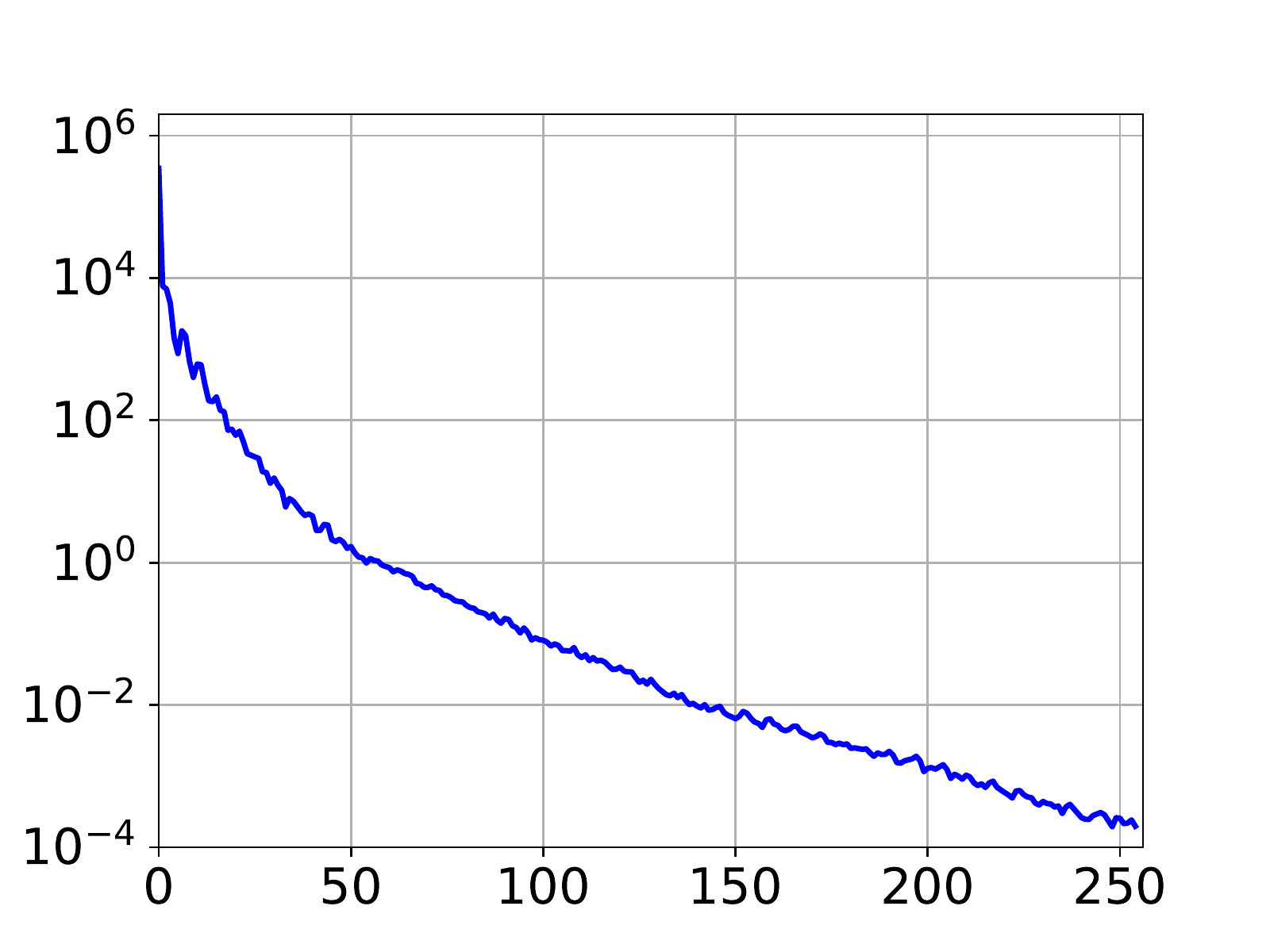}};
    \node (ib_1) at (3.1,-0.05) {$n_0$};
    \node[rotate=90] (ib_1) at (-0.05,2.15) {$|\Phi_{n_0}|$};

    
    \node (ib_1) at (3.7,4.25) {};
    \node (ib_2) at (3.7,0.4) {};
    \path [draw=black!100, thick, dashed] (ib_1) -- (ib_2);

    \node[rotate=-90] (ib_3) at (3.525,2.32) {\scriptsize $\alpha = 3/5$}; 

    \node (ib_1) at (4.68,4.25) {};
    \node (ib_2) at (4.68,0.4) {};
    \path [draw=black!100, thick, dashed] (ib_1) -- (ib_2);

    \node[rotate=-90] (ib_3) at (4.505,2.05) {\scriptsize $\alpha = 4/5$}; 


    \node[rotate=-90] (ib_3) at (5.475,1.7) {\scriptsize $\alpha = 1$}; 

    \node (ib_3) at (3.2,4.35) {Geopotential}; 
    \end{tikzpicture}
    \label{fig:galewsky_with_bump_test_case_geopotential_spectrum}
}
  \subfigure[]{
    \begin{tikzpicture}
    \node[anchor=south west,inner sep=0] at (0,0){\includegraphics[scale=0.385]{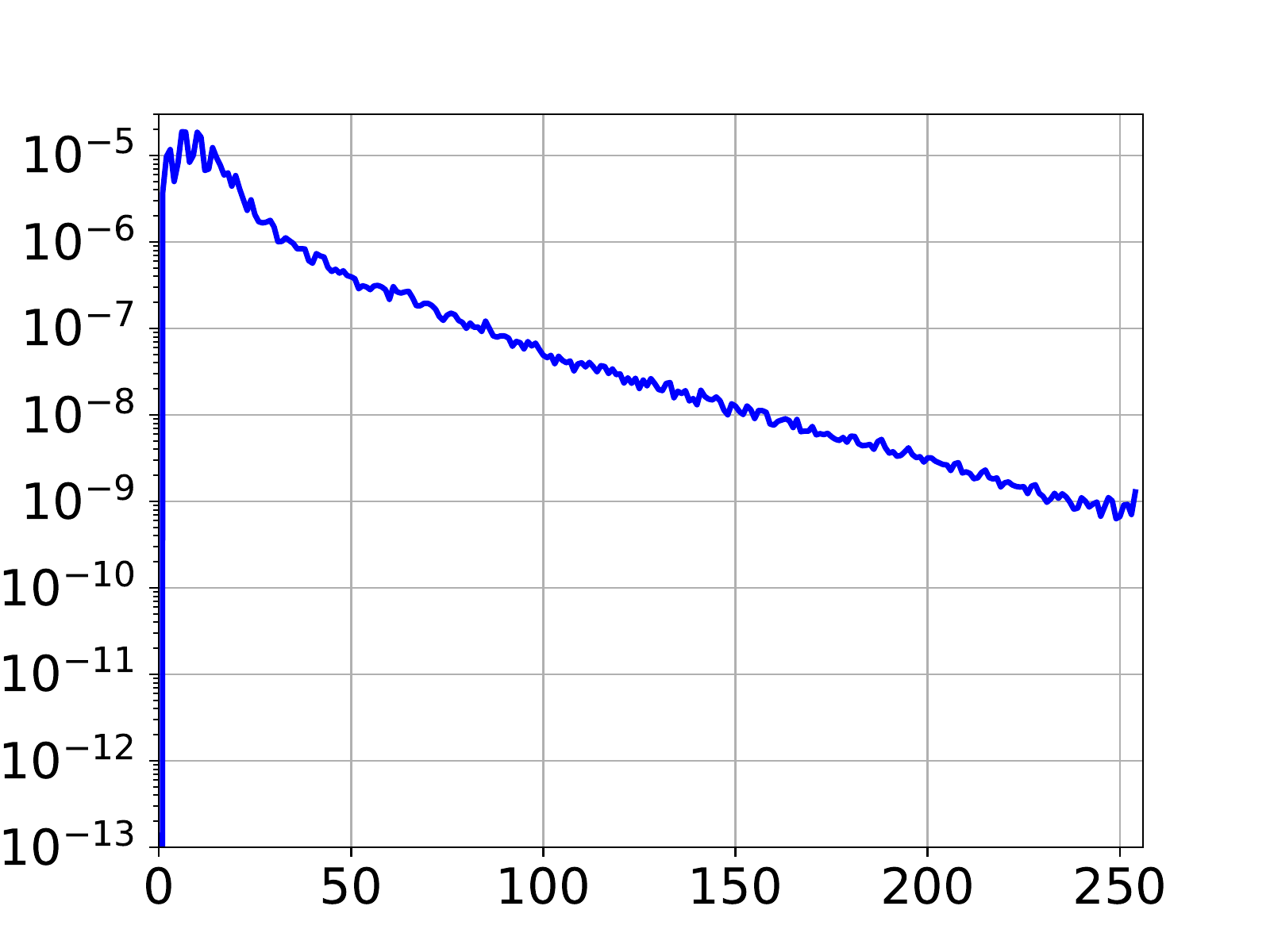}};
    \node (ib_1) at (3.1,-0.05) {$n_0$};
    \node[rotate=90] (ib_1) at (-0.25,2.15) {$|\zeta_{n_0}|$};

    
    \node (ib_1) at (3.7,4.25) {};
    \node (ib_2) at (3.7,0.4) {};
    \path [draw=black!100, thick, dashed] (ib_1) -- (ib_2);

    \node[rotate=-90] (ib_3) at (3.525,1.92) {\scriptsize $\alpha = 3/5$}; 

    \node (ib_1) at (4.68,4.25) {};
    \node (ib_2) at (4.68,0.4) {};
    \path [draw=black!100, thick, dashed] (ib_1) -- (ib_2);

    \node[rotate=-90] (ib_3) at (4.505,1.65) {\scriptsize $\alpha = 4/5$}; 


    \node[rotate=-90] (ib_3) at (5.475,1.3) {\scriptsize $\alpha = 1$}; 

    \node (ib_3) at (3.2,4.35) {Vorticity}; 
    \end{tikzpicture}
    \label{fig:galewsky_with_bump_test_case_vorticity_spectrum}
}

\vspace{-0.4cm}
\caption{\label{fig:galewsky_with_bump_test_case_spectrum} 
  Unstable barotropic wave: max-spectrum of the geopotential field in
  Fig.~\oldref{fig:galewsky_with_bump_test_case_geopotential_spectrum} and of
  the vorticity field in Fig.~\oldref{fig:galewsky_with_bump_test_case_vorticity_spectrum}
  obtained after $144 \, \text{hours}$ using SDC(5,8) with $R_f =256$ and a time step size
  $\Delta t_{\textit{ref}} = 60 \, \text{s}$. The quantities on the $y$-axis are
  defined as $|\Phi_{n_0}| = \max_{r} |\Phi^r_{n_0}|$ and
  $|\Phi_{n_0}| = \max_{r} |\Phi^r_{n_0}|$, respectively. For each coarsening ratio $\alpha$, a
  vertical dashed line shows the fraction of the spectrum that is truncated during
  the coarsening step in MLSDC-SH and PFASST-SH.
}
\end{figure}

The vortices generated with PFASST-SH  are  in
Fig.~\oldref{fig:galewsky_with_bump_test_case_vorticity_field_and_error_visc_1e5}
for a resolution given by $R_f = 256$, 
along with the error in the vorticity field with respect to a reference solution
computed with SDC($M_f + 1=5$, $N_S = 8$). This reference solution relies on a time step size
$\Delta t_{\textit{ref}} = 60 \, \text{s}$ and the same spatial resolution.
Figure~\oldref{fig:galewsky_with_bump_test_case_vorticity_field_and_error_visc_1e5}
shows that PFASST-SH can accurately match the vorticity field generated with
SDC($M_f + 1 = 5$, $N_S = 8$), with small errors concentrated in the vortical structures of the flow.
The spectrum of the reference is shown in
Fig.~\oldref{fig:galewsky_with_bump_test_case_spectrum}. Comparing the vorticity
spectrum with those of Figs.~\oldref{fig:gaussian_bump_test_case_spectrum} and
\oldref{fig:rossby_haurwitz_wave_test_case_spectrum} highlights the
amplification
of large spectral coefficients corresponding to high-frequency features in this
test case.

\begin{figure}[ht!]
\centering
\subfigure[]{
\begin{tikzpicture}
\node[anchor=south west,inner sep=0] at (0,0.045){\includegraphics[scale=0.385]{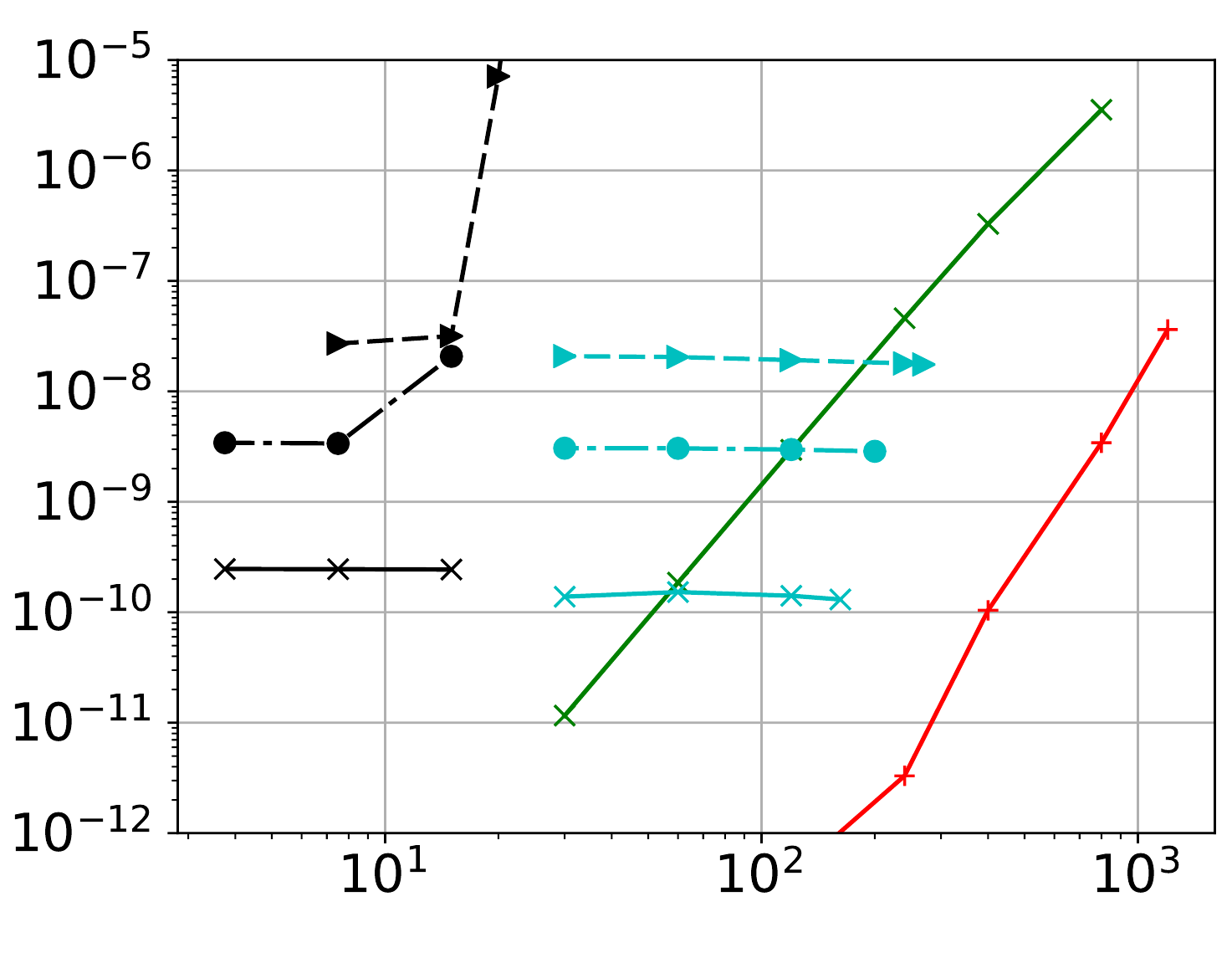}};
\node (ib_1) at (3.2,-0.05) {$\Delta t$ $(s)$};
\node[rotate=90] (ib_1) at (-0.1,2.45) {Normalized $L_{\infty}$-error};
\node (ib_1) at (3.3,4.45) {$R_{\textit{norm}} = 256$};

\node (ib_1) at (3.1,3.) {\scriptsize $\alpha = 3/5$};
\node (ib_1) at (2.88,2.53) {\scriptsize $\alpha = 4/5$};
\node (ib_1) at (2.75,1.855) {\scriptsize $\alpha = 1$};

\end{tikzpicture}
\label{fig:galewsky_with_bump_test_case_accuracy_geopotential_256_visc_1e5}
}
\hspace{-0.5cm}
\subfigure[]{
\begin{tikzpicture}
\node[anchor=south west,inner sep=0] at (0,0.05){\includegraphics[scale=0.385]{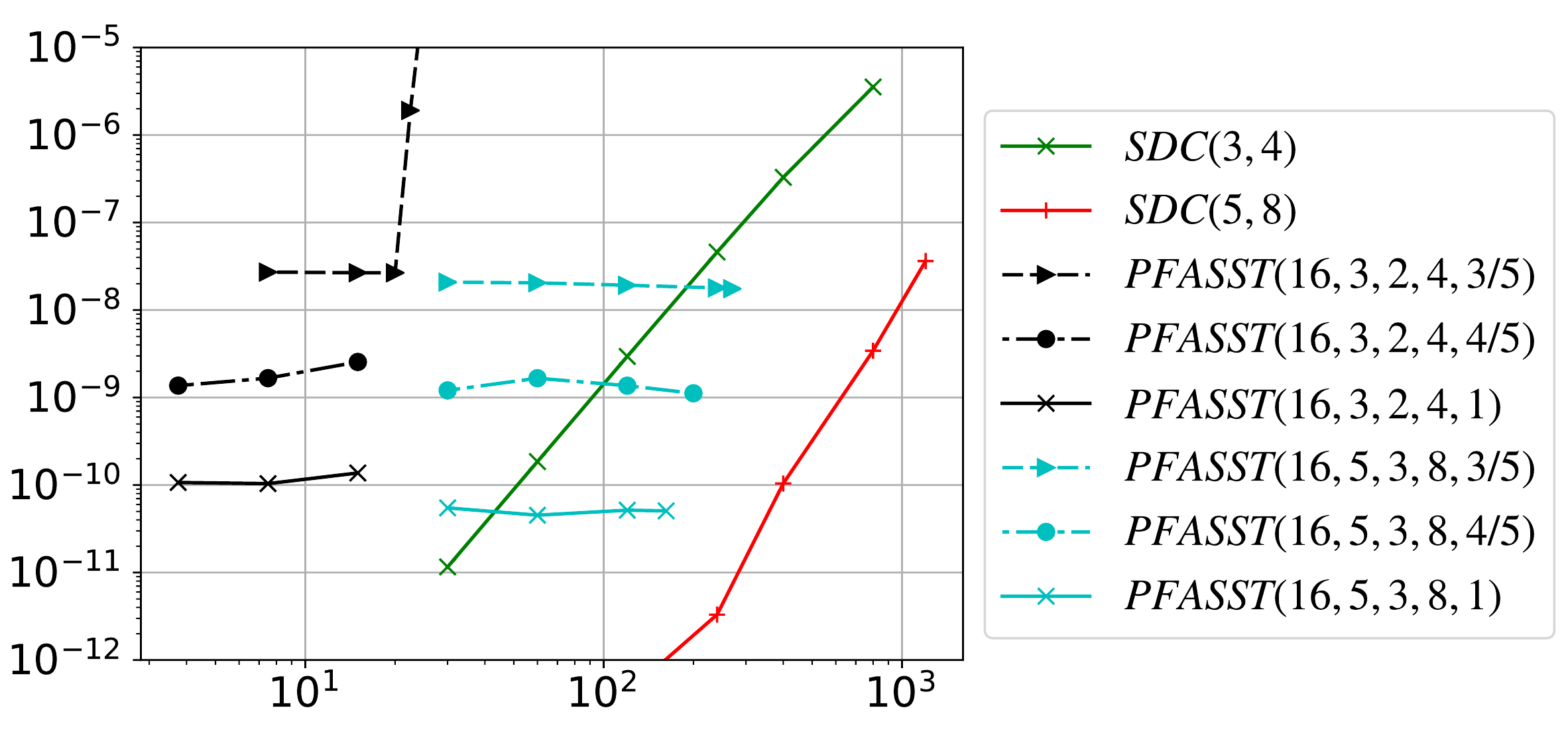}};
\node (ib_1) at (3.2,-0.05) {$\Delta t$ $(s)$};
\node (ib_1) at (3.1,4.45) {$R_{\textit{norm}} = 64$};
\end{tikzpicture}
\label{fig:galewsky_with_bump_test_case_accuracy_geopotential_32_visc_1e5}
}
\vspace{-0.45cm}
\caption{\label{fig:galewsky_with_bump_test_case_accuracy_geopotential_field_visc_1e5}
  Unstable barotropic wave: $L_{\infty}$-norm of the error in the geopotential field after
  $144 \, \text{hours}$ as a function of time step size. For both SDC and PFASST-SH, the norms
  are computed with $R_{\text{norm}} = 256$ in 
  \oldref{fig:galewsky_with_bump_test_case_accuracy_geopotential_256_visc_1e5},
  and $R_{\text{norm}} = 64$ in
  \oldref{fig:galewsky_with_bump_test_case_accuracy_geopotential_32_visc_1e5}.
  The fourth (respectively,
  fifth) parameter of the PFASST-SH schemes refers to the number of iterations (respectively,
  the coarsening ratio in space). 
}
\end{figure}

We now proceed to the analysis of the evolution of the error norm
as a function of time step size in  
Fig.~\oldref{fig:galewsky_with_bump_test_case_accuracy_geopotential_field_visc_1e5}.
We focus on PFASST($n_{\textit{ts}} = 16$, $M_f + 1 = 3$, $M_c + 1 = 2$,
$N_{\textit{PF}} = 4$, $\alpha$) and PFASST($n_{\textit{ts}} = 16$,
$M_F + 1 = 5$, $M_c + 1 = 3$, $N_{\textit{PF}} = 8$, $\alpha$) with
four and eight iterations, respectively.
We have chosen these values of $N_{\textit{NF}}$ to be
consistent with the parameters used in the other test cases.
For this example characterized by the presence of large spectral
coefficients associated with high-frequency modes, PFASST-SH is
accurate only if the spatial coarsening  ratio is close to one.
But, for a small value of $\alpha$, most of the spectrum is
represented on the coarse level, which imposes a more severe
stability restriction on the time step size. This is particularly
the case for PFASST($n_{\textit{ts}} = 16$, $M_f + 1 = 3$,
$M_c + 1 = 2$, $N_{\textit{PF}} = 4$, $\alpha$), as this scheme is
only stable for $\Delta t \leq 15 \, \text{s}$. This is almost
two orders of magnitude smaller than the largest time step sizes
achieved by the serial schemes.

With more SDC nodes on the fine
and coarse levels, PFASST($n_{\textit{ts}} = 16$, $M_F + 1 = 5$,
$M_c + 1 = 3$, $N_{\textit{PF}} = 8$, $\alpha$) can still achieve
stable time step sizes as large as $240 \, \text{s}$ for
$\alpha = 3/5$, which is only one order of magnitude smaller
than the largest time step size taken by SDC($M_f + 1 = 5$,
$N_S = 8$). However, we see that reducing the time step size
for a fixed number of eight iterations in PFASST($n_{\textit{ts}} = 16$,
$M_F + 1 = 5$, $M_c + 1 = 3$, $N_{\textit{PF}} = 8$, $\alpha$) does not
reduce the error. 
We show later in this section that for the Galewsky benchmark test,
PFASST-SH performs better with a larger number of iterations.

\begin{figure}[ht!]
\centering
\subfigure[]{
\begin{tikzpicture}
\node[anchor=south west,inner sep=0] at (0,0.045){\includegraphics[scale=0.385]{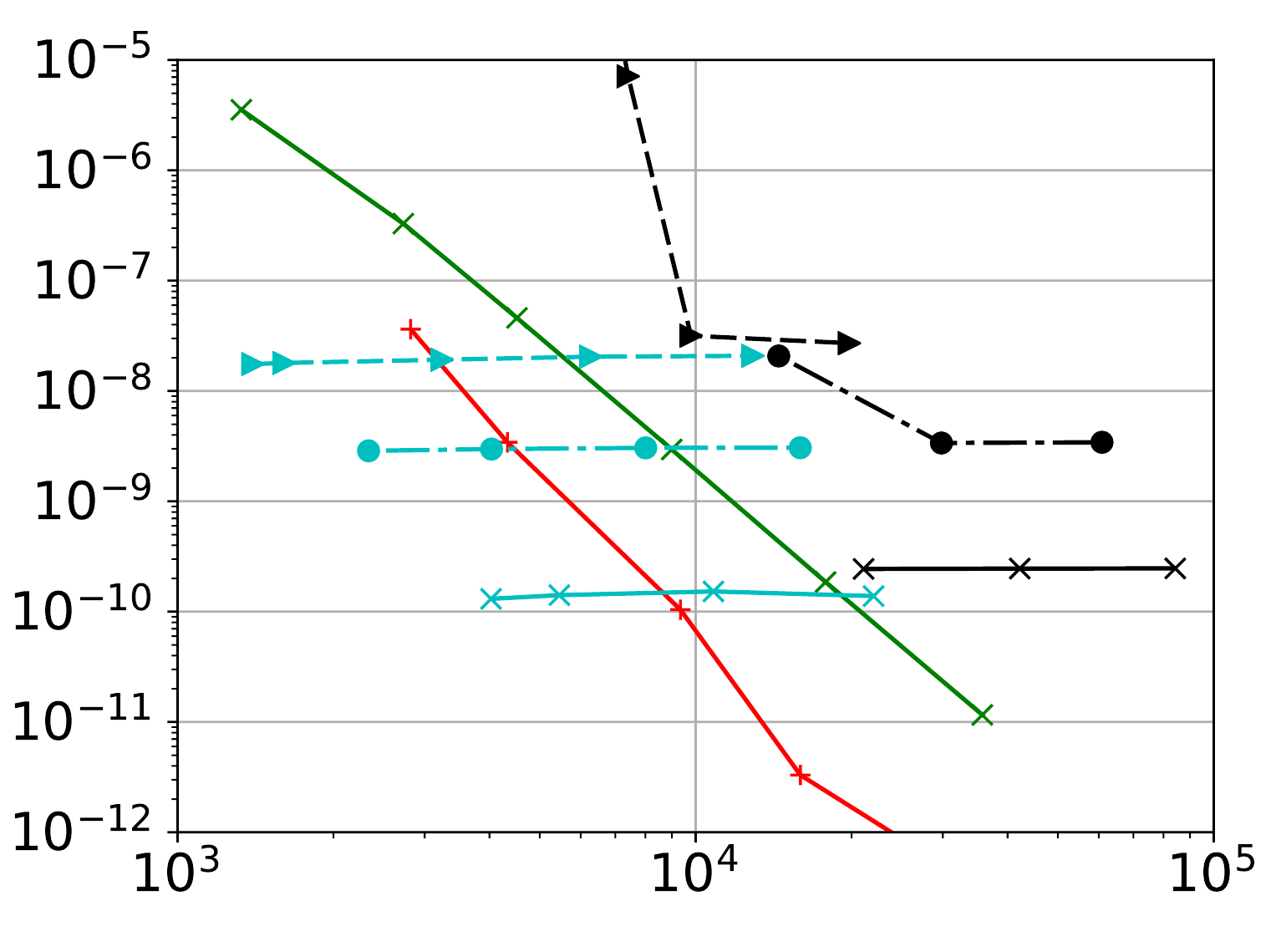}};
\node (ib_1) at (3.2,-0.05) {Wall-clock time $(s)$};
\node[rotate=90] (ib_1) at (-0.1,2.45) {Normalized $L_{\infty}$-error};
\node (ib_1) at (3.3,4.4) {$R_{\textit{norm}} = 256$};

\node (ib_1) at (1.4,2.95) {\scriptsize $\alpha = 3/5$};
\node (ib_1) at (1.65,2.53) {\scriptsize $\alpha = 4/5$};
\node (ib_1) at (2.45,1.875) {\scriptsize $\alpha = 1$};

\end{tikzpicture}
\label{fig:galewsky_with_bump_test_case_computational_cost_geopotential_256_visc_1e5}
}
\hspace{-0.5cm}
\subfigure[]{
\begin{tikzpicture}
\node[anchor=south west,inner sep=0] at (0,0.05){\includegraphics[scale=0.385]{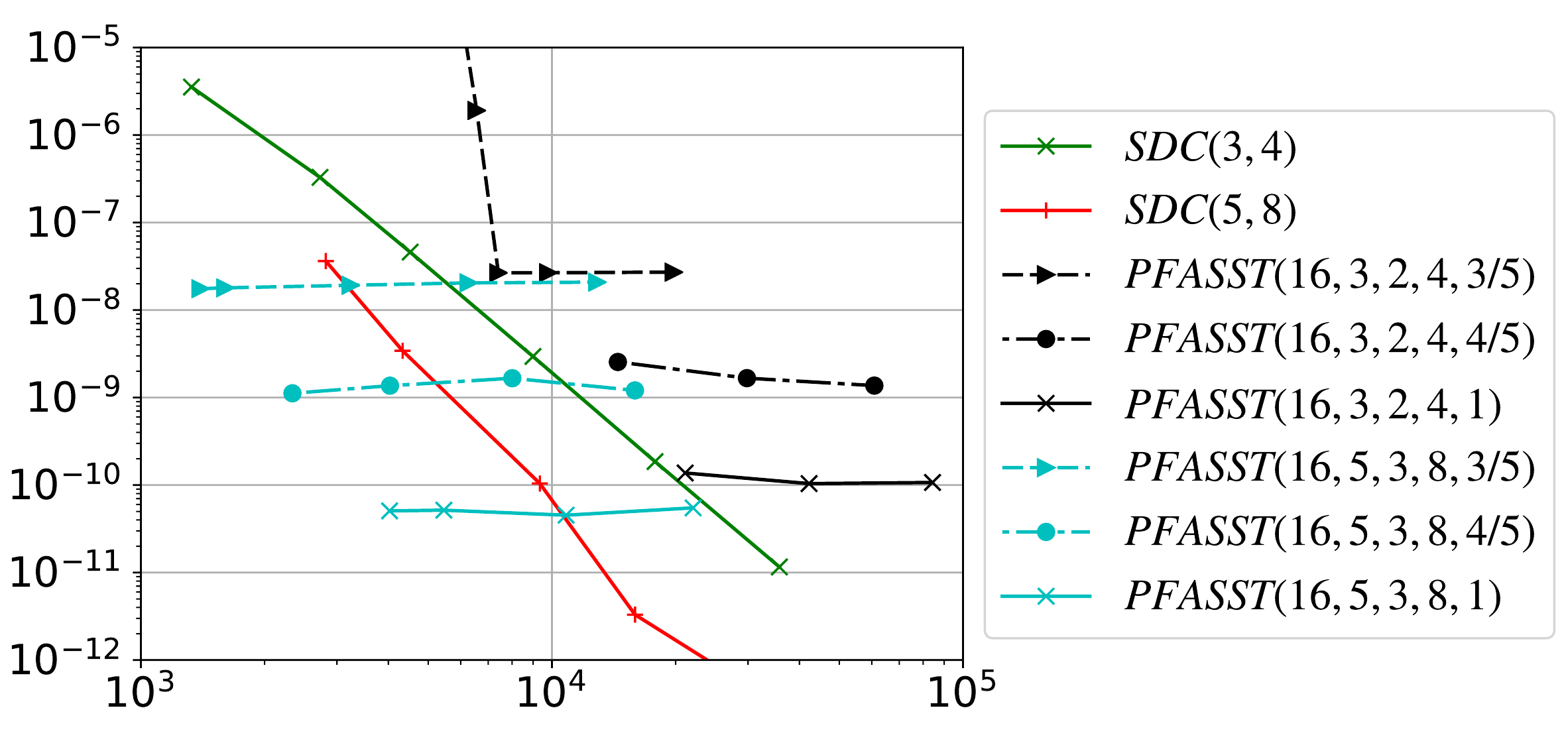}};
\node (ib_1) at (3.2,-0.05) {Wall-clock time $(s)$};
\node (ib_1) at (3.1,4.45) {$R_{\textit{norm}} = 64$};
\end{tikzpicture}
\label{fig:galewsky_with_bump_test_case_computational_cost_geopotential_32_visc_1e5}
}
\vspace{-0.45cm}
\caption{\label{fig:galewsky_with_bump_test_case_computational_cost_geopotential_field_visc_1e5} 
  Unstable barotropic wave: $L_{\infty}$-norm of the error in the geopotential field after
  $144 \, \text{hours}$ as a function of wall-clock time. For both SDC and PFASST-SH, the norms are
  computed $R_{\text{norm}} = 256$ in
\oldref{fig:galewsky_with_bump_test_case_computational_cost_geopotential_256_visc_1e5},
  and $R_{\text{norm}} = 64$ in
  \oldref{fig:galewsky_with_bump_test_case_computational_cost_geopotential_32_visc_1e5}.
  The fourth (respectively,
  fifth) parameter of the PFASST-SH schemes refers to the number of iterations (respectively,
  the coarsening ratio in space). For each scheme, the leftmost point corresponds to the
  largest stable time step. The figure shows that for their largest
  stable time step, the PFASST(16,5,3,8,$\cdot$) schemes reduce the computational cost compared to
  the serial SDC schemes, while the PFASST(16,3,2,4,$\cdot$) are consistently most expensive than their
  serial reference.
  }
\end{figure}

Figure~\oldref{fig:galewsky_with_bump_test_case_computational_cost_geopotential_field_visc_1e5}
illustrates the computational cost of PFASST-SH by showing the error norm
as a function of wall-clock time. Due to its limited stability, PFASST($n_{\textit{ts}} = 16$,
$M_f + 1 = 3$, $M_c + 1 = 2$, $N_{\textit{PF}} = 4$, $\alpha$) is
not competitive for this test case, as it is consistently slower than the serial
SDC schemes for all time step sizes. Conversely, for their largest stable time step,
the PFASST($n_{\textit{ts}} = 16$, $M_f + 1 = 5$, $M_c + 1 = 3$, $N_{\textit{PF}} = 8$,
$\alpha$) schemes are more efficient than both SDC($M_f + 1 = 5$, $N_S = 8$) and
SDC($M_f + 1 = 3$, $N_S = 4$). 
The speedup is slightly larger when the error norm
only accounts for the low-frequency spectral coefficients ($R_{\textit{norm}} = 64$).
This is because for $\alpha < 1$, the high-frequency features cannot be resolved
by PFASST-SH when the number of iterations is small compared to the number of time steps.

\begin{figure}[ht!]
\centering
\subfigure[]{
\begin{tikzpicture}
  \node[anchor=south west,inner sep=0] at (0,0.15){\includegraphics[scale=0.3825]{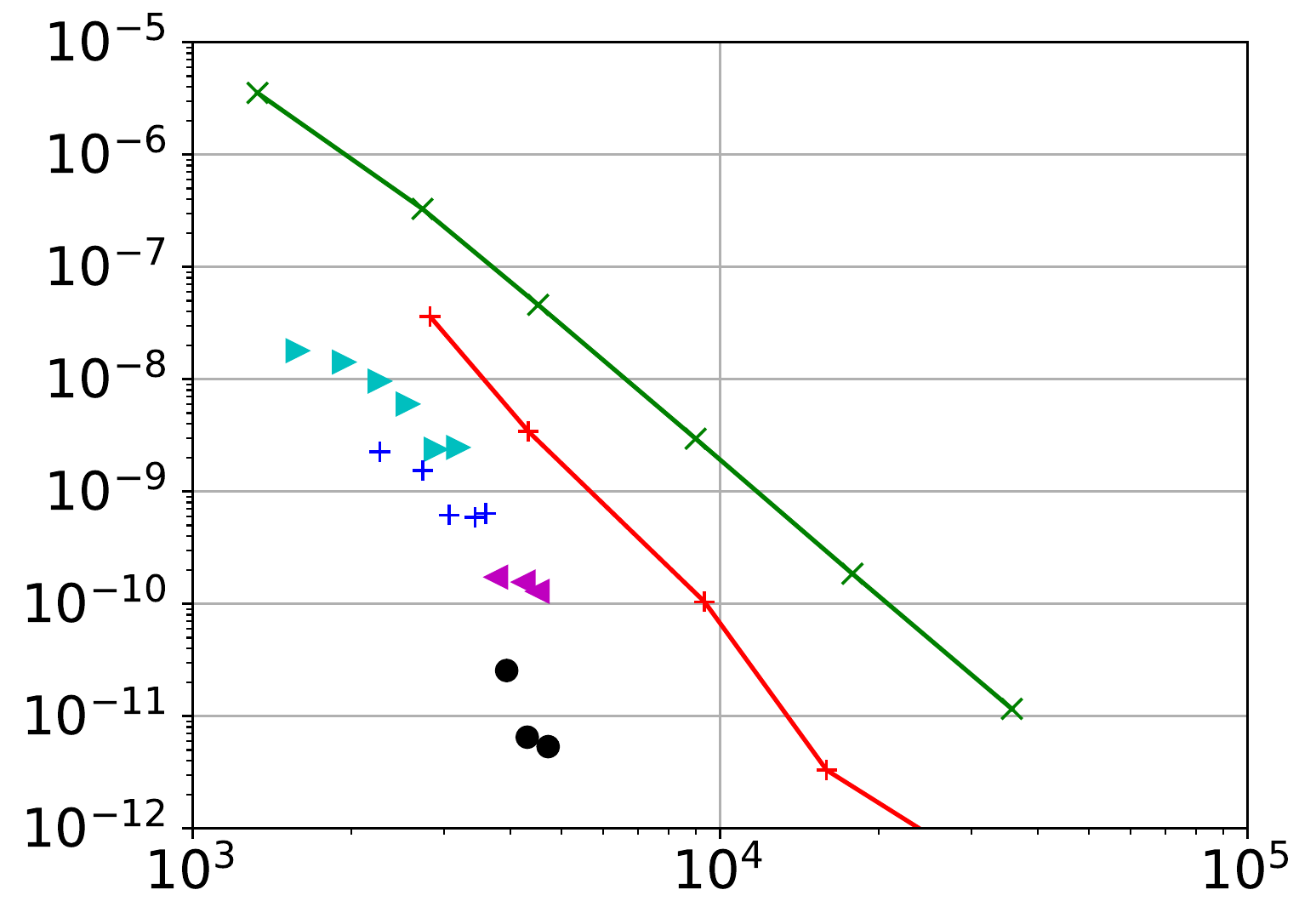}};
  \node (ib_1) at (2.715,0.9715) {\tiny \textcolor{black}{18}};
  \node (ib_1) at (2.1,1.3) {\tiny \textcolor{black}{14}};
  \node (ib_1) at (2.66,1.65) {\tiny \textcolor{magenta}{18}};
  \node (ib_1) at (2.05,1.764) {\tiny \textcolor{magenta}{14}};
  \node (ib_1) at (2.41,2.025) {\tiny \textcolor{blue}{18}};
  \node (ib_1) at (1.52,2.31) {\tiny \textcolor{blue}{10}};
  \node (ib_1) at (2.275,2.325) {\tiny \textcolor{cyan}{18}};
  \node (ib_1) at (1.17,2.77) {\tiny \textcolor{cyan}{8}};
\node (ib_1) at (3.2,-0.05) {Wall-clock time $(s)$};
\node[rotate=90] (ib_1) at (-0.05,2.35) {Normalized $L_{\infty}$-error};
\node (ib_1) at (3.3,4.4) {$R_{\textit{norm}} = 256$};
\end{tikzpicture}
\label{fig:galewsky_with_bump_test_case_computational_cost_geopotential_256_iters_all}
}
\hspace{-0.6cm}
\subfigure[]{
\begin{tikzpicture}
\node[anchor=south west,inner sep=0] at (0,0.05){\includegraphics[scale=0.3825]{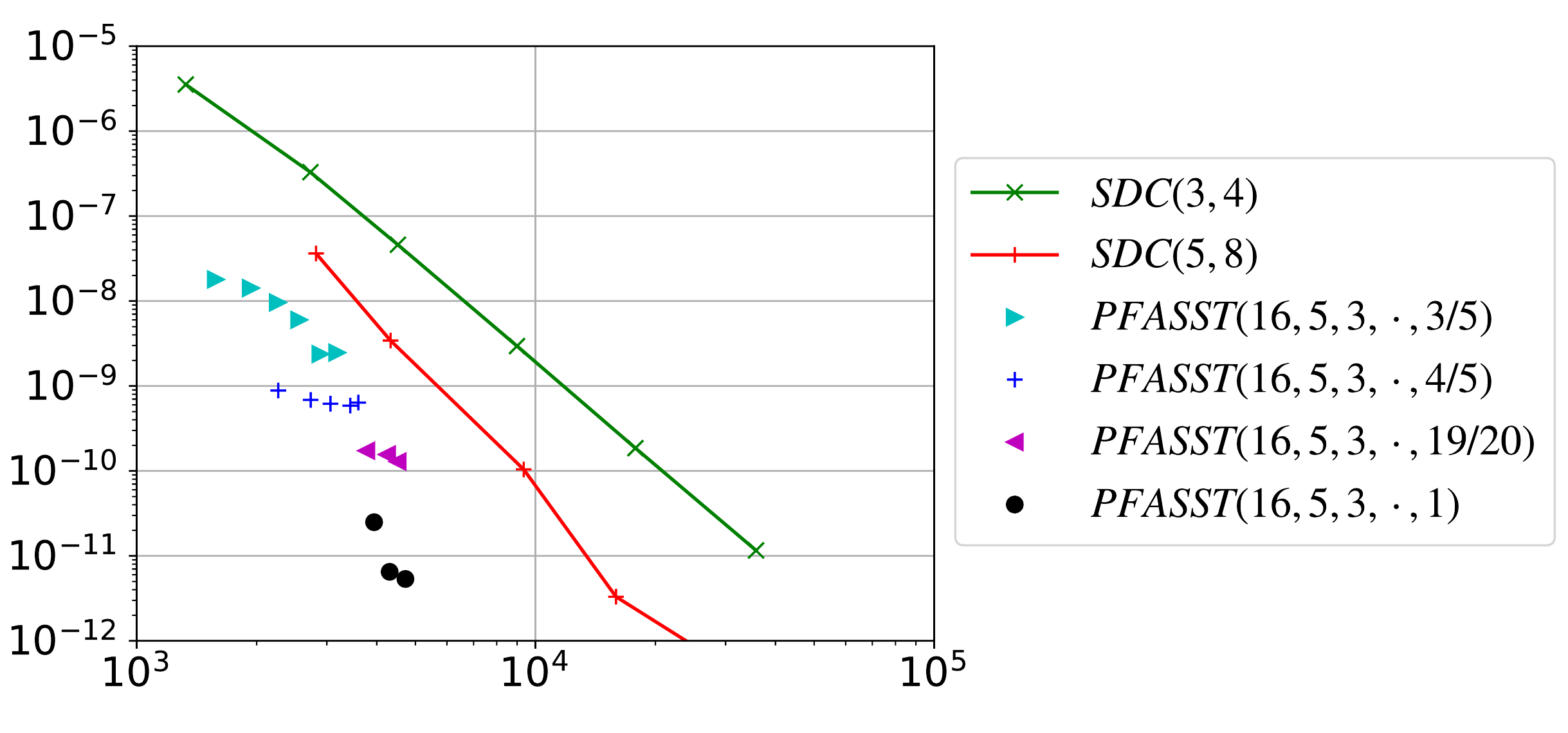}};
\node (ib_1) at (3.2,-0.05) {Wall-clock time $(s)$};
\node (ib_1) at (3.1,4.425) {$R_{\textit{norm}} = 64$};
\end{tikzpicture}
\label{fig:galewsky_with_bump_test_case_computational_cost_geopotential_32_iters_all}
}
\vspace{-0.45cm}
\caption{\label{fig:galewsky_with_bump_test_case_computational_cost_geopotential_field_iters_all} 
  Unstable barotropic wave: $L_{\infty}$-norms of the error in the geopotential field after
  $144 \, \text{hours}$ as a function of wall-clock time. The norms are
  computed $R_{\text{norm}} = 256$ in
  \oldref{fig:galewsky_with_bump_test_case_computational_cost_geopotential_256_iters_all}, 
  and $R_{\text{norm}} = 64$ in
  \oldref{fig:galewsky_with_bump_test_case_computational_cost_geopotential_32_iters_all}.
  We fix the time step size to $\Delta t = 240 \, \text{s}$ and increase the number
  of iterations. In
  \oldref{fig:galewsky_with_bump_test_case_computational_cost_geopotential_256_iters_all},
  for each spatial coarsening ratio, we annotate the data point corresponding to the
  stable PFASST-SH scheme with the smallest number of iterations, and the data point
  corresponding to PFASST(16,5,3,18,$\alpha$). We see that achieving stability
  for this time step size requires at least eight iterations for $\alpha = 3/5$,
  but at least 14 iterations for $\alpha = 1$.
}
\end{figure}

In Figs.~\oldref{fig:galewsky_with_bump_test_case_accuracy_geopotential_field_visc_1e5} and
\oldref{fig:galewsky_with_bump_test_case_computational_cost_geopotential_field_visc_1e5},  
the error norm obtained with the PFASST-SH schemes stagnates at a
relatively large magnitude. Here,
we study the sensitivity of PFASST-SH to the number
of iterations, $N_{\textit{PF}}$ using the same methodology as in
Fig.~\oldref{fig:gaussian_bump_test_case_computational_cost_geopotential_field}
to show that this limitation can be overcome by increasing
$N_{\textit{PF}}$. In
Fig.~\oldref{fig:galewsky_with_bump_test_case_computational_cost_geopotential_field_iters_all},
the PFASST-SH results are obtained with a fixed time step size
$\Delta t = 240 \, \text{s}$, a spatial coarsening ratio
$\alpha \in \{3/5, 4/5, 19/20, 1 \}$ and different values of the
number of iterations, $N_{\textit{PF}}$. We now focus on
PFASST($n_{\textit{ts}} = 16$, $M_f + 1 = 5$, $M_c + 1 = 3$,
$N_{\textit{PF}}$, $\alpha$) with $N_{PF} \geq 8$ as this family of
schemes yields consistent speedups compared to the serial SDC
schemes.

Figure~\oldref{fig:galewsky_with_bump_test_case_computational_cost_geopotential_field_iters_all}
shows that increasing the number of PFASST-SH
iterations reduces the error in the spectral coefficients for both
$R_{\textit{norm}} = 256$ and $R_{\textit{norm}} = 64$
until a threshold is reached at approximately
$N_{\textit{PF}} = 16$. We note that for this challenging test case,
the threshold is much higher than for the Gaussian test case
of Section~\oldref{subsection_propagation_of_a_gaussian_dome}.
This error reduction for $N_{\textit{PF}} < 16$
takes place because increasing $N_{\textit{PF}}$
allows to progressively propagate more information on
high-frequency modes to the last processor.

In terms of
computational cost, performing more PFASST-SH iterations only
slightly increases the wall-clock time. This is because in 
PFASST-SH, the cost of the parallel fine sweeps is
amortized over multiple iterations. Therefore, unlike Parareal,
PFASST-SH can still achieve a speedup when the number of iterations,
$N_{\textit{PF}}$, is close or equal to the number of processors,
$n_{\textit{ts}}$. In fact, we see in
Fig.~\oldref{fig:galewsky_with_bump_test_case_computational_cost_geopotential_field_iters_all}
that increasing $N_{\textit{PF}}$ leads to an increase in the speedup
compared to the serial schemes. 

\begin{table}[!ht]
\centering
\scalebox {0.9}{
             \begin{tabular}{lccccc}
           \\ \toprule 
           Integration         &    $L_{\infty}$-error       & Time step & Wall-clock & Observed & Observed  \\
           scheme              &    $R_{\textit{norm}} = 64$ & size (s)  & time (s)   & speedup $S^{\textit{obs}}_{S}$ & speedup $S^{\textit{obs}}_{\textit{ML}}$ \\\toprule
           SDC(5,8)             &    $3.3 \times 10^{-12}$            &        240           &  15,928      &  -   &   -     \\
           MLSDC(5,3,4,4/5)     &    $4.2 \times 10^{-12}$            &        60            &  46,127      & 0.3  &   -     \\
           PFASST(16,5,3,16,1)  &    $6.5 \times 10^{-12}$            &        240           &   4,311      & 3.7  &  10.7
           \\
           \bottomrule 
         \end{tabular}}
\caption{\label{tab:galewsky_with_bump_timings}
  Unstable barotropic wave: $L_{\infty}$-error for $R_{\textit{norm}} = 64$, time step size,
  and wall-clock time for the serial and parallel SDC-based schemes.
  $S^{\textit{obs}}_S$  denotes the observed speedup achieved with
  MLSDC-SH and PFASST-SH with respect to SDC, while
  $S^{\textit{obs}}_{\textit{ML}}$ denotes the speedup achieved
  with PFASST-SH with respect  to MLSDC-SH.
}
\end{table}

The speedups observed for this experiment are 
in Table~\oldref{tab:galewsky_with_bump_timings}. We highlight
that by running the PFASST-SH simulations for different numbers
of iterations and spatial coarsening ratios, we 
found a near-optimal set of parameters for this
challenging test case. In the best configuration, PFASST-SH uses
no coarsening in space ($\alpha = 1$) and performs 16 iterations
($N_{\textit{PF}} = n_{\textit{ts}}$). This yields an accurate
parallel-in-time scheme with a good stability. Specifically,
for an error norm of $3 \times 10^{-12}$,
PFASST($n_{\textit{ts}} = 16$, $M_f + 1 = 5$,
$M_c + 1 = 3$, $N_{\textit{PF}} = 16$, $\alpha = 1$) can 
take the same time step size of $240 \, \text{s}$ as
SDC($M_f + 1 = 5$, $N_S = 8$), leading to a speedup of 3.7.
For this error norm, MLSDC($M_f + 1 = 5$, $M_c + 1 = 3$,
$N_{\textit{ML}} = 4$, $\alpha = 4/5$) is more expensive
than the single-level SDC($M_f + 1 = 5$, $N_S = 8$), as shown
in \cite{hamon2018multi}. 

\begin{figure}[ht!]
\centering
\subfigure[]{
\begin{tikzpicture}
  \node[anchor=south west,inner sep=0] at (0,0.05){\includegraphics[scale=0.3825]{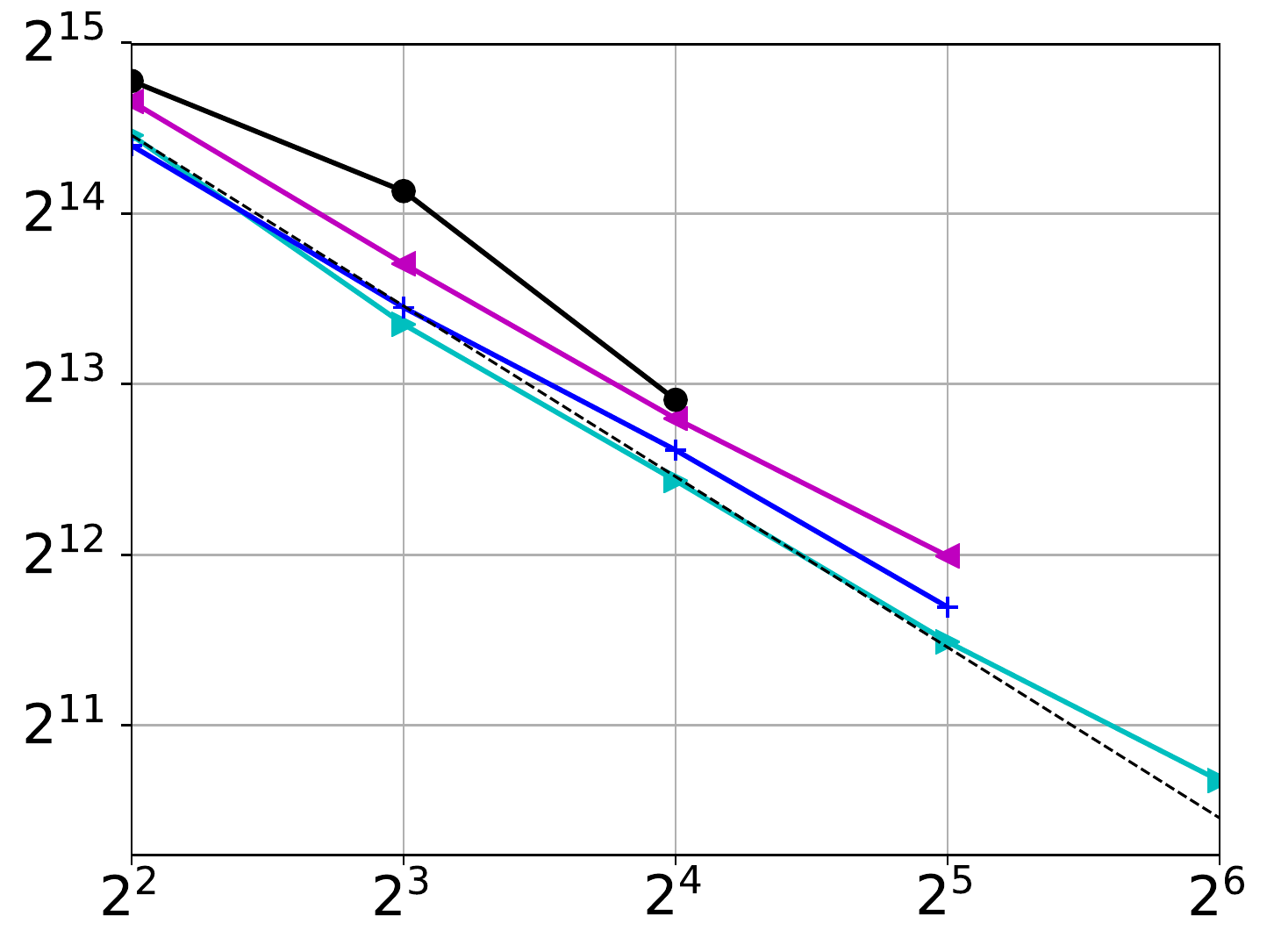}};
\node (ib_1) at (3.2,-0.05) {Number of processors};
\node[rotate=90] (ib_1) at (-0.1,2.35) {Wall-clock time};
\end{tikzpicture}
\label{fig:galewsky_with_bump_test_case_computational_cost_geopotential_total_orhun}
}
\hspace{-0.6cm}
\subfigure[]{
\begin{tikzpicture}
\node[anchor=south west,inner sep=0] at (0,0.05){\includegraphics[scale=0.3825]{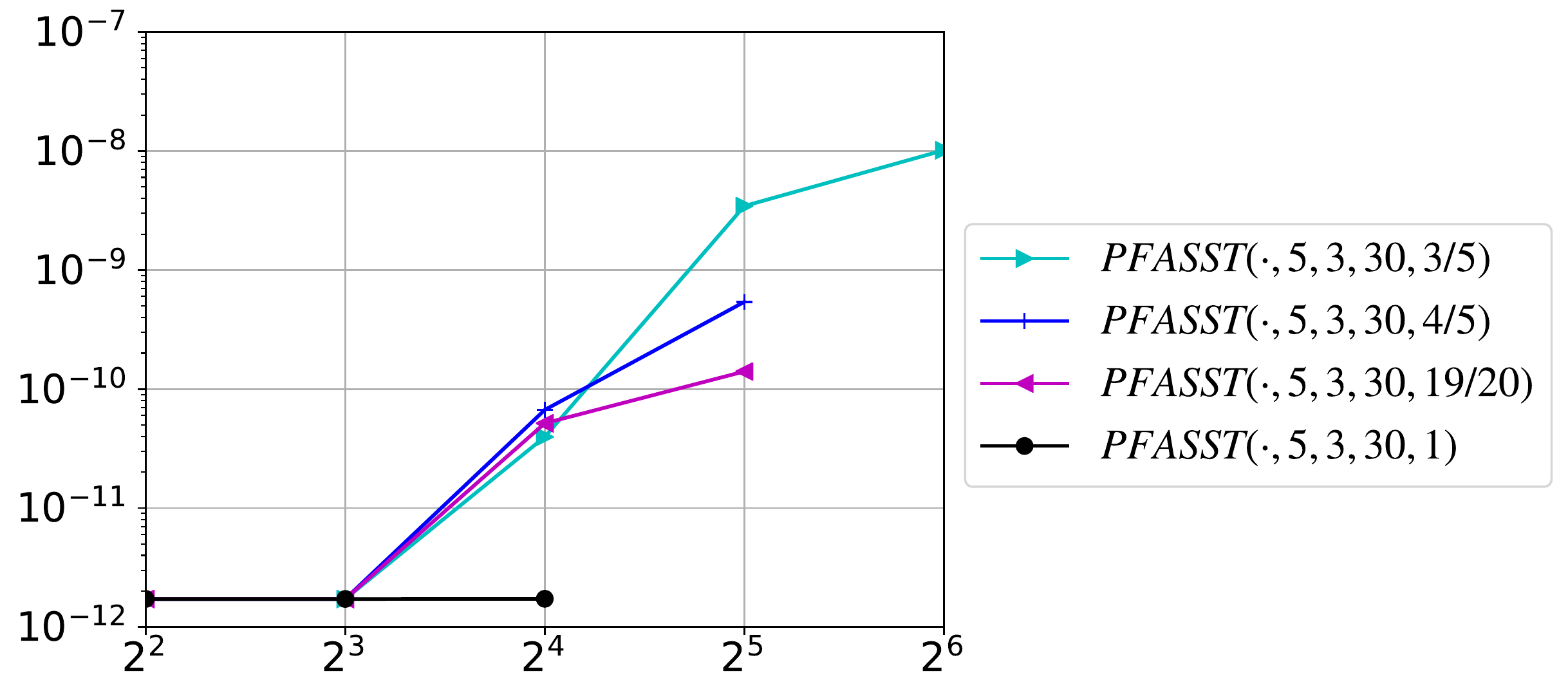}};
\node (ib_1) at (3.2,-0.05) {Number of processors};
\node[rotate=90] (ib_1) at (-0.1,2.35) {Normalized $L_{\infty}$-error};
\end{tikzpicture}
\label{fig:galewsky_with_bump_test_case_computational_cost_geopotential_total_marion}
}
\vspace{-0.45cm}
\caption{\label{fig:galewsky_with_bump_test_case_computational_cost_geopotential_total}
  Unstable barotropic wave: wall-clock time as a function of the number of processors
  in
  \oldref{fig:galewsky_with_bump_test_case_computational_cost_geopotential_total_orhun},
  and $L_{\infty}$-norm of the error as a function of the number of processors
  in
  \oldref{fig:galewsky_with_bump_test_case_computational_cost_geopotential_total_marion}.
  The norm of the error is computed with $R_{\textit{norm}} = 256$. As in the
  previous figures, these results are obtained for the simulation of
  $144 \, \text{hours}$  of propagation. For each value of the spatial coarsening coefficient,
  the right-most data point corresponds to the stable PFASST-SH simulation with the largest
  number of processors (i.e., using more processors with the same number of iterations
  would make PFASST-SH unstable).
}
\end{figure}  

Finally, we run a strong scalability study to
assess the evolution of the accuracy and efficiency of PFASST-SH
when the number of processors is increased from four to 64. The
results of this  study run
with the same resolution ($R_f = 256$) and
physical parameters as in the previous paragraphs are in
Fig.~\oldref{fig:galewsky_with_bump_test_case_computational_cost_geopotential_total}.
We simulate 144 hours with a time step of
$\Delta t = 225 \, s$.
Figure~\oldref{fig:galewsky_with_bump_test_case_computational_cost_geopotential_total_orhun}
shows that we can significantly reduce the wall-clock time by
increasing the number of processors, and therefore, the number
of time steps solved in parallel. Specifically, the wall-clock time of
PFASST($n_{\textit{ts}} = 16$, $M_f + 1 = 5$, $M_c + 1 = 3$, 
$N_{\textit{PF}} = 30$, $\alpha = 3/5$) is reduced by a factor
of 7.8 (respectively, 13.7) when the number of processors is
increased from four to 32 (respectively, 64). 
This is an encouraging result
although it does not prove that PFASST-SH can improve the scalability
of global spectral models -- this would require performing the same
study on a large-scale problem on which \textit{parallelism in space}
is saturated. But,
  Fig.~\oldref{fig:galewsky_with_bump_test_case_computational_cost_geopotential_total_marion} shows 
that increasing the number of processors causes a deterioration
of the error as PFASST-SH becomes unstable. This is a result
that will be improved in future work with the research directions
discussed in the conclusion.

\section{\label{section_conclusion}Conclusion and future work}

We have proposed an iterative, multilevel, parallel-in-time
integration method for the Shallow-Water Equations (SWE) on the
rotating sphere. Our approach combines the global Spherical
Harmonics (SH) transform with the Parallel Full Approximation
Scheme in Space and Time (PFASST) to construct a robust and
efficient numerical scheme for the nonlinear wave-propagation
problems arising from the SWE. The method computes multiple
time steps in parallel by applying a sequence of concurrent
fine corrections to iteratively improve an initial condition
propagated on the coarse level. A key feature of this algorithm
is the coarsening and interpolation procedure designed in
\cite{hamon2018multi} to accurately transfer the approximate
solution between levels without introducing spurious modes.

We have studied the properties of PFASST-SH using a suite
of standard nonlinear test cases designed for the development
of dynamical cores for weather and climate simulations.
These numerical tests illustrate the importance
of the spatial coarsening ratio, $\alpha$, for the accuracy
and computational cost of the scheme.
Aggressive spatial coarsening yields a PFASST-SH scheme
that only partially resolves the small-scale features of the solution, and
that is therefore inexpensive but also relatively inaccurate compared
to the serial single-level and multi-level SDC schemes.
Conversely, milder spatial coarsening leads
to a more costly, but significantly more accurate scheme
that captures the full spectrum of the solution. 
Our results show that PFASST-SH remains robust and
efficient for the challenging test case proposed in
\cite{galewsky2004initial} and characterized by the progressive
amplification of large high-frequency modes in the spectrum of
the solution. For this numerical experiment, we have demonstrated
that our parallel-in-time algorithm can both fully represent
the fine solution of the problem and achieve a 
speedup compared to the serial single-level and multi-level
schemes. 

Future work includes combining PFASST with Semi-Lagrangian (SL)
methods, which have 
exhibited excellent stability
properties in the context of a Parareal-based predictor-corrector
scheme \citep{schmitt_2018_sl_parareal}. This would aim at
increasing the largest stable time step size taken by 
PFASST to further reduce the  cost of the
parallel-in-time integration method. Another research
direction consists in
exploiting the potential of exponential integrators in PFASST,
as done before using a Parareal approach \citep{gander2013paraexp}
and a serial SDC approach \citep{buvoli2015class}.
Exponential integration could advantageously replace the
backward-Euler step in the PFASST iteration to reduce the dispersion
errors generated by the current scheme.
Finally, future work also includes
extending these results to the full atmospheric dynamics 
-- i.e., accounting for the horizontal \textit{and} vertical
components -- to demonstrate the applicability of PFASST
in dynamical cores and target  realistic, larger-scale
weather and climate simulations.

\section{Acknowledgements}
The work of Fran\c cois Hamon and Michael Minion was supported
by the U.S. Department of Energy, Office of Science, Office of
Advanced Scientific Computing Research, Applied Mathematics
program under contract number DE-AC02005CH11231. 
The simulations were performed at the National Energy Research
Scientific Computing Center (NERSC), a DOE Office of Science
User Facility supported by the U.S. DOE under Contract 
No. DE-AC02-05CH11231. 
Martin Schreiber received funding from NCAR for a research
stay in summer 2017 at the Mesa Labs which also led to
contributions to the present paper.

\bibliography{biblio}

\begin{thebibliography}{}

\bibitem[B{\'e}nard, 2019]{benard2019numerical}
B{\'e}nard, P. (2019).
\newblock Numerical investigation of {R}ossby waves for nonlinear shallow-water
  equations on the sphere.
\newblock {\em Quarterly Journal of the Royal Meteorological Society}.

\bibitem[Bolten et~al., 2017]{bolten2016multigrid}
Bolten, M., Moser, D., and Speck, R. (2017).
\newblock A multigrid perspective on the parallel full approximation scheme in
  space and time.
\newblock {\em Numerical Linear Algebra with Applications}, 24(6):e2110--n/a.

\bibitem[Bolten et~al., 2018]{bolten2017asymptotic}
Bolten, M., Moser, D., and Speck, R. (2018).
\newblock Asymptotic convergence of the parallel full approximation scheme in
  space and time for linear problems.
\newblock {\em Numerical Linear Algebra with Applications}, 25(6):e2208.

\bibitem[Bourke, 1972]{bourke1972efficient}
Bourke, W. (1972).
\newblock An efficient, one-level, primitive-equation spectral model.
\newblock {\em Monthly Weather Review}, 100(9):683--689.

\bibitem[Bourlioux et~al., 2003]{bourlioux2003high}
Bourlioux, A., Layton, A.~T., and Minion, M.~L. (2003).
\newblock High-order multi-implicit spectral deferred correction methods for
  problems of reactive flow.
\newblock {\em Journal of Computational Physics}, 189(2):651--675.

\bibitem[Brandt, 1977]{brandt1977multi}
Brandt, A. (1977).
\newblock Multi-level adaptive solutions to boundary-value problems.
\newblock {\em Mathematics of computation}, 31(138):333--390.

\bibitem[Burrage, 1997]{burrage1997parallel}
Burrage, K. (1997).
\newblock Parallel methods for {O}{D}{E}s.
\newblock {\em Advances in Computational Mathematics}, 7(1-2):1--31.

\bibitem[Butcher, 1997]{butcher1997order}
Butcher, J.~C. (1997).
\newblock Order and stability of parallel methods for stiff problems.
\newblock {\em Advances in Computational Mathematics}, 7(1):79--96.

\bibitem[Buvoli, 2015]{buvoli2015class}
Buvoli, T. (2015).
\newblock A class of exponential integrators based on spectral deferred
  correction.
\newblock {\em arXiv preprint arXiv:1504.05543}.

\bibitem[Christlieb et~al., 2009]{christlieb2009comments}
Christlieb, A., Ong, B., and Qiu, J.-M. (2009).
\newblock Comments on high-order integrators embedded within integral deferred
  correction methods.
\newblock {\em Communications in Applied Mathematics and Computational
  Science}, 4(1):27--56.

\bibitem[Dutt et~al., 2000]{dutt2000spectral}
Dutt, A., Greengard, L., and Rokhlin, V. (2000).
\newblock Spectral deferred correction methods for ordinary differential
  equations.
\newblock {\em BIT Numerical Mathematics}, 40(2):241--266.

\bibitem[Emmett and Minion, 2012]{emmett2012toward}
Emmett, M. and Minion, M.~L. (2012).
\newblock Toward an efficient parallel in time method for partial differential
  equations.
\newblock {\em Communications in Applied Mathematics and Computational
  Science}, 7(1):105--132.

\bibitem[Emmett and Minion, 2014]{emmett2014efficient}
Emmett, M. and Minion, M.~L. (2014).
\newblock Efficient implementation of a multi-level parallel in time algorithm.
\newblock In {\em Domain Decomposition Methods in Science and Engineering XXI},
  pages 359--366. Springer.

\bibitem[Evans et~al., 2010]{evans2010accuracy}
Evans, K.~J., Taylor, M.~A., and Drake, J.~B. (2010).
\newblock Accuracy analysis of a spectral element atmospheric model using a
  fully implicit solution framework.
\newblock {\em Monthly Weather Review}, 138(8):3333--3341.

\bibitem[Falgout et~al., 2014]{falgout2014parallel}
Falgout, R.~D., Friedhoff, S., Kolev, T.~V., MacLachlan, S.~P., and Schroder,
  J.~B. (2014).
\newblock Parallel time integration with multigrid.
\newblock {\em SIAM Journal on Scientific Computing}, 36(6):C635--C661.

\bibitem[Farhat and Chandesris, 2003]{farhat2003time}
Farhat, C. and Chandesris, M. (2003).
\newblock Time-decomposed parallel time-integrators: theory and feasibility
  studies for fluid, structure, and fluid--structure applications.
\newblock {\em International Journal for Numerical Methods in Engineering},
  58(9):1397--1434.

\bibitem[Galewsky et~al., 2004]{galewsky2004initial}
Galewsky, J., Scott, R.~K., and Polvani, L.~M. (2004).
\newblock An initial-value problem for testing numerical models of the global
  shallow-water equations.
\newblock {\em Tellus A}, 56(5):429--440.

\bibitem[Gander, 1999]{gander1999waveform}
Gander, M.~J. (1999).
\newblock A waveform relaxation algorithm with overlapping splitting for
  reaction diffusion equations.
\newblock {\em Numerical Linear Algebra with Applications}, 6(2):125--145.

\bibitem[Gander and G{\"u}ttel, 2013]{gander2013paraexp}
Gander, M.~J. and G{\"u}ttel, S. (2013).
\newblock {P}{A}{R}{A}{E}{X}{P}: A parallel integrator for linear initial-value
  problems.
\newblock {\em SIAM Journal on Scientific Computing}, 35(2):C123--C142.

\bibitem[Gardner et~al., 2018]{gardner2018implicit}
Gardner, D.~J., Guerra, J.~E., Hamon, F.~P., Reynolds, D.~R., Ullrich, P.~A.,
  and Woodward, C.~S. (2018).
\newblock Implicit--explicit ({I}{M}{E}{X}) {R}unge--{K}utta methods for
  non-hydrostatic atmospheric models.
\newblock {\em Geoscientific Model Development}, 11(4):1497.

\bibitem[Gelb and Gleeson, 2001]{gelb2001spectral}
Gelb, A. and Gleeson, J.~P. (2001).
\newblock Spectral viscosity for shallow water equations in spherical geometry.
\newblock {\em Monthly Weather Review}, 129(9):2346--2360.

\bibitem[Giraldo, 2005]{giraldo2005semi}
Giraldo, F.~X. (2005).
\newblock Semi-implicit time-integrators for a scalable spectral element
  atmospheric model.
\newblock {\em Quarterly Journal of the Royal Meteorological Society},
  131(610):2431--2454.

\bibitem[Giraldo et~al., 2013]{giraldo2013implicit}
Giraldo, F.~X., Kelly, J.~F., and Constantinescu, E.~M. (2013).
\newblock Implicit-explicit formulations of a three-dimensional nonhydrostatic
  unified model of the atmosphere ({N}{U}{M}{A}).
\newblock {\em SIAM Journal on Scientific Computing}, 35(5):B1162--B1194.

\bibitem[G{\"o}tschel and Minion, 2019]{Gotschel2019-ka}
G{\"o}tschel, S. and Minion, M.~L. (2019).
\newblock An efficient parallel-in-time method for optimization with parabolic
  {P}{D}{E}s.
\newblock {\em arXiv preprint arXiv:1901.06850}.

\bibitem[Hack and Jakob, 1992]{hack1992description}
Hack, J.~J. and Jakob, R. (1992).
\newblock {\em Description of a global shallow water model based on the
  spectral transform method}.
\newblock National Center for Atmospheric Research.

\bibitem[Hamon et~al., 2018]{hamon2018concurrent}
Hamon, F.~P., Day, M.~S., and Minion, M.~L. (2018).
\newblock Concurrent implicit spectral deferred correction scheme for
  low-{M}ach number combustion with detailed chemistry.
\newblock {\em Combustion Theory and Modelling}, 0(0):1--31.

\bibitem[Hamon et~al., 2019]{hamon2018multi}
Hamon, F.~P., Schreiber, M., and Minion, M.~L. (2019).
\newblock Multi-level spectral deferred corrections scheme for the shallow
  water equations on the rotating sphere.
\newblock {\em Journal of Computational Physics}, 376:435--454.

\bibitem[Haut and Wingate, 2014]{haut2014asymptotic}
Haut, T.~S. and Wingate, B.~A. (2014).
\newblock An asymptotic parallel-in-time method for highly oscillatory
  {P}{D}{E}s.
\newblock {\em SIAM Journal on Scientific Computing}, 36(2):A693--A713.

\bibitem[Jia et~al., 2013]{jia2013spectral}
Jia, J., Hill, J.~C., Evans, K.~J., Fann, G.~I., and Taylor, M.~A. (2013).
\newblock A spectral deferred correction method applied to the shallow water
  equations on a sphere.
\newblock {\em Monthly Weather Review}, 141(10):3435--3449.

\bibitem[Kanamitsu et~al., 1983]{kanamitsu1983description}
Kanamitsu, M., Tada, K., Kudo, T., Sato, N., and Isa, S. (1983).
\newblock Description of the {J}{M}{A} operational spectral model.
\newblock {\em Journal of the Meteorological Society of Japan. Ser. II},
  61(6):812--828.

\bibitem[Layton and Minion, 2004]{layton2004conservative}
Layton, A.~T. and Minion, M.~L. (2004).
\newblock Conservative multi-implicit spectral deferred correction methods for
  reacting gas dynamics.
\newblock {\em Journal of Computational Physics}, 194(2):697--715.

\bibitem[Lions et~al., 2001]{lions2001resolution}
Lions, J.-L., Maday, Y., and Turinici, G. (2001).
\newblock R{\'e}solution d'{E}{D}{P} par un sch{\'e}ma en temps parar{\'e}el.
\newblock {\em Comptes Rendus de l'Acad{\'e}mie des Sciences-Series
  I-Mathematics}, 332(7):661--668.

\bibitem[Lock et~al., 2014]{lock2014numerical}
Lock, S.-J., Wood, N., and Weller, H. (2014).
\newblock Numerical analyses of {R}unge--{K}utta implicit--explicit schemes for
  horizontally explicit, vertically implicit solutions of atmospheric models.
\newblock {\em Quarterly Journal of the Royal Meteorological Society},
  140(682):1654--1669.

\bibitem[Lott et~al., 2015]{lott2015algorithmically}
Lott, P.~A., Woodward, C.~S., and Evans, K.~J. (2015).
\newblock Algorithmically scalable block preconditioner for fully implicit
  shallow-water equations in {C}{A}{M}-{S}{E}.
\newblock {\em Computational Geosciences}, 19(1):49--61.

\bibitem[Minion, 2003]{minion2003semi}
Minion, M.~L. (2003).
\newblock Semi-implicit spectral deferred correction methods for ordinary
  differential equations.
\newblock {\em Communications in Mathematical Sciences}, 1(3):471--500.

\bibitem[Minion, 2011]{minion2011hybrid}
Minion, M.~L. (2011).
\newblock A hybrid parareal spectral deferred corrections method.
\newblock {\em Communications in Applied Mathematics and Computational
  Science}, 5(2):265--301.

\bibitem[Rivier et~al., 2002]{rivier2002efficient}
Rivier, L., Loft, R., and Polvani, L.~M. (2002).
\newblock An efficient spectral dynamical core for distributed memory
  computers.
\newblock {\em Monthly Weather Review}, 130(5):1384--1396.

\bibitem[Robert et~al., 1972]{robert1972implicit}
Robert, A., Henderson, J., and Turnbull, C. (1972).
\newblock An implicit time integration scheme for baroclinic models of the
  atmosphere.
\newblock {\em Monthly Weather Review}, 100(5):329--335.

\bibitem[Ruprecht, 2018]{ruprecht2018wave}
Ruprecht, D. (2018).
\newblock Wave propagation characteristics of parareal.
\newblock {\em Computing and Visualization in Science}, 19(1-2):1--17.

\bibitem[Ruprecht and Speck, 2016]{ruprecht2016spectral}
Ruprecht, D. and Speck, R. (2016).
\newblock Spectral deferred corrections with fast-wave slow-wave splitting.
\newblock {\em SIAM Journal on Scientific Computing}, 38(4):A2535--A2557.

\bibitem[Schaeffer, 2013]{schaeffer2013efficient}
Schaeffer, N. (2013).
\newblock Efficient spherical harmonic transforms aimed at pseudospectral
  numerical simulations.
\newblock {\em Geochemistry, Geophysics, Geosystems}, 14(3):751--758.

\bibitem[Schmitt et~al., 2018]{schmitt_2018_sl_parareal}
Schmitt, A., Schreiber, M., Peixoto, P., and Sch{\"a}fer, N. (2018).
\newblock A numerical study of a semi-{L}agrangian {P}arareal method applied to
  the viscous {B}urgers equation.
\newblock {\em Computing and Visualization in Science}, 19(1-2):45--57.

\bibitem[Schreiber and Loft, 2018]{schreiber2018sph}
Schreiber, M. and Loft, R. (2018).
\newblock A parallel time integrator for solving the linearized shallow water
  equations on the rotating sphere.
\newblock {\em Numerical Linear Algebra with Applications}, page e2220.

\bibitem[Schreiber et~al., 2017]{schreiber2017beyond}
Schreiber, M., Peixoto, P.~S., Haut, T., and Wingate, B. (2017).
\newblock Beyond spatial scalability limitations with a massively parallel
  method for linear oscillatory problems.
\newblock {\em The International Journal of High Performance Computing
  Applications}, pages 1--21.

\bibitem[Schreiber et~al., 2019]{schreiber_2019_expnonlinearswe_sphere}
Schreiber, M., Schaeffer, N., and Loft, R. (2019).
\newblock Exponential integrators with parallel-in-time rational approximations
  for shallow-water equations on the rotating sphere.
\newblock {\em Parallel Computing}.

\bibitem[Smolarkiewicz et~al., 2014]{smolarkiewicz2014consistent}
Smolarkiewicz, P.~K., K{\"u}hnlein, C., and Wedi, N.~P. (2014).
\newblock A consistent framework for discrete integrations of soundproof and
  compressible {P}{D}{E}s of atmospheric dynamics.
\newblock {\em Journal of Computational Physics}, 263:185--205.

\bibitem[Speck et~al., 2014]{speck2014space}
Speck, R., Ruprecht, D., Emmett, M., Bolten, M., and Krause, R. (2014).
\newblock A space-time parallel solver for the three-dimensional heat equation.
\newblock {\em Parallel Computing: Accelerating Computational Science and
  Engineering (CSE)}, 25:263--272.

\bibitem[Speck et~al., 2012]{speck2012massively}
Speck, R., Ruprecht, D., Krause, R., Emmett, M., Minion, M.~L., Winkel, M., and
  Gibbon, P. (2012).
\newblock A massively space-time parallel {N}-body solver.
\newblock In {\em Proceedings of the International Conference on High
  Performance Computing, Networking, Storage and Analysis}, page~92. IEEE
  Computer Society Press.

\bibitem[Swarztrauber, 2004]{swarztrauber2004shallow}
Swarztrauber, P.~N. (2004).
\newblock Shallow water flow on the sphere.
\newblock {\em Monthly Weather Review}, 132(12):3010--3018.

\bibitem[Temperton, 1991]{temperton1991scalar}
Temperton, C. (1991).
\newblock On scalar and vector transform methods for global spectral models.
\newblock {\em Monthly Weather Review}, 119(5):1303--1307.

\bibitem[Thuburn and Li, 2000]{thuburn2000numerical}
Thuburn, J. and Li, Y. (2000).
\newblock Numerical simulations of rossby--haurwitz waves.
\newblock {\em Tellus A}, 52(2):181--189.

\bibitem[Ullrich and Jablonowski, 2012]{ullrich2012operator}
Ullrich, P. and Jablonowski, C. (2012).
\newblock Operator-split {R}unge--{K}utta--{R}osenbrock methods for
  nonhydrostatic atmospheric models.
\newblock {\em Monthly Weather Review}, 140(4):1257--1284.

\bibitem[Wedi et~al., 2013]{wedi2013fast}
Wedi, N.~P., Hamrud, M., and Mozdzynski, G. (2013).
\newblock A fast spherical harmonics transform for global nwp and climate
  models.
\newblock {\em Monthly Weather Review}, 141(10):3450--3461.

\bibitem[Weiser, 2015]{weiser2015faster}
Weiser, M. (2015).
\newblock Faster {S}{D}{C} convergence on non-equidistant grids by {D}{I}{R}{K}
  sweeps.
\newblock {\em BIT Numerical Mathematics}, 55(4):1219--1241.

\bibitem[Weller et~al., 2013]{weller2013runge}
Weller, H., Lock, S.-J., and Wood, N. (2013).
\newblock Runge--{K}utta {I}{M}{E}{X} schemes for the horizontally
  explicit/vertically implicit ({H}{E}{V}{I}) solution of wave equations.
\newblock {\em Journal of Computational Physics}, 252:365--381.

\bibitem[Williamson et~al., 1992]{williamson1992standard}
Williamson, D.~L., Drake, J.~B., Hack, J.~J., Jakob, R., and Swarztrauber,
  P.~N. (1992).
\newblock A standard test set for numerical approximations to the shallow water
  equations in spherical geometry.
\newblock {\em Journal of Computational Physics}, 102(1):211--224.

\end{thebibliography}

\end{document}